\magnification=1200

\tolerance=500

\font \twelvebf=cmbx12

\def \om {{\omega_{\cC/T}}}

\hfuzz=3pt

\def \a {\alpha}

\def \avd {anneau de valuation discr\`ete}

\def \Coker {{\rm Coker\,}}

\def\dem{\noindent {\bf D\'emonstration.}\enspace \nobreak }

\def \de {\delta}

\def \e {\epsilon}

\def \equi {\Longleftrightarrow}

\def \expl#1 {\medbreak \noindent {\bf Exemple #1.}\enspace }
\def \expls#1 {\medbreak \noindent {\bf Exemples #1.}\enspace }

\def \cExt {{{\cal E}\it xt\,}}
\def \Ext {{\rm Ext\,}}

\def \f {\varphi }
\def \fl {\rightarrow }
\def \lf {\leftarrow }
\def \Fl#1{\smash{\mathop{\longrightarrow}\limits^{#1}}}

\def \g {\gamma}
\def \G {\Gamma}
\def \Gr{{ \rm Gr}}

\def \Hom { {\rm Hom\,}}
\def \Homgr {{\rm Homgr\,}}

\def \id #1 {<\!\! #1 \!\! >}
\def \Im {{ \rm Im\,}}

\def \impl {\Rightarrow}

\def \Ker { {\rm Ker\,}}

\def \la {\lambda}

\def \lign {\hfil \break } 

\def \Vf#1{\Big\downarrow \rlap{$\vcenter{\hbox{$\scriptstyle#1$}}$}}

\def \Mod {{\rm Mod}\,}

\def \ov {\overline}

\def \page#1{\leaders\hbox to 5 mm{\hfil.\hfil}\hfill\rlap{\hbox to
10mm{\hfill#1}}\par}

\def \ps {{\it psi}}

\def \qis {{\it qis}}

\def \rema#1 {\medbreak \noindent {\bf Remarque #1.}\enspace }
\def \remas#1 {\medbreak \noindent {\bf Remarques #1.}\enspace }

\def \s {\sigma}

\def \seFl#1{\smash{\mathop{\searrow}\nolimits\!\!\!^{#1}}}
\def \sh {{^\sharp}}
\def \Spec {{\rm Spec\,}}

\def \T {\otimes }
\def \tarte#1 {\medbreak \noindent {\it #1.}\medbreak}
\def \ti {\times}
\def \titre#1{\medbreak \noindent {\bf #1.}\medbreak}

\outer \def \th #1. #2\par{ \medbreak 
\noindent {\bf#1. \enspace} {\sl#2 }\par
\ifdim \lastskip< \medskipamount \removelastskip \penalty55 \medskip \fi}
\def \Tor {{\rm Tor \,}}
\def \ul {\underline}

 \def \vide {\emptyset}
\def  \Vf#1{\Big\downarrow \rlap{$\vcenter{\hbox{$\scriptstyle#1$}}$}}
\def \vf {\downarrow}

\def \wi {\widetilde}
\def \wh {\widehat}

\def \Lpont {{L^{\vphantom L}}_{\hbox {\bf .}}}
\def \Fpont {{F^{\vphantom F}}_{\hbox {\bf .}}}
\def \Mpont {{M^{\vphantom M}}_{\hbox {\bf .}}}

\def \Ppont {{P^{\vphantom P}}_{\hbox {\bf .}}}

\def \Xpont {{X^{\vphantom T}}_{\hbox {\bf .}}}

\def \Lpontpr {{L^{\vphantom L}}_{\hbox {\bf .} >r}}
\def \Lpontmr {{L^{\vphantom L}}_{\hbox {\bf .} \leq r}}
\def \Lpontps {{L^{\vphantom L}}_{\hbox {\bf .} >s}}

\def \Qpont {{Q^{\vphantom T}}_{\hbox {\bf .}}}
\def \upont {{u^{\vphantom u}}_{\hbox {\bf .}}}
\def \Kpont {{K^{\vphantom T}}_{\hbox {\bf .}}}
\def \tas {{\hbox {\bf .}}}

\def \cC {{\cal C}}
\def \cD {{\cal D}}
\def \cE {{\cal E}}
\def \cL {{\cal L}}
\def \cP {{\cal P}}
\def \cF {{\cal F}}
\def \cG {{\cal G}}

\def \cI {{\cal I}}
\def \cJ {{\cal J}}
\def \cO {{\cal O}}
\def \cK {{\cal K}}
\def \cM {{\cal M}}
\def \cN {{\cal N}}

\def \cS {{\cal S}}
\def \cT {{\cal T}}

\def \bP {{\bf P}}

\def \bR {{\bf R}}

\def \bZ {{\bf Z}}
\def \bN {{\bf N}}

\def \sN{{\sl N}}

\def \om {{\omega_{\cC/T}}}
\def \Om {{\Omega}}

\def \Id {{\rm Id}}

\def \sE {{\sl E}}
\def \sN {{\sl N}}

\null
\vskip 6 cm
\centerline {\twelvebf TRIADES ET}
\vskip 1 cm
\centerline {\twelvebf FAMILLES DE COURBES GAUCHES}

\vskip 4 cm
\centerline {\twelvebf Robin HARTSHORNE}
\vskip 1 cm
\centerline {\twelvebf Mireille MARTIN-DESCHAMPS}
\vskip 1 cm
\centerline {\twelvebf Daniel PERRIN}

\vfill\eject

\pageno -1

{\bf \centerline {TABLE DES MATI\`ERES}}

\vskip  1.5 cm

\noindent  Introduction  \page  {1}

\vskip 0.5 cm

\noindent 0. Notations et pr\'eliminaires  \page {8}

\indent a) Modules gradu\'es \page {8}

\indent b) Modules duaux \page {9}

\indent c) Foncteurs \page{9}

\indent d) Dualit\'e de Grothendieck et graduation \page{10}

\indent e) Courbes \page{12}

\vskip 0.5 cm
\noindent 1. La d\'efinition des triades \page{12}

\indent a) Complexes et foncteurs associ\'es \page{12}

\indent b) Pseudo-isomorphismes \page{13}

\indent c) Triades \page{14}

\indent d) Triades majeures \page{15}

\indent e) \'Etude des pseudo-isomorphismes de triades  \page{15}

\indent f) Un ``lemme de Verdier'' \page{16}

\indent g) Triades duales 
\page{19}

\indent h) Triades minimales, triades \'el\'ementaires \page{20}
 
\indent i) Exemples \page{22}

\indent j) Annexe : g\'en\'eralisation  \page{23}
\vskip 0.5 cm
\noindent 2. Faisceaux triadiques \page{24}

\indent a) Faisceaux triadiques \page{25}

\indent b) Triades et faisceaux triadiques \page{26}

\indent c) Pseudo-isomorphismes \page{27}

\indent d) Le lemme de Verdier inverse pour les faisceaux triadiques \page{29}

\indent e) Le cas d'un \avd\ \page{31}

\vskip 0.5 cm

\noindent 3. Courbes et triades : les th\'eor\`emes de Rao \page{31}

\indent a) R\'esolutions de type \sE\ et \sN\ triadiques  \page{31}

\indent b) Triade associ\'ee \`a une famille de courbes \page{34}

\indent c) Le th\'eor\`eme de Rao pour les triades : fibres de $\Psi_A$ \page{36}

\indent d) Le th\'eor\`eme de Rao pour les triades : image de $\Psi_A$ \page{36}

\indent e) Le cas de la liaison impaire \page{37}

\vskip 0.5 cm

\noindent 4. Courbes, faisceaux et triades : dictionnaire \page{38}

\indent  a) Propri\'et\'es de triades et des foncteurs triadiques \page{38}

\indent b) Le dictionnaire courbes-faisceaux-triades \page{40}

\vskip 0.5 cm

\noindent 5. Construction de triades, application \`a la construction de
familles de courbes \page{42}

\indent  a)  Construction d'une triade \`a partir de son c\oe ur, de son conoyau et d'une
 extension de longueur $2$ de ces modules  \page{43}

\indent b) Le cas d'un \avd \page{44}

\indent c)  Sous-quotient associ\'e \`a une triade \page{46}

\indent d) Construction de triades \`a partir de d\'eformations de sous-quotients : analyse des conditions
n\'eces\-saires 
\page{47}

\indent e) Construction de triades \`a partir de sous-quotients : la construction
triviale
\page{49}

\indent f) Calcul des familles de courbes obtenues \`a partir des triades triviales
\page{51}

\indent g) Construction de triades modulaires : variation autour du module $H$
\page{53}

\indent h) Constructions de triades \`a partir d'un sous-quotient : bis
\page{54}

\indent i) L'exemple des courbes de degr\'e $4$ et genre $0$ : construction de la triade
\page{56}

\indent j) L'exemple des $(4,0)$ : construction de la famille de courbes
\page{57}

\vskip 0.5 cm

\noindent Bibliographie \page {59}
\vfill\eject

\pageno=1

\centerline{\bf Introduction}

\vskip 1 cm

Cet article  concerne la classification des courbes gauches,
c'est-\`a-dire  l'\'etude du sch\'ema de Hilbert $H_{d,g}$
des courbes  (localement  Cohen-Macaulay et \'equidimensionnelles), de degr\'e $d$
et genre arithm\'etique $g$, de $\bP^3$ (espace projectif de dimension $3$  sur un
corps
$k$ alg\'ebrique\-ment clos). Il s'inscrit dans le
programme
\'enonc\'e voil\`a d\'ej\`a plusieurs ann\'ees par deux d'entre  nous dans
[MDP1]. La voie d'acc\`es au sch\'ema de Hilbert que nous privil\'egions est
l'utilisation du module de Rao. Si on pose $R = k[X,Y,Z,T]$, le module de Rao d'une
courbe
$C$  est le $R$-module gradu\'e de longueur finie :
$$M_C =
\bigoplus_{n \in
\bZ} H^1 \cJ_C (n).$$ 

L'utilisation de ce module a conduit \`a stratifier le sch\'ema de Hilbert
$H_{d,g}$ par les sous-sch\'emas \`a cohomologie constante $H_{\g,\rho}$ sur
lesquels les dimensions $h^i \cJ_C (n)$ des espaces de cohomologie sont constantes et
d\'etermin\'ees par les fonctions $\g$ (li\'ee \`a la postulation $h^0 \cJ_C
(n)$)  et
$\rho(n) = h^1 \cJ_C(n)$, cf. [MDP1] VI, (la sp\'ecialit\'e $h^2 \cJ_C (n) = h^1
\cO_C (n)$ est d\'etermin\'ee par les deux autres). Sur le sch\'ema
$H_{\g,\rho}$ on a un  morphisme
$\Phi :H_{\g,\rho}
\fl E_{\rho}$ (\`a valeurs dans le foncteur des structures de modules de dimensions
indiqu\'ees par la fonction $\rho$), qui
\`a une courbe
$C$ associe son module de Rao $M_C$. L'\'etude de $H_{d,g}$ est alors
d\'ecompos\'ee en trois \'etapes : l'\'etape du bas qui consiste \`a
\'etudier $E_{\rho}$ et a \'et\'e abord\'ee notamment dans [MDP2] (cf. aussi [G]),
l'\'etape interm\'ediaire qui revient \`a \'etudier $\Phi$ et qui a \'et\'e
r\'esolue dans [MDP1] VII : $\Phi$ est lisse, irr\'eductible et on conna\^\i t la
dimension de ses fibres, et enfin l'\'etape du haut qui 
\'etudie le recollement des $H_{\g,\rho}$ pour obtenir des informations globales  sur
$H_{d,g}$.

C'est de cette derni\`ere \'etape dont il est question ici. L'objectif est
 d'analyser les sp\'ecialisations dans le sch\'ema de Hilbert
$H_{d,g}$ ; pr\'ecis\'ement, si $H_0$ (resp. $H$) est une composante irr\'eductible
de $H_{\g_0, \rho_0}$ (resp. $H_{\g, \rho}$), on cherche \`a quelles conditions
$H_0$ est adh\'erente \`a $H$ (i.e. $H_0 \subset \ov {H}$), voire faiblement
adh\'erente (i.e. $H_0 \cap \ov {H} \neq \vide$).

 Cette derni\`ere condition
signifie simplement qu'il existe une famille de courbes param\'etr\'ee par un
anneau de valuation discr\`ete $A$, dont le point sp\'ecial 
$C_0$ est dans $H_0$  et dont le point g\'en\'erique $C$ est dans $H$.

Une question fondamentale, mais sans doute difficile,  est de donner des
conditions n\'ecessaires et suffisantes \`a l'existence de telles familles.
Il y a, en tous cas,  deux conditions n\'ecessaires.
La premi\`ere est une condition   de  semi-continuit\'e  sur les
dimensions des espaces de cohomologie : si
$C_0$ est une sp\'ecialisation de
$C$ on a les in\'egalit\'es $h^i
\cJ_C (n)
\leq h^i
\cJ_{C_0} (n)$ pour tout $i$ et tout $n$. La deuxi\`eme concerne les structures des
modules de Rao : on montre (cf. 5.9 ci-dessous) que le module
$M_C$ est  une
d\'eformation plate d'un 
sous-quotient 
 (i.e. 
un quotient d'un sous-module), de
$M_{C_0}$. Cette condition est \`a la base de notre approche
 de la question, {\it via} les modules de Rao.

\vskip 0.3 cm

Dans le cas des courbes sur un corps on sait, \`a partir d'un module de longueur finie
$M$, d\'ecrire les courbes (notamment
minimales) de la classe de biliaison qu'il d\'efinit. Cela repose essentiellement
sur le calcul de la fonction
$q$   associ\'ee \`a $M$, cf. [MDP1] IV.   Notre objectif ici est de g\'en\'eraliser ce
processus pour construire des familles de courbes.

Pour cela, il faut  disposer d'une notion
qui g\'en\'eralise au cas des familles de courbes celle de module de Rao ordinaire
(i.e. de module de longueur finie). On est ainsi \`a la recherche, pour tout anneau
noeth\'erien $A$, d'un ensemble
$\Theta (A)$   dont les \'el\'ements
 permet\-tent de d\'ecrire les familles de modules de Rao de dimensions
variables param\'etr\'ees par $A$.  On doit aussi disposer d'une application
$\Psi_A$ de
$H_{d,g}(A)$ dans $\Theta (A)$ qui associe
\`a une famille de courbes $\cC$ param\'etr\'ee par $A$ un \'el\'ement $\Psi_A (\cC)$ de
$\Theta (A)$ qui d\'ecrive la variation des modules de Rao des courbes de la famille. 
Plus pr\'ecis\'ement, si on s'inspire du r\^ole que jouent les modules de Rao
ordinaires vis \`a vis  des courbes sur un corps on attend d'une telle application les 
propri\'et\'es suivantes :

1) une propri\'et\'e de surjectivit\'e : tout \'el\'ement  de $\Theta (A)$,
v\'erifiant des conditions conve\-na\-bles est, \`a d\'ecalage pr\`es, image par
$\Psi_A$ d'une famille de courbes (dans le cas ordinaire, tout module de longueur finie
est, \`a d\'ecalage pr\`es, le module de Rao d'une courbe : ceci est l'une des
assertions du th\'eor\`eme de Rao, cf. [R]), 

2) une description des fibres de $\Psi_A$ : deux familles de courbes ayant m\^eme image
par $\Psi_A$ \`a d\'ecalage pr\`es, sont dans la m\^eme classe de biliaison ou
liaison paire (dans le cas ordinaire c'est le deuxi\`eme point du th\'eor\`eme de
Rao, cf. [R]),

3) une description des familles minimales : en compl\'ement du point 1), un
\'el\'ement  de $\Theta (A)$
\'etant donn\'e, il s'agit de d\'ecrire les familles minimales (au sens du moindre
d\'ecalage) qui lui correspondent (dans le cas ordinaire c'est le calcul des
courbes minimales effectu\'e dans [MDP1] IV via la fonction $q$) et de d\'ecrire
les autres familles (et en particulier de pr\'eciser les d\'ecalages possibles) 
\`a partir de celles-ci (c'est l'analogue pour les familles du th\'eor\`eme de
Lazarsfeld-Rao  qui montre qu'on passe des courbes minimales
aux autres par des biliaisons
\'el\'ementaires ascendantes, cf. [LR], [BBM]  ou [MDP1]  IV),

4)  une notion de dualit\'e sur les
\'el\'ements  de $\Theta (A)$ qui corresponde sur les courbes \`a la notion de liaison
impaire, (g\'en\'eralisation, l\`a encore, du th\'eor\`eme de Rao),

5) enfin, une  \'etude de l'application $\Psi$  qui g\'en\'eralise dans la mesure du
possible les assertions de lissit\'e et d'irr\'eductibilit\'e formelles du th\'eor\`eme
de l'\'etape interm\'ediaire (cf. [MDP1] VII 1.1 et 1.5). Cette \'etude, qui
n\'ecessite de munir $\Theta$ d'une structure alg\'ebrique,  passe  par la
g\'en\'eralisation au cas des familles des notions de r\'esolutions de type {\sl N} et
{\sl E}, cf. [MDP1].

\vskip  0.3 cm

L'objectif de cet article est de poser les bases de la th\'eorie en
d\'efinissant la notion de triade qui va \^etre la g\'en\'eralisation de celle de
module de Rao ($\Theta (A)$ sera l'ensemble des triades sur $A$ modulo une
certaine relation de pseudo-isomorphisme), en d\'efinissant  la fl\`eche
$\Psi_A$ et en montrant, sur un anneau local, la plupart des propri\'et\'es attendues~:
la surjectivit\'e 1), le lien avec la biliaison 2), la description des familles minimales
et le th\'eor\`eme de Lazarsfeld-Rao 3), la dualit\'e et son lien avec la liaison
impaire 4) et en introduisant les notions de faisceau triadique et de r\'esolutions
(triadiques) de type {\sl N} et {\sl  E} de  5).

Via les faisceaux triadiques qui jouent un r\^ole d'interm\'ediaire entre courbes
et triades, les d\'emonstrations des points 1)
\`a 4) ci-dessus reposent en grande partie sur nos deux articles [HMDP1] et [HMDP2].

\vskip 0.5 cm

Commen\c cons par expliquer  les id\'ees qui conduisent \`a la
notion de triade.

 Pour d\'efinir
$\Theta(A)$ et $\Psi_A$,  une id\'ee (trop) simple consiste \`a  prendre pour 
$\Theta(A)$  l'ensemble des classes d'isomorphisme de $R_A =A[X,Y,Z,T]$-modules gradu\'es 
$M_A=
\bigoplus M_{A,n}$ tels que les $M_{A,n}$ soient de type fini sur $A$ et presque tous nuls, puis 
\`a prendre comme application $\Psi_A $ celle qui associe
\`a une famille de courbes $\cC$ la famille de modules
$H^1_*\cJ_\cC= \bigoplus_{nÊ\in \bZ} H^1(\bP_A, \cJ_{\cC}(n))$. 

Cette formule permet de d\'efinir, sur le sous-sch\'ema
 des courbes \`a cohomologie cons\-tante, le morphisme \'evoqu\'e plus haut $\Phi: H_{\g,\rho} \fl
E_{\rho}$,  mais elle n'est pas valable sur
$H_{d,g}$ tout entier. La raison fondamentale de ce fait est que si $\cC$ est une famille de
courbes  qui n'est pas 
\`a sp\'ecialit\'e constante 
\footnote{$(^1)$}{La premi\`ere approche de ce type
de questions, semble avoir
\'et\'e effectu\'ee, dans ce cas particulier de la sp\'ecialit\'e constante, par
Ballico et Bolondi, cf. [BB].} 
 le $A$-module $H^1_*\cJ_\cC$ ne commute pas au changement de base, c'est-\`a-dire
que le module de Rao d'une fibre $\cC \T_A k(t)$ n'est pas le module  $H^1_* \cJ_C
\T_A k(t)$, voir l'exemple de la famille $\cC_2$ ci-dessous. On ne peut donc se
contenter pour d\'ecrire la variation du module de Rao dans une famille de courbes de
consid\'erer un
$R_A$-module global
$M_A$  et ses fibres $M_A \T_A k(t)$.

L'id\'ee que nous proposons  pour surmonter cette difficult\'e s'inspire de ce que
l'on fait pour prouver les th\'eor\`emes de cohomologie et changement de base, cf. par
exemple [AG] III \S 12. Elle consiste, dans un premier temps, 
\`a regarder les {foncteurs} $V: Q \mapsto V(Q)$, de la cat\'egorie des $A$-modules
dans la cat\'egorie des $R_A$-modules gradu\'es (en pensant notamment comme
$A$-modules aux corps r\'esiduels $k(t)$ pour $t \in T$). Bien
s\^ur, un $A$-module $M_A$ fournit un cas particulier de tel foncteur, celui
qui associe \`a $Q$ le module gradu\'e $M_A\T_A Q$, mais plus g\'en\'eralement,
\`a  une famille de courbes
$\cC$ quelconque est  associ\'e   le foncteur de Rao $V_\cC : Q
\mapsto  \bigoplus_{n \in \bZ} H^1(\bP_A, \cJ_{\cC} (n) \T_A Q)$ qui d\'ecrit en
particulier la variation du module de Rao en les points de $T$.

Bien entendu, si l'on ne fait pas d'hypoth\`eses compl\'ementaires, cette notion de
foncteur est
 trop grossi\`ere. Nous allons  examiner quelques
exemples de familles de courbes et, pour d\'eterminer les foncteurs associ\'es,  tenter
de comprendre comment varient leurs r\'esolutions.

Les deux premiers exemples concernent le sch\'ema $H_{6,3}$. Ce sch\'ema contient une
composante dont la courbe g\'en\'erale $C$  est ACM (arithm\'etiquement de
Cohen-Macaulay c'est-\`a-dire dont le module de Rao est nul). 

 Il contient aussi deux autres sch\'emas $H_{\g,\rho}$ qui sont tous deux form\'es de
courbes  de la classe de liaison de la r\'eunion de deux droites
disjointes (i.e. les courbes dont le module de Rao est
\'egal \`a
$k= R/(X,Y,Z,T)$ en un unique degr\'e). Le premier est le sch\'ema $H_1$ dont la
courbe g\'en\'erale
$C_1$ est une courbe de bidegr\'e
$(4,2)$ sur une quadrique, de module de Rao $k(-2)$ (le module \'egal \`a $k$ mais concentr\'e en
degr\'e
$2$). Le deuxi\`eme est le sch\'ema  $H_2$  qui a pour courbe g\'en\'erale $C_2$ la
r\'eunion d'une quartique plane et de deux droites qui la coupent chacune en un
point, avec comme module de Rao $k(-1)$.

Il y a une famille
  $\cC_1$,
\`a sp\'ecialit\'e constante, param\'etr\'ee par un \avd\ d'uniformisante $a$ qui
joint
$C$ et
$C_1$. Cette famille admet la r\'esolution ``de type {\sl N}'' (cf. [MDP1] VII 2.6) 
suivante :
$$ 0 \fl R_A(-4)^3 \fl [R_A(-3)^4 \oplus R_A(-2) \Fl {d_1} R_A(-2)] \fl I_{\cC_2}
\fl 0$$
avec $d_1= (X,Y,Z,T,a)$ (on  note $[A \Fl {u} B]$ le noyau de $u$). On voit aussit\^ot que pour
$a$ inversible, c'est-\`a-dire au point g\'en\'erique de $\Spec A$, le facteur $R_A(-2)$ se
simplifie, donnant la courbe ACM
$C$. Dans ce cas, la variation du module de Rao est d\'ecrite par le module global 
$H^1_*\cJ_{\cC_1}=M_A=
\Coker d_1 = R_A/(a,X,Y,Z,T)(-2)$.

Il y a aussi une famille
$\cC_2$, 
\`a postulation constante, qui joint $C$ et $C_2$. Cette fois, la famille admet la
r\'esolution ``de type {\sl E}'' suivante (avec le m\^eme $d_1$) :
$$ 0 \fl R_A(-5) \Fl{^td_1} R_A(-4)^4 \oplus R_A(-5) \fl R_A(-3)^4 \oplus R_A(-4)
\fl I_{\cC_2} \fl 0$$
(cf. [MDP1] VII 2.1) et on voit, l\`a encore, la simplification  du terme $R_A(-5)$ au point
g\'en\'erique.

 Cependant, comme cette famille v\'erifie $M_A=H^1_*\cJ_{\cC_2}=0$ alors
que certaines courbes de la famille sont dans la classe de deux droites, donc
ont un module de Rao non nul, il est clair que la variation du module de Rao ne peut \^etre
d\'ecrite par le module $M_A$.

\vskip 0.3 cm

Dans les deux exemples pr\'ec\'edents la famille $\cC_1$ (resp. $\cC_2$) admet une
r\'esolution  \`a trois termes, de type {\sl N} (resp.  {\sl E})\footnote
{$(^2)$}{mais pas l'inverse !} et le passage de la courbe sp\'eciale \`a la courbe
g\'en\'erique correspond \`a une simplification
\'evidente de la r\'esolution. L'exemple suivant montre que les choses peuvent \^etre
plus complexes lorsque la famille n'est ni \`a sp\'ecialit\'e, ni \`a postulation
constante.

Consid\'erons le sch\'ema $H_{4,0}$. Il poss\`ede deux composantes, toutes deux de
dimension $16$, l'une
$H_1$ dont la courbe g\'en\'erale $C_1$ est une courbe de bidegr\'e $(3,1)$ sur une
quadrique, de module de Rao $k(-1)$, l'autre $H_2$ dont la courbe g\'en\'erale $C_2$
est r\'eunion disjointe d'une cubique plane et d'une droite et a un module de Rao
 du type de $R/(X,Y,Z,T^3)$ de dimensions $1,1,1,$ en degr\'es $0,1,2$.\footnote
{$(^3)$}{On notera que $k(-1)$ est un sous-quotient de ce module.} Voici les
r\'esolutions de type {\sl E} de ces courbes :
$$ 0 \fl R(-5) \fl R(-4)^4 \fl  R(-2)\oplus R(-3)^3 \fl I_{C_1} \fl 0$$
$$ 0 \fl R(-6) \fl R(-3) \oplus R(-5)^3 \fl R(-2)^2 \oplus R(-4)^2 \fl I_{C_2} \fl
0.$$
Bien que les composantes $H_1$ et $H_2$ v\'erifient les conditions n\'ecessaires sur la
cohomologie et le module de Rao pour que
$H_2 \cap \ov {H_1}$ soit non vide, il n'est nullement
\'evident, au vu de ces r\'esolutions, pas plus d'ailleurs qu'avec celles de type {\sl N}, de
dire s'il peut exister une famille de courbes  de point
g\'en\'erique dans $H_1$ et de point sp\'ecial dans $H_2$.  En particulier, \`a
l'inverse des exemples ci-dessus, aucune simplification n'est apparente. 
Pourtant, nous verrons au \S5 qu'il existe bien une telle famille.

On  comprend un peu mieux
ce ph\'enom\`ene si on note qu'il y a dans l'id\'eal de $C_2$ une
\'equation de degr\'e
$2$ de plus que dans celui de $C_1$ et si on   supprime cette \'equation en {\bf d\'esaturant}
l'id\'eal
$I_{C_2}$. Par exemple, si $C_2$ est r\'eunion de $(X,Y)$ et de $(Z,T^3)$ \footnote
{$(^4)$}{Cet exemple est choisi parce que le calcul est facile, mais les courbes de $H_2$
r\'eunions disjointes d'une cubique et d'une droite ne sont pas dans l'adh\'erence de $H_1$.
Il faut utiliser des structures multiples, cf. 5.22, mais le
calcul de la r\'esolution du d\'esatur\'e est identique.}, son id\'eal est
$(XZ,YZ,XT^3,YT^3)$ et, si on remplace $YZ$ par les \'equations $XYZ,
Y^2Z,YZ^2,YZT$, ce qui ne change pas $C_2$, on obtient un id\'eal $J$ non satur\'e
dont la r\'esolution a {\bf quatre} termes est :
$$0 \fl R(-6) \fl R(-5)^4 \oplus R(-6) \fl R(-4)^6 \oplus R(-5)^3 \fl R(-2) \oplus
R(-3)^3 \oplus R(-4)^2
\fl J \fl 0$$
et l\`a, on voit, apr\`es simplifications d'un $6$, de trois $5$ et de deux
$4$ appara\^\i tre  les chiffres de la r\'esolution de $I_{C_1}$ dans celle de $J$
(bien entendu, cf. la note
$(^4)$, ce calcul ne prouve pas l'existence d'une famille mais il donne un indice
num\'erique favorable).

\vskip 0.3 cm

L'id\'ee fondamentale que nous voulons \'eclairer par cet exemple c'est qu'on ne
saurait comprendre les familles g\'en\'erales de courbes (celles qui ne sont ni \`a
sp\'ecialit\'e, ni
\`a postulation constante)  si l'on ne tol\`ere pas cette op\'eration de
d\'esatura\-tion (cf. 2.12), avec comme cons\'equence in\'eluctable  l'apparition de
r\'esolutions
\`a quatre termes des id\'eaux. 

\`A cet \'egard, nous d\'efinissons ci-dessous des r\'esolutions de type \sE\
``cotriadiques'' du faisceau d'id\'eaux $\cJ_\cC$, cf. 3.1. Il s'agit de r\'esolutions
de la forme $0 \fl
\cE \fl \cF \fl \cJ_\cC \fl 0$ o\`u $\cE$ est localement libre et d\'efini par une
suite exacte $0 \fl \cL_5 \fl \cL_4 \fl \cL_3 \fl \cE \fl 0$ avec $\cF$ et les $\cL_i$
dissoci\'es (i.e., dans le cas d'un anneau local, somme directe finie de faisceaux 
inversibles). Nous utilisons aussi la notion duale de r\'esolution de type
\sN\ ``triadique''  de la forme :
$0 \fl \cP \fl \cN \fl \cJ_\cC \fl 0$, o\`u $\cN$ est d\'efini par une suite exacte
$0 \fl \cN \fl \cL_1 \fl \cL_0 \fl \cL_{-1} \fl 0$ avec $\cP$ et les $\cL_i$
dissoci\'es. Par rapport \`a une
r\'esolution de type  \sE\ (resp. \sN) sur un corps, on notera la pr\'esence du terme
suppl\'ementaire  $\cL_5$  (resp. $\cL_{-1}$).

Nous montrons dans ce travail que si $\cC$ est une famille de courbes sur un anneau local
$A$ il existe effectivement des r\'esolutions de type \sE\ et \sN\ triadiques (cf. 3.1).
L'int\'er\^et de ces r\'esolutions est de  r\'epondre \`a la question de la nature du
foncteur
$V_\cC$ d\'efini ci-dessus et qui d\'ecrit la variation du module de Rao. On montre en
effet par un petit calcul  cohomologique que, si 
$\Lpont = (L_1 \fl L_0 \fl L_{-1})$ est le complexe de $R_A$-modules
gradu\'es associ\'e
\`a une r\'esolution de type \sN, on a $V_\cC(Q) = h_0 (\Lpont \T Q)$, o\`u
$h_0$ d\'esigne l'homologie en degr\'e $0$ du complexe. On notera que c'est
aussi cette homologie qui appara\^\i t dans les th\'eor\`emes de cohomologie et  changement
de base. 

Le complexe $\Lpont$ ci-dessus est appel\'e une {\bf triade}. C'est cette notion  qui
g\'en\'eralise celle de module de Rao et l'application
$\Psi_A $ consiste donc  \`a remplacer une famille de courbes, c'est-\`a-dire, via les
r\'esolutions de type \sN, un complexe
\`a quatre termes, par une triade, c'est-\`a-dire un complexe 
\`a trois termes.

\vskip 0.5 cm

Examinons maintenant le plan de ce travail et les r\'esultats obtenus.

\vskip 0.2 cm

Le paragraphe 1 est consacr\'e \`a mettre en place et \`a \'etudier la notion de
triade. 

Une triade est un complexe  de $R_A$-modules gradu\'es $\Lpont =(L_1 \Fl {d_1} L_0 \Fl
{d_0} L_{-1})$ avec notamment des hypoth\`eses de finitude sur ses groupes d'homologie,
cf. 1.10 pour une d\'efinition pr\'ecise. Un cas  particulier important est celui
 des triades majeures (cf. 1.13).

On d\'efinit le foncteur  $V$  associ\'e \`a une triade (cf. 1.3) : c'est le
foncteur $ h_0 (\Lpont \T \tas)$ (cf. aussi 1.36).

On d\'efinit ensuite un
pseudo-isomorphisme entre deux triades   (en abr\'eg\'e un \ps, cf. 1.7) : c'est   un
morphisme de complexes
 qui induit  un isomorphisme sur les foncteurs
$V$ et un monomorphisme  sur les foncteurs $h_{-1} (\Lpont \T \tas)$ et on dit que deux
triades sont pseudo-isomorphes si elles sont reli\'ees par une cha\^\i ne de \ps. On
\'etablit le lemme  1.19 qui permet de se limiter \`a des cha\^\i
nes
\`a deux maillons (\`a la mani\`ere de Verdier, cf. [V]).

Deux
triades pseudo-isomorphes d\'efinissent le m\^eme foncteur, mais la
r\'eciproque est fausse (cf. 1.35.c). 
 Si on est sur un
anneau de valuation discr\`ete toute triade est pseudo-isomorphe \`a une unique triade
majeure  dite \'el\'ementaire (cf. 1.32), caract\'eris\'ee par le fait que son conoyau
$h_{-1} (\Lpont)$ est de torsion. 

On d\'efinit  en 1.26 la notion de triade duale (d\'efinie \`a \ps\ pr\`es).
\vskip 0.3 cm

 Le paragraphe 2 est consacr\'e \`a l'\'etude des faisceaux triadiques. Il s'agit
(cf. 2.3) des  faisceaux
$\cN$ qui s'ins\`erent dans une suite exacte $ 0 \fl \cN \fl \cL_1 \fl \cL_0 \fl
\cL_{-1}
\fl 0$, avec les $\cL_i$ dissoci\'es. Il y a une \'equivalence de
cat\'egories entre la cat\'egorie des triades majeures, \`a homotopie pr\`es, et celle
des faisceaux triadiques, la correspondance envoyant une triade
$\Lpont$ sur le faisceau $\cN$ associ\'e \`a son noyau $N= h_1(\Lpont)$. On montre
(l\`a encore \`a l'aide d'un lemme de Verdier 2.11), que les \ps\ de triades et ceux de
faisceaux (d\'efinis comme dans [HMDP1]) se correspondent (cf. 2.14). Dans le cas d'un
\avd, on montre que les triades \'el\'ementaires correspondent aux faisceaux
extravertis minimaux de [HMDP1], cf. 2.15.

\vskip 0.3 cm

Le paragraphe 3 fait le lien entre  familles de
courbes, triades et faisceaux triadiques. On commence par
montrer, par deux m\'ethodes (l'une par d\'esaturation et liaison l'autre par troncature \`a partir  du complexe 
$\bR\G_*\cJ_{\cC}$),  l'existence de
r\'esolutions de type
\sN\ triadiques (i.e. dont le faisceau $\cN$ est triadique). On d\'efinit ensuite (cf.
3.6) une triade de Rao associ\'ee \`a une famille de courbes comme  la triade associ\'ee 
au faisceau $\cN$ d'une
r\'esolution de type
\sN\ triadique de
$\cC$. Cette  triade  est bien
d\'efinie
\`a
\ps\ pr\`es.  Le foncteur associ\'e n'est autre que le foncteur de Rao $V_\cC$
d\'efini ci-dessus. 
\vskip 0.2 cm

 On obtient ainsi la
fl\`eche
$\Psi_A$
\'evoqu\'ee plus haut qui associe
\`a une famille de courbes $\cC$ la classe de triades de ses triades de Rao
$\cT(\cC)$ (pour la relation de pseudo-isomorphisme) et on prouve, \`a l'aide de
[HMDP1], le r\'esultat suivant :

\th {Th\'eor\`eme 3.9}. (Th\'eor\`eme de Rao pour les triades)  On suppose $A$ local
\`a corps r\'esiduel infini.  Soient
$\cC$ et
$\cC'$ deux familles plates de courbes param\'etr\'ees par  $A$ et soient $\Lpont$ et
$\Lpont'$ des triades de Rao de $\cC$ et $\cC'$. Alors,
$\cC$ et
$\cC'$ sont dans la m\^eme classe de biliaison si et seulement si les triades $\Lpont$ et
$\Lpont'$  sont pseudo-isomorphes,  \`a d\'ecalage pr\`es (i.e. s'il existe un entier
$h$ tel que l'on ait
$\cT(\cC) = \cT(\cC') (h)$).

Puis,
gr\^ace aux r\'esultats de [HMDP2]  exprim\'es en termes de triades, on
 montre la g\'en\'eralisation de l'assertion de surjectivit\'e (toujours \`a
d\'ecalage pr\`es) du th\'eor\`eme de Rao (i.e. de l'application $\Psi_A$) :

\th {Th\'eor\`eme 3.10}. On suppose $A$ local \`a corps r\'esiduel infini.  Soit
$\Lpont$ une triade. Il existe une famille de courbes
$\cC$ telle que la triade de Rao de $\cC$ soit pseudo-isomorphe \`a $\Lpont$, \`a
d\'ecalage pr\`es.

Plus pr\'ecis\'ement, le th\'eor\`eme suivant donne exactement l'image de $\Psi_A$, i.e.
les d\'ecalages possibles. Il contient aussi la g\'en\'eralisation du th\'eor\`eme de
Lazarsfeld-Rao~:

\th {Th\'eor\`eme 3.11}. On suppose $A$ local \`a corps r\'esiduel infini. Soit $\Lpont$
une triade, $\cT$ sa classe de pseudo-isomorphisme, $\cN$ un faisceau triadique
associ\'e et soit
$\cN_0$ le faisceau extraverti minimal (unique) de la classe de pseudo-isomorphisme de
$\cN$ (cf. [HMDP1] 2.14). 
Soit  $q=
q_{\cN_0}$ la fonction $q$ de $\cN_0$ (cf. [HMDP2] 2.4) et soit  $h_0= \sum_{n \in \bZ}
nq(n) +
\deg
\cN_0$. \lign 1) Il
existe une famille de courbes 
$\cC_0$ et une r\'esolution 
$ 0 \fl \cP_0 \fl \cN_0 \fl \cJ_{\cC_0} (h_0) \fl 0$ avec $\cP_0$ dissoci\'e. La
classe de triade associ\'ee \`a $\cC_0$ est \'egale \`a $\cT(-h_0)$.  \lign 
2) Si $\cC_1$ est
une famille de courbes telle que $\cT(\cC_1) = \cT (-h)$,   on a
$h
\geq h_0$.  Si
$d$ et
$g$ (resp.
$d_0$ et
$g_0$) sont respectivement le degr\'e et le genre de $\cC_1$ (resp. $\cC_0$) on a $d \geq
d_0$ et
$g
\geq g_0$.\lign
3) R\'eciproquement, pour tout $h \geq h_0$ il existe une famille de courbes $\cC_1$ avec
$\cT(\cC_1) = \cT(-h)$. \lign 
4) Si de plus on a
$h=h_0$,
$\cC_1$ et $\cC_0$ sont jointes par une d\'eformation  \`a cohomologie uniforme 
(cf. [HMDP2] 2.8) et triade constante. 
\lign On dit que
$\cC_0$ est une {\bf famille minimale} de courbes. \lign
Les autres familles  de
courbes de la classe de biliaison s'obtiennent \`a partir de $\cC_0$ par des biliaisons
\'el\'ementaires  suivies d'une d\'eformation \`a cohomologie uniforme et triade
constante.
 
Enfin, on a le lien entre dualit\'e et liaison impaire :

\th {Th\'eor\`eme 3.13}.  On suppose $A$ local \`a corps r\'esiduel infini. Soient $\cC$
et $\cC'$ deux familles plates de courbes param\'etr\'ees par  $A$, et soient $\Lpont$ et
$\Lpont'$ des triades de Rao de $\cC$ et $\cC'$. Alors, $\cC$ et
$\cC'$ sont li\'ees par un nombre impair de liaisons \'el\'ementaires si et
seulement si les triades
$\Lpont$ et
$\Lpont'$ sont duales \`a \ps\ pr\`es et \`a  d\'ecalage pr\`es, i.e. s'il existe un
entier
$h$ tel que
$\Lpont (h)$ soit  pseudo-isomorphe \`a  une triade duale de $\Lpont'$.

\vskip 0.3 cm

Le paragraphe 4 \'etablit un dictionnaire reliant courbes, faisceaux et triades. On y
 caract\'erise les triades dites modulaires (celles qui sont  d\'efinies par un 
$R_A$-module) et on montre qu'elles correspondent aux familles de courbes \`a
sp\'ecialit\'e constante. On d\'efinit aussi la notion duale de
triade repr\'esentable qui correspond aux familles de courbes \`a postulation constante.

\vskip 0.3 cm

 Le paragraphe 5 est essentiellement constitu\'e d'exemples. On y reprend,
comme fil conducteur, le troisi\`eme des exemples \'etudi\'es ci-dessus : les courbes
de degr\'e $4$ et genre $0$. 

On  montre 
(cf. 5.8) que si on a une triade  sur un \avd, la valeur du foncteur $V$ au point
g\'en\'erique  est une d\'eformation plate d'un sous-quotient   de la valeur de $V$ au
point ferm\'e. Il en r\'esulte  (cf. 5.9) qu'on a ce m\^eme  r\'esultat pour  les modules
de Rao aux points g\'en\'erique et ferm\'e d'une famille de courbes.

On montre, r\'eciproquement, (cf. 5.15) que si on a deux modules gradu\'es $M$ et
$M_0$ de longueur finie tels que
$M$ soit un sous-quotient de
$M_0$ il existe une triade dite triviale, param\'etr\'ee par un \avd, dont les
valeurs du foncteur associ\'e au point ferm\'e (resp. g\'en\'erique) sont $M_0$
(resp.
$M$).  Ce r\'esultat permet de construire de nombreux exemples de triades et de
d\'eterminer les familles minimales associ\'ees gr\^ace \`a l'algorithme
d\'evelopp\'e dans [HMDP2] \S 3 (cf. 5.17, 5.18).

Cependant on constate dans de nombreux cas (et notamment celui qui est cens\'e mener \`a
une famille
$(4,0)$) que cette construction
 n'est
pas optimale, au sens o\`u elle ne donne pas des familles de  degr\'e assez petit. On
reprend donc la construction pr\'ec\'edente avec plus de soin en montrant comment
construire de ``meilleures'' triades. Cette m\'ethode conduit en
particulier
\`a la construction d'une famille de courbes de $H_{4,0}$ et \`a la preuve de la
connexit\'e du sch\'ema de Hilbert
$H_{4,0}$ (cf. 5.21).

Plus g\'en\'eralement, les applications \`a l'\'etude du sch\'ema de Hilbert
des courbes de la notion de triade semblent prometteuses, car elles donnent un
moyen syst\'ematique de construire des familles de courbes dans tous les cas et en
particulier le cas non trivial o\`u ces familles ne sont ni \`a postulation, ni \`a
sp\'ecialit\'e constante. Ces constructions jouent notamment un r\^ole essentiel  dans
les probl\`emes de connexit\'e, cf. par exemple [AA].

\titre {0. Notations et pr\'eliminaires}

 Si $A$ est un anneau noeth\'erien, on note $\bP^3_A$ l'espace projectif
de dimension $3$ sur $A$ et  $R_A$ l'anneau $A[X,Y,Z,T]$. Si $\cF$ est un faisceau
coh\'erent sur
$\bP^3_A$ on note
 $H^i \cF$ le $A$-module $ H^i(\bP^3_A, \cF)$. 
On pose $\G_* \cF = \bigoplus _{n\in \bZ} \G
\cF (n)$ et $H^i_* \cF=\bigoplus _{n\in \bZ} H^i
\cF (n)$ pour $i \geq 0$. Lorsque $A=k$ est un corps on pose $R=
R_k$, $\bP^3=\bP^3_k$ et  on note
$h^i
\cF$ la dimension de $H^i \cF$. 

\tarte {a) Modules gradu\'es}

Soit $A$ un anneau noeth\'erien. On d\'esigne par
${\Mod}_A$  la cat\'egorie des $A$-modules et par ${\Gr}_{R_A}$ la cat\'egorie des
$R_A$-modules gradu\'es, avec comme morphismes les homomorphismes de degr\'e z\'ero.
Ce sont des cat\'egories ab\'eliennes.

Si $M$ est un $R_A$-module gradu\'e
le module d\'ecal\'e $M(d)$ est le module $M$ muni de la graduation d\'efinie par
$M(d)_n = M_{n+d}$.

Si $M$ est un $R_A$-module gradu\'e et $r$ un entier, 
$M_{>r} = \bigoplus _{n>r}M_n$
est un sous-module gradu\'e de $M$.  Le
quotient $M/M_{>r}$ est not\'e
$M_{\leq r}$.   C'est aussi un $R_A$-module gradu\'e,  dont les composantes de
degr\'e $\leq r$ sont isomorphes \`a celles de $M$ et qui est nul en degr\'e
$>r$.  Les foncteurs $M \mapsto M_{>r}$ et $M \mapsto M_{\leq r}$
de
$\Gr_{R_A}$ dans elle-m\^eme sont appel\'es foncteurs de troncature et ils sont exacts.

Un $R_A$-{\bf module 
dissoci\'e} est un module de la forme
$$F = \bigoplus_{i=1}^r M_i \T_A  R_A(-n_i), $$
o\`u les $ n_i$ sont des entiers et les $M_i$ des $A$-modules projectifs de type fini.
Si $A$ est local tout module dissoci\'e est libre. On v\'erifie que les
modules dissoci\'es sont des objets projectifs dans la cat\'egorie $\Gr_{R_A}$.

De m\^eme, on introduit la notion de {\bf faisceau dissoci\'e}
sur $\bP^3_A$ : c'est un faisceau de la forme 
$$ \cF = \bigoplus_{i=1}^r M_i \T_A \cO_{\bP_A}(-n_i) $$
o\`u les $ n_i$ sont des entiers et les $M_i$ des $A$-modules projectifs de type fini. 
Si le module $F$ est dissoci\'e le faisceau associ\'e  $\widetilde F$ l'est aussi ; si
le faisceau
$\cF$ est dissoci\'e le module $H^0_* \cF$ l'est aussi.

Si $M$ est un $R_A$-module gradu\'e et $Q$ un $A$-module, on a  une structure de
$R_A$-module gradu\'e sur $M\T_A Q$ en posant $(M \T_A Q)_n= M_n \T_A Q$. Sous ces
m\^emes hypoth\`eses on pose 
$$ \Homgr_{A} (M,Q) = \bigoplus_{n \in \bZ} \Hom_A(M_{-n}, Q)$$
et on d\'efinit une structure de $R_A$-module gradu\'e sur  $\Homgr_{A} (M,Q)$
en posant
 $(\la.f) (x) = f(\la.x)$.

Si $M$ et $N$ sont des $R_A$-modules gradu\'es et $d$ un entier on note
$\Hom^d_{R_A}(M,N)$ l'espace des homomorphismes de $R_A$-modules, homog\`enes de degr\'e
$d$ et  on pose
$$\Homgr_{R_A}(M,N) = \bigoplus_{d \in \bZ} \Hom^d_{R_A}(M,N).$$

\tarte {b) Modules duaux}

Nous rappelons ici les deux notions de modules duaux d'un $R_A$-module gradu\'e $M$, les
modules
$M^\vee$ et $M^*$ qu'on prendra garde \`a ne pas confondre. 

Si $M$ est un $R_A$-module gradu\'e, on d\'efinit son dual (sur $R_A$) comme le module
 $M^\vee = \Homgr_{R_A} (M, R_A)$. Si $M$ est un $R_A$-module de type fini, $M^\vee$
est simplement  le module $\Hom_{R_A} (M, R_A)$. Si
$M$ est libre gradu\'e,
$M= \bigoplus_{i=1}^r R_A(-n_i)$, on a 
$M^\vee = \bigoplus_{i=1}^r R_A(n_i)$.

\vskip 0.3 cm

Si $M$ est un $R_A$-module gradu\'e on d\'efinit le module dual (sur $A$)
$M^*$ de
$M$ comme $\Homgr_{A} (M,A)$. On a donc
$M^* = \bigoplus_{n \in \bZ} (M_n)^*$, avec $(M_n)^*=\Hom_A (M_n, A)$ ; ce module est
gradu\'e  par 
$(M^*)_{n} =
\Hom_A (M_{-n}, A)$ et sa structure de $R_A$-module
d\'efinie par
$\la f (x) = f (\la x)$. 

Un cas particulier de module dual est le module $R_A^*$ qui est gradu\'e en degr\'es
$\leq 0$. Pr\'ecis\'ement, $({R_A}^*)_{-n}$ est un $A$-module libre de dimension
${n+3 \choose 3}$. On notera les formules $(R_A(d))^* = R_A^*(-d)$ et 
$(M_{\leq r})^* = (M^*)_{\geq -r}$.

\vskip 0.3 cm

L'application $M \mapsto M^*$ d\'efinit un foncteur contravariant de la
cat\'egorie
$\Gr_{R_A}$ dans elle-m\^eme. De plus, si l'on se restreint aux
$R_A$-modules gradu\'es dont  les composantes de tous degr\'es sont des $A$-modules
projectifs de type fini, le foncteur $M \mapsto M^*$  est 
exact et involutif. 

\vskip 0.3 cm

La proposition suivante pr\'ecise les liens entre les modules duaux et la
cohomologie sur l'espace projectif $\bP^3_A$ dans le cas dissoci\'e :

\th {Proposition 0.1}. \lign
1)  Soit $L$ un $R_A$-module dissoci\'e, $\cL$ le faisceau sur $\bP^3_A$ associ\'e. On
a $  H^0_* (\bP^3_A, \cL) = L$. \lign
2) Soit $\cL$ un faisceau dissoci\'e et soit $L=H^0_* (\bP^3_A, \cL)$.
On a $H^0_* (\bP^3_A, \cL^\vee) =L^\vee$.

\tarte {c) Foncteurs}

On consid\`ere ici, pour l'essentiel, des foncteurs covariants de la cat\'egorie
$Mod_A$   dans elle-m\^eme ou dans la cat\'egorie
$\Gr_{R_A}$. Un tel foncteur
$V : \Mod_{A} \fl \Gr_{R_A}$ est une somme directe  de
foncteurs $V_n : \Mod_A \fl \Mod_A$ et est muni d'une op\'eration gradu\'ee de $R_A$
(c'est-\`a-dire, pour tout $A$-module $Q$, de morphismes $R_A \times V(Q) \fl V(Q)$
homog\`enes de degr\'e
$0$ et  fonctoriels en
$Q$).  Si $V$ est un tel
foncteur, le foncteur d\'ecal\'e $V(h)$ est d\'efini par la formule 
$[V(h) (Q)]_n= V(Q)_{n+h}$.

Les foncteurs consid\'er\'es seront toujours suppos\'es
$A$-lin\'eaires.

Nous utiliserons la th\'eorie des foncteurs coh\'erents due \`a Auslander (cf. [H]).
Un foncteur $A$-lin\'eaire $V : \Mod_A \fl \Mod_A$ est dit {\bf coh\'erent } s'il est
conoyau d'un morphisme de foncteurs 
$\Hom_A(M, \tas) \fl \Hom_A(N, \tas)$ avec $M,N$ deux $A$-modules de type fini.

Pour les foncteurs coh\'erents on a une th\'eorie de la dualit\'e : si $V$ s'\'ecrit
comme conoyau d'un morphisme $\Hom_A(M, \tas) \fl \Hom_A(N, \tas)$ on d\'efinit le
foncteur dual $V^*$ comme le noyau du morphisme 
$N \T_A \tas \fl M\T_A \tas$. C'est un foncteur coh\'erent et l'op\'eration $^*$ est
involutive : on a
$V^{**} = V$.

Si  $V =
\bigoplus _{d \in \bZ} V_d : \Mod_A \fl \Gr_{R_A}$ est un foncteur tel que, pour
tout degr\'e $d$,  $V_d$ soit un foncteur coh\'erent, on d\'efinit  le
dual de $V$ par la formule 
$V^* = \bigoplus_{d \in \bZ} (V_{-d})^*$. Ce foncteur est \`a valeurs dans $\Gr_{R_A}$.

\th {Proposition 0.2}. \lign
1) Avec les notations pr\'ec\'edentes, si $V$ est exact \`a
droite (c'est-\`a-dire de la forme $V(A) \T_A \tas)$, on a un isomorphisme de
$R_A$-modules gradu\'es  : 
$V^* (A) 
\simeq (V(A))^*$.
\lign 
2) Si $\cF$ est un faisceau coh\'erent localement libre sur $\bP^r_A$ on a un
isomorphisme  de
$R_A$-modules gradu\'es  : 
$(H^r_* \cF)^* \simeq H^0_* \cF^\vee (-r-1)$. \lign
3) Soit $\cL$ un faisceau dissoci\'e sur $\bP^3_A$ et soit $L=H^0_* (\bP^3_A,
\cL)$. On a  un isomorphisme  de
$R_A$-modules gradu\'es  :  $ H^3_* (\bP^3_A, \cL) \simeq (L^\vee
(-4))^*$.

\dem 1) Posons $M = V(A)$. On a donc $V_n (A) = M_n$. Par d\'efinition on a $(V^*)_n =
(V_{-n})^*$, donc $(V^*)_n (A) = \Hom_A (M_{-n}, A)$, cf. [H] 4.1. On en d\'eduit
l'isomorphisme cherch\'e.

2) On applique le point 1) au foncteur $H^r_*(\cF \T_A \tas)$ qui est exact \`a droite
et la conclusion vient de [H] 7.6.

3) r\'esulte de 2) et de 0.1.2.

\tarte {d) Dualit\'e de Grothendieck et graduation}

Dans ce paragraphe nous donnons une variante gradu\'ee du th\'eor\`eme de dualit\'e de
Grothendieck. Pour tout ce qui concerne les cat\'egories d\'eriv\'ees on renvoie \`a
[RD].

Posons $X= \bP^r_A$ et $Y = \Spec A$ et soit $f : X \fl Y$ la projection canonique. 
Soit $\Mod_X$ la cat\'egorie des $\cO_X$-modules et $\cD(X)$ la cat\'egorie d\'eriv\'ee
de $ \Mod_X$. On note $\cD^+(X)$ (resp. $\cD^-(X)$) la sous-cat\'egorie correspondant
aux complexes born\'es inf\'erieurement (resp. sup\'erieurement) et
$\cD_{qc} (X)$ la sous-cat\'egorie des complexes \`a cohomologie quasi-coh\'erente. De
m\^eme on note $\cD (A), \cD^+(A),...$ et $\cD(R_A),...$ les cat\'egories d\'eriv\'ees de
$\Mod_A$ et $\Gr_{R_A}$.  

Comme $\Mod _X$ a
assez d'injectifs, le foncteur
$\G_*$ de
$\Mod_X$ dans
$\Gr_{R_A}$ admet un foncteur d\'eriv\'e 
$\bR \G_* : \cD^+(X) \fl \cD^+(R_A)$. Si $\cJ$ est un objet injectif de $\Mod_X$ il en
est de m\^eme de
$\cJ(n)$ et on en d\'eduit la formule
$(\bR\G_* \cF)_n =
\bR
\G \cF(n)$.

\vskip 0.3 cm

Soient $\cF$ et $\cG \in \Mod_X$. On pose 
$$\Hom_{X,*} (\cF, \cG) = \bigoplus_{n \in \bZ} \Hom_X(\cF, \cG (n))= \bigoplus_{n \in
\bZ} \Hom_X(\cF(-n), \cG).$$ On peut calculer le foncteur d\'eriv\'e de  $\Hom_{X,*}
(\tas\, ,\tas)$ en utilisant, par exemple, des
r\'esolutions injectives de la seconde variable et on obtient 
$$\bR \Hom_{X,*} (\tas\, ,\tas) : \cD^-(X) \times \cD^+(X) \fl \cD^+(R_A).$$
Ce foncteur est compatible avec la graduation et redonne donc le foncteur $\bR\Hom$
usuel en degr\'e $n$.

On consid\`ere aussi le foncteur d\'efini au paragraphe  {\it a)}
$$\Homgr_{A} (\tas\, , \tas) : \Gr_{R_A} \ti \Mod_A \fl  \Gr_{R_A}$$
et son d\'eriv\'e $\bR \Homgr_{A} (\tas\, , \tas)$. On a un morphisme de foncteurs 
$$\Hom_{X,*} (\tas \, , \tas) \fl \Homgr_{A} (\G_* (\tas) \, , \G (\tas))$$
et, au niveau des cat\'egories d\'eriv\'ees, on  a
le lemme suivant :

\th {Lemme 0.3}. Soient $\cF^\tas \in \cD^-(X)$ et $ \cG^\tas \in \cD^+(X)$. Il existe
un morphisme fonctoriel
$$\f : \bR \Hom_{X,*}(\cF^\tas, \cG^\tas) \fl \bR\Homgr_{A} (\bR \G_* (\cF^\tas), \bR
\G (\cG^\tas)).$$

\dem On notera que nous avons utilis\'e $\bR \G_*$ pour  $\cF^\tas$ mais $\bR \G$ pour
$\cG^\tas$. La preuve est analogue \`a  celle de [RD] II, 5.5, p. 103. On repr\'esente
$\cF^\tas$ par un complexe de faisceaux acycliques pour $\G_*$ (par exemple
des faisceaux flasques) et $\cG^\tas$ par un complexe de faisceaux injectifs. Comme
$\cG^\tas$ est injectif on a $\bR
\Hom_{X,*}(\cF^\tas,
\cG^\tas)  = \Hom_{X,*}(\cF^\tas, \cG^\tas) $ et ce dernier complexe s'envoie dans 
$\Homgr_{A} (\G_* \cF^\tas, \G \cG^\tas)$ qui n'est autre que 
$\Homgr_{A} (\bR \G_* (\cF^\tas), \bR
\G (\cG^\tas))$ car $\cF^\tas$ est acyclique et $\cG^\tas$ injectif. En composant
avec le morphisme canonique $\Homgr \fl \bR \Homgr$  on obtient $\f$.

\vskip 0.3 cm

On peut maintenant d\'efinir le morphisme de dualit\'e. Soit $\omega = \Omega^r_{X/Y}=
\cO_{\bP_A} (-r-1)$ le faisceau dualisant relatif. On pose, pour $G^\tas \in \cD^+(A)$,
$f\sh (G^\tas) = f^*(G^\tas) \T \omega [r]$, cf. [RD] p. 145. Pour $\cF^\tas \in
\cD^-(X)$ et $G^\tas \in \cD^+(A)$ on d\'efinit le morphisme de dualit\'e 
$$ \theta : \bR \Hom_{X,*} (\cF^\tas, f\sh G^\tas) \fl \bR \Homgr_{A} ( \bR \G_*
(\cF^\tas), G^\tas)$$
comme compos\'e du morphisme  $\f$ d\'efini en 0.3 et du morphisme trace  $t : \bR
\G f\sh (G^\tas) \fl G^\tas$, cf. [RD] III,  4.3, p. 155.
On a alors le th\'eor\`eme suivant :

\th {Th\'eor\`eme 0.4}. Pour tout $ \cF^\tas \in \cD^-_{qc} (X)$ et tout $G^\tas \in
\cD^+(A)$ le morphisme $\theta$ est un isomorphisme.

\dem 
Le th\'eor\`eme r\'esulte de [RD] III, 5.1, p.
161 appliqu\'e en chaque degr\'e.

\th {Corollaire 0.5}. Soit $\cF$ un faisceau localement libre de type fini sur
$\bP^3_A$. On a un isomorphisme fonctoriel dans $\cD^+(R_A)$ :
$$ \theta : \bR\G_* \cF^\vee (-4)\, [3] \simeq \bR \Homgr_{A} (\bR \G_* \cF, A).$$

\dem On prend $G^\tas = A$. On a alors $f\sh A = \omega\, [3]= \cO_{\bP}(-4) \, [3]$ par
d\'efinition de
$f\sh$. Par ailleurs, comme $\cF$ est localement libre de type fini on a $\Hom_X(\cF,
\cG) =
\G (\cF^\vee \T \cG)$, d'o\`u le r\'esultat.

\tarte {e) Courbes}

Sur un corps $k$ on appelle courbe un sous-sch\'ema de $\bP^3_k$ de dimension $1$, sans
composante ponctuelle (immerg\'ee ou non), c'est-\`a-dire localement de
Cohen-Macaulay.   Cette notion est stable par extension du corps de base, cf. [M],
remarque p. 182.

Si  $T$  est un sch\'ema, une famille de courbes   $\cC$  sur
$T$ (on dira  simplement une courbe de $\bP^3_T$), est un sous-sch\'ema
ferm\'e de
$\bP^3_T$, plat sur
$T$, et dont les fibres sont des courbes au sens pr\'ec\'edent. On note $\cO_{\cC}$
le faisceau structural de $\cC$, $\cJ_{\cC}$ le faisceau d'id\'eaux qui d\'efinit $\cC$
dans $\bP^3_T$  et $I_\cC$
son id\'eal satur\'e : $I_\cC = \bigoplus_{n \in \bN} H^0 \cJ_\cC(n)$.

\titre {1. La d\'efinition des triades} 
 
On reprend les notations du paragraphe 0 : $A$ d\'esigne  un anneau noeth\'erien et
$R_A$ l'anneau de polyn\^omes
 $A[X,Y,Z,T]$.
 
 \tarte {a) Complexes et foncteurs associ\'es}

Nous consid\'erons dans toute la suite des 
complexes
$\Lpont = (L_1 \Fl {d_1} L_0 \Fl {d_0} L_{-1})$  dans la cat\'egorie $\Gr_{R_A}$.  Le
groupe  d'homologie d'indice
$i$ d'un  tel complexe est not\'e
$h_i (\Lpont)$. C'est un
$R_A$-module gradu\'e.

\th {D\'efinition 1.1}. Soit $\Lpont = (L_1 \Fl {d_1} L_0 \Fl {d_0}
L_{-1})$  un complexe.
Les $R_A$-modules gradu\'es $N=h_1(\Lpont)$, $H=h_0(\Lpont)$, $C=h_{-1}(\Lpont)$ sont
appel\'es respectivement  {\bf noyau, c\oe ur} et {\bf conoyau} du complexe.  \lign
Le complexe  {\bf dual} $\Lpont^*$ est le complexe 
$L_{-1}^* \Fl {d_0^*} L_0^* \Fl {d_1^*} L_1^*$, cf. \S\ 0. Si $u : \Lpont \fl \Lpont'$
est un morphisme on note $u^*$ le morphisme induit sur les complexes duaux.\lign Le
complexe {\bf d\'ecal\'e}
$\Lpont (n)$ est le complexe d\'efini par les modules d\'ecal\'es, cf. \S\ 0 : 
$$L_1(n) \Fl {d_1(n)}
L_0(n)
\Fl {d_0(n)} L_{-1}(n).$$
Les complexes {\bf tronqu\'es} $\Lpontmr$ et $\Lpontps$ sont les complexes obtenus en
tronquant les mo\-du\-les $L_i$, cf. \S\ 0.

\vskip 0.2 cm

 On notera que la fl\`eche $d_0$ se factorise par $E= \Coker d_1$ en
une fl\`eche $d$. On a ainsi une suite exacte
$$0 \fl H \fl E \Fl{d} L_{-1} \fl C \fl 0,$$
qui fournit un
\'el\'ement de $\Ext^2_{R_A} (C,H)$.

\vskip 0.2 cm

\th {D\'efinition 1.2}.  Soit $\Lpont = (L_1 \Fl {d_1} L_0 \Fl {d_0} L_{-1})$ un
 complexe de
$R_A$-modules gradu\'es. On pose $L_i = \bigoplus_{n \in \bZ} L_{i,n}$. On dit que
$\Lpont$ est un {\bf complexe ad\'equat} si les $L_{i,n}$ sont des $A$-modules 
 plats et de type fini pour tout $i$ et tout $n$ (donc projectifs).

On notera que l'ensemble des complexes ad\'equats est stable par les op\'erations de
troncature, de somme directe et de dualit\'e et que la dualit\'e est involutive.

\th {D\'efinition 1.3}. Soit $\Lpont = (L_1 \Fl {d_1} L_0 \Fl {d_0}
L_{-1})$  un complexe ad\'equat.  Pour $i= 1,0,-1$ on consid\`ere le foncteur $h_i
(\Lpont
\T_A
\tas)$ de  
$\Mod_A$ dans $\Gr_{R_A}$ qui \`a un $A$-module $Q$ associe
$ h_i (\Lpont
\T_A Q)$. On
appelle {\bf foncteur associ\'e} \`a $\Lpont$ et on note $V_{\Lpont}$ 
 le foncteur $h_0 (\Lpont \T_A \tas)$. 

 Les foncteurs  $h_1
(\Lpont
\T_A
\tas)$ et   $h_{-1} (\Lpont \T_A
\tas)$ d'un complexe ad\'equat s'\'echangent par dualit\'e au
sens du paragraphe 0 {\it c)}, cf. [H] 4.3. Pr\'ecis\'ement $h_1 ((\Lpont^*)_n \T \tas)$
est dual de $h_{-1}( (\Lpont)_{-n} \T_A \tas)$. De m\^eme  on a 
$V_{(\Lpont^*)} = (V_{\Lpont})^*$ avec renversement de la graduation.

 \vskip 0.2 cm

Le lemme suivant est imm\'ediat :

\th {Lemme 1.4}. Soit $\Lpont = (L_1 \Fl {d_1} L_0 \Fl {d_0}
L_{-1})$  un complexe ad\'equat, $E= \Coker d_1$, $d : E\fl L_{-1}$ la fl\`eche d\'eduite
de
$d_1$ (cf. ci-dessus), $V$ le foncteur associ\'e. Alors, pour tout
$A$-module $Q$, on a un isomorphisme fonctoriel en
$Q$ :
$\Ker (d\T Q) \simeq V(Q)$.

\th {Lemme 1.5}. Soient $\Lpont = (L_1 \Fl {d_1} L_0 \Fl {d_0}
L_{-1})$  un complexe ad\'equat, $H$ et $C$ le c\oe ur et le conoyau de $\Lpont$, $V$ le
foncteur associ\'e.  Soit $Q$ un $A$-module
quelconque. On a une suite exacte fonctorielle
en $Q$ :
$$ \Tor_2^A(C,Q) \fl H\T_A Q \fl V(Q) \fl \Tor^A_1(C,Q) \fl 0.$$

\dem (cf. aussi [H] 3.6) On pose $K = \Im d_0$ et $E = \Coker d_1$. On a les suites exactes $0 \fl
H
\fl E
\fl K
\fl 0$ et
$0 \fl K \fl L_{-1} \fl C \fl 0$. En tensorisant par $Q$   on obtient les suites exactes
$$ \Tor^A_1(K,Q) \Fl {u} H \T_A Q \fl E \T_A Q \fl K \T_A Q \fl 0    $$
et, comme $L_{-1}$ est plat, 
$$ 0 \fl \Tor^A_1(C,Q) \fl K \T_A Q \fl L_{-1} \T_A Q \fl C \T_A Q \fl 0.    $$
 On pose $M_Q= \Coker u$ 
 et on applique le lemme du serpent au diagramme suivant
$$ \matrix {0 & \fl & M_Q & \fl & E \T_A Q & \fl & K \T_A Q & \fl & 0 \cr
&& \vf && \parallel && \vf \cr
0 & \fl & V(Q) & \fl & E \T_A Q & \Fl {d\T Q} & L_{-1} \T_A Q  \cr
}$$
(en vertu de 1.4 le module $V(Q)$ est le noyau de $d \T Q$). Le r\'esultat
s'ensuit ais\'ement en tenant compte de l'\'egalit\'e $\Tor^A_1(K,Q) =
\Tor^A_2 (C,Q)$.

\tarte {b) Pseudo-isomorphismes}

On rappelle qu'un quasi-isomorphisme de complexes (un \qis) est un morphisme qui induit
un isomorphisme sur les groupes d'homologie.  Le r\'esultat suivant pr\'ecise l'effet
d'un \qis\ sur les foncteurs.  Sa preuve est identique \`a celle de la derni\`ere
assertion de [AG] III 12.3 :

\th {Proposition 1.6}. Soit $u : \Lpont \fl \Lpont'$ un morphisme de complexes. On
suppose que $u$ est un \qis\ et que les modules $L_i$ et $L'_i$ sont plats sur $A$.
Alors, $u$ induit un isomorphisme de foncteurs $h_i( \Lpont \T_A \tas) \fl  h_i( \Lpont'
\T_A \tas)$ pour tout $i$.

La notion  de pseudo-isomorphisme, qui va jouer jouer un r\^ole essentiel dans ce qui
suit, appara\^\i t alors comme  une forme affaiblie de celle de 
\qis\ :

\th {D\'efinition 1.7}. Soit $u : \Lpont \fl \Lpont'$ un morphisme de complexes
ad\'equats. On dit que $u$ est un {\bf pseudo-isomorphisme} (en abr\'eg\'e un \ps) s'il
induit  un isomorphisme de foncteurs :
$h_0( \Lpont \T_A \tas) \fl  h_0( \Lpont' \T_A \tas)$ et  un monomorphisme de foncteurs
: $h_{-1}( \Lpont \T_A \tas)
\fl  h_{-1}( \Lpont' \T_A \tas)$.

Le r\'esultat suivant est imm\'ediat :

\th {Proposition 1.8}. \lign
a) Le compos\'e de deux \ps\  est un \ps,Ê\lign
b) si $g \circ f$ et $g$  sont des \ps, il en est de m\^eme de $f$,
\lign c) un morphisme homotope \`a un \ps\ est un \ps, \lign
d) si $u$ est un \qis, $u$ est 
un
\ps.

\th {Proposition 1.9}. Soit $\Lpont$ un complexe ad\'equat.\lign   On suppose que les foncteurs 
$h_0( \Lpont \T_A \tas)$ et $h_{-1}( \Lpont \T_A \tas) $ sont born\'es sup\'erieurement,
i.e. qu'on a $h_i( L_n \T_A \tas)=0 $ pour $i= 0,-1$ et $n > r$. Alors, la projection
$\Lpont \fl \Lpontmr$ est un \ps\  qui induit un isomorphisme sur les conoyaux et un 
\'epimorphisme sur les foncteurs $h_1$.

\dem  La suite exacte de complexes
$ 0 \fl \Lpontpr  \fl \Lpont  \Fl{p} \Lpontmr \fl 0 $ est scind\'ee comme suite de
$A$-modules, donc reste exacte quand on tensorise par un $A$-module $Q$. 
Comme l'homologie de $\Lpontpr\T_A \tas $ en degr\'es $0$ et $-1$  est nulle,
la suite  d'homologie associ\'ee donne le r\'esultat.

\tarte {c) Triades}

\th {D\'efinition 1.10}. Soit $\Lpont = (L_1 \Fl {d_1} L_0 \Fl {d_0} L_{-1})$ un
 complexe ad\'equat. \lign
On dit que $\Lpont$ est une
{\bf triade}  si  le noyau
$N=h_1(\Lpont)$ est un
$R_A$-module de type fini et si
le c\oe ur  $H=h_0(\Lpont)$ et le conoyau
$C=h_{-1}(\Lpont)$ sont des
$A$-modules de type fini.

\th {D\'efinition 1.11}.
Un morphisme (resp. un pseudo-isomorphisme) de triades  n'est rien d'autre qu'un
morphisme (resp. un pseudo-isomorphisme)
de complexes. 
On dit que deux triades  $\Lpont$ et $\Lpont '$ sont
\'equivalentes pour la relation de pseudo-isomorphisme (ou
simplement sont pseudo-isomorphes) s'il existe une cha\^\i ne de
\ps\ 
 qui les joint, dans laquelle les $\Lpont^{(i)}$ sont des triades :
$$\matrix {&&\Lpont^{(1)}&&&&\cdots&&&&\Lpont^{(n)} \cr
&\swarrow&& \searrow&&\swarrow&&\searrow&& \swarrow&& \searrow \cr
\Lpont&&&&\Lpont^{(2)} &&&&\Lpont^{(n-1)} & &&&\Lpont' \cr
}$$

\remas{1.12} \lign
1) Nous verrons au paragraphe {\it j)} comment on peut donner une d\'efinition plus
g\'en\'erale des triades et des \ps\ dans le cas de complexes non n\'ecessairement
ad\'equats (i.e.  sans hypoth\`eses de platitude).  Nous n'avons pas souhait\'e
introduire ce raffinement au d\'ebut de ce texte afin de ne pas trop l'alourdir. 
Toutes les triades envisag\'ees dans
cet article, \`a l'exception de celles du \S \ {\it j)}, seront des complexes
ad\'equats. \lign
2) L'existence d'un
 pseudo-isomorphisme entre deux triades implique que les les foncteurs associ\'es
sont isomorphes, mais la r\'eciproque est inexacte comme on le verra plus loin (cf.
1.35 c). \lign
3) Si $\Lpont$ est une triade les foncteurs $h_{-1}(\Lpont \T_A \tas)$ et 
$h_{0}(\Lpont \T_A \tas)$ sont born\'es (i.e., leurs composantes de degr\'es $\ll 0$ et
$\gg 0$ sont nulles). Pour $h_{-1}$ cela r\'esulte de la formule
$h_{-1}(\Lpont \T_A \tas)= C \T_A \tas$ et pour $h_0$ de 1.5. Le foncteur
$h_1$ est  born\'e inf\'erieurement. En effet, pour $n \ll 0$ on a
$N_n= H_n=C_n=0$, de sorte que la suite $0 \fl L_{1,n} \fl L_{0,n} \fl L_{-1,n} \fl 0$
est exacte et, comme $L_{-1,n}$ est plat sur $A$, elle le reste par tensorisation par un
$A$-module $Q$.

\tarte {d) Triades majeures}

\th {D\'efinition 1.13}. On appelle {\bf triade majeure} une triade $\Lpont= (L_1 \Fl
{d_1} L_0 \Fl {d_0} L_{-1})$ dans laquelle les $L_i$ sont des $R_A$-modules
dissoci\'es (donc libres si $A$ est local).

\vskip 0.3 cm

Soit $\Lpont = (L_1 \Fl {d_1} L_0 \Fl {d_0}
L_{-1})$  un complexe. Rappelons qu'une {r\'esolution libre} de
$\Lpont$ est un complexe $\Fpont = (F_i)_{i \geq -1}$   de $R_A$-modules gradu\'es
libres, avec un \qis\ $u:\Fpont \fl \Lpont$. 
Il est clair qu'une telle r\'esolution libre existe toujours (mais le complexe $\Fpont$
a, en g\'en\'eral, plus que trois termes).

\th {Proposition 1.14}. Toute triade $\Lpont$ admet une {\bf  r\'esolution majeure} :
      il existe une triade majeure $\Mpont$ et un \ps\  $u : \Mpont \fl
\Lpont$ qui induit, de plus, un isomorphisme sur 
les conoyaux et un \'epimorphisme sur les foncteurs $h_1$.

\dem En vertu de [AG] III 12.3 il existe une r\'esolution  $\Fpont$ de
$\Lpont$ avec les $F_i$ libres de type fini sur $R_A$. On obtient $\Mpont$ en ne gardant
que les termes de degr\'es
$-1, 0, 1$~:
$\Mpont = (F_1
\fl F_0
\fl F_{-1})$.

\rema {1.15} Attention, l'emploi du mot r\'esolution en 1.14 est un abus de langage~:
une r\'esolution majeure n'est pas une r\'esolution projective du complexe au sens
usuel  car la fl\`eche $u$ n'est pas un
\qis\ (elle n'est que surjective sur les noyaux). Cependant on passe ais\'ement d'une
notion 
\`a l'autre. Si on a une r\'esolution libre de type fini on obtient une
r\'esolution majeure en tronquant comme  ci-dessus. Inversement, si on a une
r\'esolution majeure on la prolonge en une r\'esolution projective par le processus
suivant  : on a
$u :
\Mpont
\fl
\Lpont$ qui induit un morphisme surjectif $\wh u : N' \fl N$ sur les noyaux. Soit $Q$
le noyau de
$\wh u$ et soit  $\cdots \fl F_3 \fl F_2 \fl Q \fl 0$ une r\'esolution libre de $Q$.
Alors le complexe $\Fpont =(\cdots \fl F_3 \fl F_2 \fl M_1 \fl M_0 \fl M_{-1})$ donne
une r\'esolution projective de $\Lpont$ avec un diagramme commutatif:
$$\matrix{ \Mpont & \fl & \Fpont \cr &\seFl{u}& \vf \cr &&\Lpont \cr
}$$ On dit que $\Fpont' = (\cdots
\fl F_3
\fl F_2)
$ est une {queue} pour la triade $\Mpont$.

\tarte {e) \'Etude des pseudo-isomorphismes de triades}

Nous donnons d'abord une caract\'erisation  des \ps\ :

\th {Proposition 1.16}. Un morphisme de triades $f : \Lpont \fl \Lpont '$ est un \ps\
si et seulement s'il v\'erifie les propri\'et\'es suivantes :  \lign
a) $f$ induit un isomorphisme $  H \simeq H'$ sur les c\oe urs, \lign
b) $f$ induit une injection $C \fl C'$ des conoyaux et le quotient $D = C'/C$ est
plat sur $A$. \lign
Si de plus $f$ induit une surjection $N \fl N'$ sur les noyaux, $f$ est un \ps\  
qui induit un \'epimorphisme sur les foncteurs $h_1$.

\dem
Supposons d'abord que $f$ soit un \ps.  On en
d\'eduit  que
$H= V(A)$ et
$H'= V'(A)$ sont isomorphes et que $f$ induit une injection $C \fl C'$ dont on note $D$
le conoyau. Comme
$V$ et
$V'$ sont isomorphes, il r\'esulte de 1.5 que
$\Tor_1^A (C, \tas) \fl \Tor_1^A (C', \tas)$ est un isomorphisme et, avec la suite
exacte 
$$ \Tor_1^A (C, \tas)  \fl \Tor_1^A (C', \tas)  \fl \Tor_1^A (D, \tas)  \fl C \T \tas
\fl C' \T \tas
$$
 et le monomorphisme de foncteurs   $ C \T \tas
\fl C' \T \tas$,
on en d\'eduit
 $\Tor_1^A (D, \tas) =0$, donc $D$ plat sur $A$.  

Inversement, si on a  la suite exacte $0 \fl C \fl C' \fl D \fl 0$
avec $D$ plat sur $A$ le morphisme $C \T \tas \fl C'\T \tas$ est un
monomorphisme et   $\Tor_i^A (C, \tas)  \fl \Tor_i^A (C', \tas) $ est un
isomorphisme pour $i>0$. Si, de plus, on a
$H \simeq H'$ on conclut, avec 1.5, que $V$ et $V'$ sont isomorphes, donc que $f$ est
un \ps.

Si de plus $f$ induit une surjection $N \fl N'$, montrons  que pour un
$A$-module 
$Q$ le morphisme
 $h_1(\Lpont\T_A Q) \fl h_1(\Lpont'\T_A Q)$ est surjectif. C'est vrai pour $Q=A$ par
hypoth\`ese, donc aussi pour $Q$ libre  sur $A$. Dans le
cas g\'en\'eral on prend une r\'esolution $0 \fl G \fl F \fl Q \fl 0$ avec $F$ libre.
Comme les $L_i$ sont plats on en d\'eduit une suite exacte de complexes :
$$0 \fl\Lpont \T G \fl \Lpont \T F \fl \Lpont \T Q \fl 0$$
et de m\^eme pour $\Lpont'$. En d\'eroulant l'homologie on a le diagramme commutatif
de suites exactes :
$$ \matrix {h_1(\Lpont \T F)& \fl& h_1(\Lpont \T Q)& \fl& h_0 (\Lpont \T G) & \fl &
h_0(\Lpont \T F) \cr
\Vf{\a}&& \Vf{\beta}&& \Vf{\g}&& \Vf {\de}\cr 
h_1(\Lpont' \T F)& \fl& h_1(\Lpont' \T Q)& \fl& h_0 (\Lpont' \T G) & \fl &
h_0(\Lpont' \T F) \cr}$$
Le morphisme $\a$ est surjectif car $F$ est libre et $\g$ et $\de$ sont des
isomorphismes car $u$ est un \ps. Il en r\'esulte que $\beta$ est surjectif.

\vskip 0.3 cm

\tarte {f) Un ``lemme de Verdier''}

Les propositions qui suivent ont pour but d'\'etablir un lemme de calcul de
fractions  pour la relation de pseudo-isomorphisme  1.11  du type de celui de
Verdier (cf. [V]) qui permet de se ramener \`a des
cha\^\i nes de \ps\ de longueur $2$.

\th {Proposition 1.17}. Soient $\Lpont, \Lpont', \Lpont''$ des triades et $u : \Lpont
\fl \Lpont'$ et $v : \Lpont
\fl \Lpont''$ des morphismes. On suppose que $u$ induit un isomorphisme $H \simeq H'$
sur les c\oe urs et une injection $C \fl C'$ sur les conoyaux. On pose $D = C'/C$.
Alors, il existe une triade $\Lpont'''$ et des morphismes $u' : \Lpont''
\fl \Lpont'''$ et $v' : \Lpont'
\fl \Lpont'''$ avec $u'v= v'u$ \`a homotopie pr\`es et tels que $u'$ induise un
isomorphisme $H'' \simeq H'''$ sur les c\oe urs, une injection $C'' \fl C'''$ sur les
conoyaux, avec le m\^eme quotient $C'''/C'' \simeq D$. On a donc le diagramme suivant :
$$\matrix {\Lpont & \Fl {v}& \Lpont'' \cr
\Vf{u}&&\Vf{u'} \cr
\Lpont'&\Fl {v'}&\Lpont''' \cr
}$$

\dem Quitte \`a rajouter des queues (form\'ees de $R_A$-modules de type fini) comme en
1.15 on peut supposer que 
$\Lpont, \Lpont', \Lpont'' $ sont des complexes d\'efinis pour $i \geq -1$ avec comme
seuls modules d'homologie non nuls les c\oe urs et les conoyaux. Comme la cat\'egorie
des complexes  est triangul\'ee ([RD] I 2) on peut  ins\'erer
le morphisme
$u$ dans un triangle $\Lpont \Fl{u} \Lpont' \fl \Xpont \fl \Lpont[-1]$ o\`u le
complexe $\Xpont$ est le c\^one de $u$, donc est ad\'equat. On obtient une suite
exacte d'homologie
$$ \cdots h_0 \Lpont \fl h_0 \Lpont' \fl h_0 \Xpont \fl h_{-1} \Lpont \fl h_{-1}
\Lpont'  \fl h_{-1} \Xpont \fl h_{-2} \Lpont
$$
dont on d\'eduit qu'on a $h_{-1} \Xpont= C'/C=D$ et $h_i \Xpont = 0$ pour $i \neq -1$.

On consid\`ere ensuite le morphisme compos\'e $\Xpont \fl \Lpont[-1] \fl \Lpont''[-1]$
que l'on ins\`ere dans un triangle $\Lpont'' \fl \Lpont''' \fl \Xpont \fl
\Lpont''[-1]$ avec $\Lpont'''$ ad\'equat. La suite exacte longue analogue pour ce
triangle montre que le noyau de $\Lpont'''$ est nul et que
$u'$ induit un isomorphisme des c\oe urs
$H'' \simeq H'''$ et une injection des conoyaux $C'' \fl C'''$ avec pour quotient $D$.

Enfin, l'existence du morphisme $v'$ commutant, \`a homotopie pr\`es, avec les autres
r\'esulte de l'axiome (TR3) des cat\'egories triangul\'ees ([RD] p. 21). Il ne reste
plus qu'\`a couper les queues des complexes obtenus pour retrouver les triades
cherch\'ees (la finitude des noyaux r\'esulte du fait que les queues sont form\'ees de
$R_A$-modules de type fini).

\th {Corollaire 1.18}. On reprend les notations de 1.17 et on suppose que $u$ est un
\ps. Alors $u'$ est un \ps. De plus, si $v$ est aussi un \ps, il en est de m\^eme de
 $v'$.

\dem Si $u$ est un \ps\ le module $D$ est plat sur $A$ (cf. 1.16), donc $u'$ est aussi
un \ps.  Si $v$ est un \ps, le compos\'e $u'v$ aussi et donc aussi $v'u$ qui lui est
homotope (cf. 1.8). On v\'erifie enfin que $v'$ induit un
isomorphisme $H' \simeq H'''$ et une injection $C' \fl C'''$ avec $C'''/C' \simeq
C''/C$ donc plat sur $A$ puisque $v$ est un \ps.

\th  {Corollaire 1.19}. (Lemme de Verdier pour les triades) Soient $\Lpont$ et
$\Lpont'$ deux triades pseudo-isomorphes. Alors il existe une triade $\Lpont''$ et des
\ps\ $\Lpont \fl \Lpont''$ et $\Lpont' \fl \Lpont''$.

\dem Par r\'ecurrence sur la longueur de la cha\^\i ne de \ps\ joignant $\Lpont$ et
$\Lpont'$ en utilisant 1.18.

\vskip 0.2 cm

Attention, on n'a pas, en g\'en\'eral, l'analogue du lemme de Verdier avec les
fl\`eches en sens inverse. Voir cependant 2.{\it d)} pour un succ\'edan\'e.

\th {Proposition 1.20}. (Lemme de factorisation) Soient $u:\Lpont \fl \Lpont''$ et
$v:\Lpont' \fl \Lpont''$ deux morphismes de triades. On suppose : \lign
a) que $v$ induit un isomorphisme $H'Ê\simeq H''$ sur les c\oe urs et une injection
$C' \fl C''$ sur les conoyaux, \lign
b) que le morphisme $C \fl C''$ induit par $u$ se factorise par $C'$, \lign
c) que $\Lpont$ est une triade majeure. \lign
Alors il existe un morphisme $w : \Lpont \fl \Lpont'$ tel que l'on ait  $u= vw$ \`a
homotopie pr\`es.

\dem Comme  dans 1.15 on peut supposer, en ajoutant des queues, que les triades sont des
complexes d\'efinis pour $i \geq -1$, avec seulement deux groupes d'homologie non
nuls. 

On compl\`ete le morphisme $v$ en un triangle $\Lpont' \fl \Lpont'' \fl \Xpont \fl
\Lpont'[-1]
$. L'hypoth\`ese a) montre que l'on a $h_{-1} \Xpont  = C''/C' = D$ et que les autres
groupes d'homologie de $\Xpont$ sont nuls. On compl\`ete ensuite le morphisme
compos\'e $\Lpont \fl \Xpont$ en un triangle $ \Lpont''' \fl \Lpont \fl \Xpont \fl
\Lpont'''[-1]$. La suite exacte d'homologie donne alors $H''' \simeq H$ et une suite
exacte~:
$$0 \fl C''' \fl C \fl D \fl h_{-2} (\Lpont'''),$$
et, puisque le morphisme $C \fl C''$ se factorise par $C'$ on voit que $C \fl D$ est
nul, donc $C''' \simeq C$.

En vertu de l'axiome (TR3) on a un morphisme $\Lpont''' \fl \Lpont'$ qui rend
commutatif, \`a homotopie pr\`es, le diagramme suivant :
$$ \matrix {  \Lpont'''& \fl& \Lpont &\fl &\Xpont &\fl&
\Lpont'''[-1] \cr
\vf&&\vf&&\parallel  &&\vf \cr
\Lpont' &\fl& \Lpont'' &\fl &\Xpont& \fl&
\Lpont'[-1]\cr
}$$
On notera que le groupe $h_{-2} \Lpont'''$ peut \^etre non nul. On remplace donc
$\Lpont'''$ par $\Lpont^{iv} = \cdots L_1''' \fl L_0''' \fl Z_{-1} (\Lpont''') \fl
0$ (o\`u $Z_{-1} (\Lpont''')$ est le noyau de $L_{-1}''' \fl L_{-2}'''$) et on a un
morphisme canonique
$\Lpont^{iv}
\fl
\Lpont'''$ qui induit un isomorphisme sur les $h_i$ pour $i \geq -1$. Alors, le
morphisme
$\Lpont^{iv} \fl
\Lpont$ est un \qis. Comme $\Lpont$ est une triade majeure, il r\'esulte de [RD] I, 4.5,
p. 41 (appliqu\'e en renversant les fl\`eches et avec des libres au lieu d'injectifs)
que ce morphisme a un inverse
\`a homotopie pr\`es. Le morphisme cherch\'e s'obtient en composant $\Lpont
\fl \Lpont^{iv} \fl \Lpont''' \fl \Lpont'$.

\vskip 0.3 cm

Avec cette proposition et  1.8  on obtient :

\th {Corollaire 1.21}. Soient $u:\Lpont \fl \Lpont''$ et
$v:\Lpont' \fl \Lpont''$ deux \ps, avec $\Lpont$ majeure. On suppose que le 
morphisme $C
\fl C''$ induit par $u$ se factorise par $C'$. Alors il existe un \ps\  $w : \Lpont \fl
\Lpont'$  v\'erifiant $u= vw$ \`a
homotopie pr\`es.

\th {Corollaire 1.22}. (Lemme de Verdier pour les triades majeures) Si $\Lpont$ et
$\Lpont'$ sont deux triades majeures pseudo-isomorphes il existe une triade majeure
$\Lpont''$ et des \ps\ $\Lpont \fl \Lpont''$ et $\Lpont' \fl \Lpont''$.

\dem En vertu du lemme de Verdier on a $\Lpont''$ qui minore les deux triades. On
regarde alors la r\'esolution majeure
$\Mpont$ de $\Lpont''$ et on conclut par 1.21.

\tarte {g) Triades duales}

Soit $\Lpont$ une triade. Le complexe dual $\Lpont^*$   n'est pas,
en g\'en\'eral, une triade. En effet, si le c\oe ur et le noyau de $\Lpont^*$ sont des
$A$-modules de type fini,  le conoyau n'est pas, en g\'en\'eral, de type fini sur $A$
(consid\'erer par exemple, la triade $\Lpont =(R_A \fl 0 \fl 0)$). Nous allons
contourner cette difficult\'e en utilisant une troncature,
mais la triade duale sera d\'efinie seulement \`a \ps\ pr\`es.

Pour \'etablir les propri\'et\'es essentielles de la dualit\'e nous aurons besoin de
deux notions suppl\'ementaires  et de deux lemmes concernant ces notions :

\th {D\'efinition 1.23}.  \lign
1) On appelle {\bf triade mineure} une triade
$\Ppont$ dans laquelle les $R_A$-modules
$P_i$ sont   de type fini sur $A$. \lign
2) Un {\bf pseudo-isomorphisme fort} de complexes ad\'equats (en abr\'eg\'e un \ps\ fort) est un
\ps\ qui induit un \'epimorphisme de foncteurs : $h_1( \Lpont \T_A \tas)
\fl  h_1( \Lpont' \T_A \tas)$.

\th {Lemme 1.24}. \lign
1) Soit $\Lpont$ une triade majeure, $p :
\Lpont \fl \Lpontmr$ la projection canonique. Alors, $\Lpontmr$ est une triade mineure
et, si $r$ est assez grand, $p$ est un \ps\  fort qui induit un isomorphisme sur les
conoyaux. Cette op\'eration de troncature
est fonctorielle.
\lign  2) Toute triade est pseudo-isomorphe \`a une triade mineure,
pr\'ecis\'e\-ment, si $\Lpont$ est une triade, il existe deux triades $\Mpont$ et
$\Ppont$ avec $\Mpont$ majeure et $\Ppont$ mineure et des \ps\ forts $u : \Mpont \fl
\Lpont$ et $v : \Mpont \fl \Ppont$.

\dem

 1) Comme $\Lpont$ est une triade majeure les modules $L_{i,n}$ sont nuls pour $n\ll 0$, donc
les modules tronqu\'es $L_{i, \leq r}$ sont de type fini sur $A$ et $\Lpontmr$ est bien une
triade mineure. Comme
$C$ et
$H$ sont de type fini sur
$A$  on a
$C_n=H_n=0$ pour $n \gg 0$, disons pour
$n > r$ et on en d\'eduit que les foncteurs $h_0(\Lpontpr \T_A \tas)$ et
$h_{-1}(\Lpontpr \T_A \tas)$ sont nuls. On conclut alors par 1.9.
Le point 2) r\'esulte de 1.14 et de 1).

\th {Lemme 1.25}. \lign
1) Le dual d'un \ps\ fort est un \ps\ fort. \lign
 2) Soit $u : \Lpont \fl \Lpont'$ un \ps\
entre deux triades. Il existe une triade $\Lpont''$ et des \ps\ forts $v : \Lpont'' \fl
\Lpont$ et
$w : \Lpont'' \fl \Lpont'$ (mais le triangle n'est pas, en g\'en\'eral, commutatif). 
\lign
 Sur l'ensemble des triades les relations d'\'equivalence engendr\'ees par les \ps\ et les \ps\
forts  sont les m\^emes.

\dem Le point 1) est clair. Pour 2), soit $N'$  le noyau de $\Lpont'$ et $\pi : F \fl N'$ une
r\'esolution libre de
$N'$. On pose 
$L''_{-1} = L_{-1}$,
$L''_0 = L_0$, $d''_0= d_0$,
$L''_1 = L_1 \oplus F$ et $d''_1 = (d_1,0)$. On
pose ensuite
$v_{-1} =
\Id$,
$v_0=
\Id$, $v_1= (\Id,0)$ et $w_{-1} = u_{-1}$, $w_0= u_0$ et $w_1= (u_1,i \pi)$ o\`u $i$ est
l'injection de $N'$ dans $L'_1$ et on a la propri\'et\'e requise.

\th {Proposition-d\'efinition 1.26}. \lign
Soit $\Lpont$ une triade, et soit $\Lpont^*$ le complexe dual. \lign
1) Il existe un entier  $r$ 
tel que, pour tout $n \leq r$, les foncteurs $h_1((\Lpont^*)_n \T_A \tas)$ et
$h_0((\Lpont^*)_n \T_A \tas)$ soient nuls. 
 Un tel entier sera dit convenable. \lign
2) Si $r$ est convenable le complexe $(\Lpont^*)_{>r}$ est une triade et l'inclusion
$j:(\Lpont^*)_{>r} \subset \Lpont^*$ est un \ps\ fort. On dit que $(\Lpont^*)_{>r}$ est {\bf une
triade duale} de
$\Lpont$.
\lign Si $r$ et $r'$ sont deux entiers convenables avec $r'<r$ l'inclusion $j_{r,r'} :
(\Lpont^*)_{>r} \subset (\Lpont^*)_{>r'}$ est un \ps\ fort. Si on a trois entiers convenables avec
$r''<r'<r$ on a $j_{r,r'} j_{r',r''} = j_{r,r''}$.

\dem Le point 1) est cons\'equence du fait que les foncteurs $h_1((\Lpont^*)_n \T_A
\tas)$ et
$h_0((\Lpont^*)_n \T_A \tas)$ sont respectivement  duaux de 
$h_{-1}((\Lpont)_{-n}
\T_A
\tas)$ et
$h_{0}((\Lpont)_{-n} \T_A \tas)$ qui sont nuls pour $n \ll 0$ et $n \gg 0$ (cf. 1.12.3). 

Comme 
$h_{1}(\Lpont \T_A \tas)$ est born\'e inf\'erieurement, cf. 1.12.3, 
$h_{-1}(\Lpont^* \T_A \tas)$ est born\'e sup\'erieurement et il est clair alors
que
$(\Lpont^*)_{>r}$ est une triade.  On consid\`ere la suite exacte de complexes 
$0 \fl (\Lpont^*)_{>r}  \Fl{j} \Lpont^* \fl (\Lpont^*)_{\leq r} \fl 0$. La suite d'homologie
associ\'ee jointe aux hypoth\`eses de nullit\'e sur $(\Lpont^*)_{\leq r}$ montre  que $j$
est un
\ps\ fort. La derni\`ere propri\'et\'e est claire.

\vskip 0.3 cm

La proposition suivante r\'esume les propri\'et\'es de la dualit\'e :

\th {Proposition 1.27}. \lign
0) Si $\Ppont$ est une triade mineure le complexe $\Ppont^*$ (qui est aussi une triade mineure)
est une triade duale de $\Ppont$.
\lign 1) Soient $\Mpont$ et $\Mpont'$ deux triades duales de $\Lpont$ et
$\Lpont'$ respectivement. Si $\Lpont$ et $\Lpont'$ sont pseudo-isomorphes il en est de
m\^eme de
$\Mpont$ et $\Mpont'$.
\lign 
2) Si $\Mpont$ est une triade duale de $\Lpont$, $\Lpont$ est pseudo-isomorphe \`a toute
triade duale de
$\Mpont$. \lign
3) Le
foncteur associ\'e \`a une duale de $\Lpont$ est le dual (cf. \S\ 0) du foncteur 
associ\'e
\`a
$\Lpont$. \lign 
4) Si $\Lpont$ est une triade, $\Mpont$ une triade majeure et si on a un morphisme de
complexes $u : \Mpont \fl \Lpont^*$ qui est un \ps, $\Mpont$ est
pseudo-isomorphe \`a toute triade duale de $\Lpont$.

\dem 
Pour le point 0), si $\Ppont$ est une triade mineure il en
est de m\^eme de $\Ppont^*$  donc, pour $r$ assez petit, on a $\Ppont^*=
(\Ppont^*)_{>r}$. On conclut  par 1.26. 

Montrons l'assertion  1). Si $\Lpont$
et
$\Lpont'$ sont pseudo-isomorphes  on peut, en vertu de 1.25.2, les joindre par une
cha\^\i ne de 
\ps\ forts entre des triades $(\Lpont)_i$. On a alors, par 1.25.1,  des \ps\ forts entre les
complexes duaux
$ ((\Lpont)_i)^*$ qui induisent, pour $r$ assez petit,  des \ps\ entre les complexes
tronqu\'es. Il en r\'esulte que les complexes tronqu\'es associ\'es \`a $\Lpont^*$ et 
$\Lpont'^*$ (c'est-\`a-dire les triades duales) sont pseudo-isomorphes.

 Soit $\Ppont$ une triade mineure pseudo-isomorphe \`a $\Lpont$ (cf. 1.24.2). Les points 0) et
1) ci-dessus montrent successivement que $\Mpont$ est pseudo-isomorphe \`a
$\Ppont^*$,  puis que toute triade duale de $\Mpont$ est pseudo-isomorphe \`a
$(\Ppont^*)^*= \Ppont$, donc aussi \`a $\Lpont$, ce qui \'etablit le point 2). 

 L'assertion sur les foncteurs
r\'esulte de [H] 4.3 (cf. 1.3) et du fait que l'inclusion $(\Lpont^*)_{>r} \subset \Lpont^*$ est un
\ps.

Montrons le point 4). Pour $r$ assez petit $(\Lpont^*)_{>r}$ est une triade duale de $\Lpont$
et, comme
$\Mpont$ est majeure,  
$u$  se factorise  en 
$\Mpont \Fl {v} (\Lpont^*)_{>r} \Fl {j} \Lpont^*$. Comme $u$ et    $j$  sont des \ps,
 on conclut par 1.8.

\tarte {h) Triades minimales, triades \'el\'ementaires}

\th {Proposition-D\'efinition 1.28}. 
1) L'ensemble des $R_A$-modules gradu\'es $C$  qui sont de type fini sur $A$ est
(partiellement) ordonn\'e par la relation :
$ C \prec C'${s'il existe} $u : C \fl C'$, $R_A$-lin\'eaire,
injectif  avec $ C'/C $ {plat sur } $A$. \lign
2) L'ensemble des triades pseudo-isomorphes \`a une
triade donn\'ee est pr\'eordonn\'e par la relation : $\Lpont \prec \Lpont'$
s'il existe un \ps\ $ u : \Lpont \fl \Lpont'$. \lign 
Si on a $\Lpont \prec \Lpont'$ et $\Lpont' \prec \Lpont$, il existe un \ps\ 
$u :
\Lpont
\fl
\Lpont'$ qui induit
un isomorphisme sur les conoyaux. \lign Les
\'el\'ements minimaux pour cette relation sont appel\'es {\bf triades minimales}. Un
\'el\'ement minimum, s'il en existe, est appel\'e  {\bf triade \'el\'ementaire}.\lign
3) L'application qui
\`a une triade associe son conoyau est  croissante pour les relations
pr\'ec\'edentes.

\dem Dans le point 1) seule l'antisym\'etrie de la relation est non triviale. Si on a
$C \prec C' \prec C$ on a une suite exacte $0 \fl C \Fl {i} C \fl M \fl 0$ avec $M$
plat sur $A$. Comme $C$ est de type fini on en d\'eduit $M=0$.

Les autres points de la proposition sont clairs.

\th {Lemme 1.29}. Soit $C$ un $R_A$-module gradu\'e, de type fini sur $A$. Alors, il
existe $C_0$ avec $C_0 \prec C$ et
$C_0$ minimal pour la relation $\prec$.

\dem En effet, si on a $C_0 \prec C$, le quotient $C/C_0$ est plat sur $A$ et la
conclusion s'obtient en raisonnant par r\'ecurrence sur la somme des rangs de ce
quotient aux points g\'en\'eriques.

\vskip 0.3 cm

Le lemme suivant assure l'existence de triades minimales :

\th {Lemme 1.30}. Soit $\Lpont= (L_1 \Fl {d_1} L_0 \Fl {d_0}
L_{-1})$ une triade, $C$ son conoyau et $C'$ un
sous-$R_A$-module de $C$. Il existe une triade $\Lpont'$, de conoyau $C'$, et un
\ps\ $u : \Lpont'
\fl
\Lpont$ qui induit un isomorphisme sur les c\oe urs et l'injection donn\'ee sur les
conoyaux. \lign
Il existe une triade minimale parmi les triades
pseudo-isomorphes \`a $\Lpont$.

\dem On peut supposer $\Lpont$ majeure, quitte \`a la remplacer par une r\'esolution. Si
$p : L_{-1}
\fl C$ est le morphisme canonique on choisit une r\'esolution libre $L'_{-1} \fl  p^{-1}
(C')$. Le morphisme $d_0$ est \`a valeurs dans $p^{-1}
(C')$ et on peut le relever en $d'_0 : L_0 \fl L'_{-1}$. On prend alors
$\Lpont'= (L_1
\Fl {d_1} L_0
\Fl{d'_0} L'_{-1})$.

Le deuxi\`eme point r\'esulte du premier et de 1.29.

\rema {1.31} En revanche, il n'existe pas en g\'en\'eral de triade \'el\'ementaire dans
une classe d'\'equivalence pour la relation de pseudo-isomorphisme (cf. 1.35.a). Nous
allons voir que c'est vrai dans le cas d'un \avd.

\th {Proposition 1.32}. Soit $A$ un \avd. \lign
1) Une triade majeure $\Mpont$ est \'el\'ementaire si et seulement si son conoyau
$C$ est de torsion. \lign
2) Soit
$\Lpont$ une triade,
 $C$ son conoyau et $C_\tau$ le sous-$R_A$-module de $C$ form\'e des \'el\'ements de
torsion sur $A$. Alors il existe une triade majeure \'el\'ementaire $\Mpont$ de conoyau
$C_\tau$ et un
\ps\ $u:\Mpont \fl \Lpont$. De plus,
cette triade est unique \`a un \ps\ pr\`es induisant un isomorphisme sur les conoyaux.

\dem 

1) Soit $\Lpont'$ une triade pseudo-isomorphe \`a $\Mpont$. Alors, il existe une triade
$\Lpont''$ et des \ps\ $\Mpont \fl \Lpont''$ et $\Lpont' \fl \Lpont''$ en vertu de
1.19. On applique alors 1.20 \`a ces morphismes. Il faut montrer que l'inclusion $C
\fl C''$ se factorise par $C'
\fl C''$ ce qui r\'esulte du fait que le quotient  $C''/C'$
est plat et que $C$ est de torsion.

2) Le lemme 1.30 assure l'existence d'une triade $\Lpont'$ de conoyau $C_\tau$ et d'un
morphisme
$u :
\Lpont' \fl \Lpont$. Comme $C/C_\tau$ est plat sur $A$, $u$ est un \ps.  Pour avoir la
triade majeure $\Mpont$ cherch\'ee il suffit de prendre une r\'esolution majeure de
$\Lpont'$. L'unicit\'e vient du point 1) et de 1.28.2.

\th {Corollaire 1.33}. Soit $A$ un \avd\ et soit $u : \Lpont \fl \Lpont'$ un morphisme
de triades qui induit un isomorphisme sur les foncteurs associ\'es. Alors $\Lpont$ et
$\Lpont'$ sont pseudo-isomorphes. Plus pr\'ecis\'ement on a un diagramme, commutatif
\`a homotopie pr\`es, o\`u les morphismes autres que $u$ sont des \ps\ :
$$\matrix { \Mpont & \Fl{v} & \Mpont' \cr 
\vf&&\vf \cr
\Lpont&\Fl {u} & \Lpont' \cr
}$$

\dem
En vertu de 1.5,
$u$ induit un isomorphisme de foncteurs
$\Tor_1^A (C,
\tas)
\fl 
\Tor_1^A (C', \tas)$. Cela implique que les sous-modules de torsion $C_\tau$ et
$C'_\tau$ sont isomorphes (cf. [H] 5.2). Si $\Mpont$ et $\Mpont'$ sont les triades
majeures
\'el\'ementaires donn\'ees par 1.32, il existe un morphisme $v$ comme annonc\'e en
vertu de 1.20.

(Attention, $u$ lui-m\^eme n'est pas n\'ecessairement un \ps\  car il ne v\'erifie pas
l'hypoth\`ese sur les conoyaux.)

\remas {1.34} \lign 
1) Soient $\Lpont$ et $\Lpont'$ deux triades et $V_L$ et $V_{L'}$ les
foncteurs  associ\'es. Attention, le corollaire 1.33 n'implique pas que si $V_L$ et
$V_{L'}$ sont isomorphes les triades soient elles-m\^emes pseudo-isomorphes. L'existence
du morphisme $u$ est essentielle, cf. 1.35.c.\lign
2)  On montre ais\'ement, \`a la mani\`ere de  [Ho], qu'une triade majeure
\'el\'ementaire est ca\-rac\-t\'eris\'ee, \`a un \ps\ induisant un isomorphisme sur les
conoyaux pr\`es, par son c\oe ur, son conoyau et par l'\'el\'ement de $\Ext^2_{R_A}
(C,H)$ associ\'e.

\tarte {i) Exemples}

\expls {1.35}

\noindent a) Sur un anneau quelconque, une classe d'\'equivalence de triades pour la relation de
pseudo-isomorphisme ne contient pas n\'ecessairement d'\'el\'ement minimum
(c'est-\`a-dire de triade \'el\'emen\-taire).

Soit $A$ un anneau local r\'egulier de dimension $2$ d'id\'eal maximal $m$. Soit $C$
le $R_A$-module \'egal \`a $m\oplus A$ en degr\'e $0$ et \`a $A$ en degr\'e $1$ avec
la multiplication par $X : m\oplus A \fl A$ donn\'ee par la matrice $(0,1)$ et les
autres multiplications nulles. Soit $C_1$ le sous-$R_A$-module $m$ de $C$, de sorte que
$C/C_1= A \oplus A$ est plat sur $A$. Soit $\psi : m \fl m\oplus A$ le morphisme
donn\'e par
$\psi (a) = (a,a)$ et soit 
$C_2$ le sous-module de
$C$ \'egal \`a $\Im \psi$ en degr\'e $0$ et $A$ en degr\'e $1$. L\`a encore $C/C_2=A$
est plat sur $A$. On consid\`ere alors les triades $\Lpont=0 \fl 0 \fl C$,
$\Lpont^{(1)}=0
\fl 0 \fl C_1$ et
$\Lpont^{(2)}=0 \fl 0 \fl C_2$. Les fl\`eches naturelles correspondant aux inclusions de
$C_1$ et
$C_2$ dans $C$ sont des \ps, donc ces triades sont pseudo-isomorphes. Les triades
$\Lpont^{(1)}$ et $\Lpont^{(2)}$ sont toutes deux minimales. En effet, c'est clair
pour $\Lpont^{(1)}$ car $C_1= m$ n'a pas de quotient plat sur $A$ non nul (car $m$
n'est pas plat puisque $A$ est de dimension  $2$). Pour $\Lpont^{(2)}$,  $C_2$
vaut $ m$ en degr\'e $0$ et
$A$  en degr\'e $1$ et la multiplication par $X$ est l'inclusion de $m$
dans
$A$, de sorte que $C_2$ n' a pas de sous-$R_A$-module dont le quotient soit plat sur
$A$ (le seul sous-$A$-module \`a quotient plat est $m$ qui n'est pas un
sous-$R_A$-module). 

\vskip 0.3 cm

 \noindent b) Sur un anneau quelconque, on peut avoir un morphisme $u$ entre deux triades
pseudo-isomorphes, qui induise un isomorphisme sur les foncteurs associ\'es, mais ne
soit ni un \ps, ni  compatible avec une \'equivalence \ps\ ``\`a la Verdier'' comme
en 1.19, ni avec une \'equivalence du type 1.33.

On reprend en effet les triades $\Lpont^{(2)}$ et $\Lpont^{(1)}$ de l'exemple a) et
on a un morphisme de triades $\Lpont^{(2)} \fl \Lpont^{(1)}$ et le morphisme des
foncteurs associ\'es est $\Tor_1^A(C_2, \tas) \fl \Tor_1^A(C_1, \tas)$ qui est un
isomorphisme en vertu de la suite exacte scind\'ee sur $A$ : $0 \fl A \fl C_2 \fl C_1
\fl 0$.

Cependant, on voit aussit\^ot que $u$ n'est pas un \ps, ni compatible avec une
\'equivalence du type 1.19 ou 1.33 (sinon $C_2 \fl C_1$ serait injective). 

\vskip 0.3 cm

 \noindent c) M\^eme sur un \avd, deux triades peuvent avoir m\^eme foncteur associ\'e
sans
\^etre pseudo-isomorphes.

On suppose que $A$ est un \avd\ d'uniformisante $a$ et de corps r\'esiduel $k$. On va
construire deux triades avec $C= k$ et $H= k(-2)$. Pour cela, on va exhiber deux
\'el\'ements de $\Ext^2_{R_A} (C,H)$, cf. 1.34.
 Le
groupe $\Ext^2_{R_A} (k,k(-2))$ se calcule avec le complexe de Koszul associ\'e \`a la
suite
$(a,X,Y,Z,T)$ (dont on num\'erote les \'el\'ements de $0$ \`a $4$) et on a $\Ext^2_{R_A}
(k,k(-2))= k(-1)^4
\oplus k^6$. Les premi\`eres composantes ont pour base les $e_{0,i}$, $i=1,\cdots,4$
et les autres  ont pour base les $e_{i,j}$ avec $1 \leq i < j \leq 4$. 

On prend d'abord l'\'el\'ement $e_{1,2}$ qui donne la triade $\Lpont = 0 \fl L_0 \Fl
{d_0} L_{-1}$ avec $L_0 = k[x] (-1)$, avec $x^2=0$,  $L_{-1}= k[y]$, avec $y^2=0$ et $d_0
(1)=
y$. On v\'erifie que cette triade a bien les $C$ et $H$ voulus. On note que cette
triade est minimale (car
$C$ est de torsion) et toute r\'esolution majeure $\Mpont$ de $\Lpont$ est
\'el\'ementaire.

On prend ensuite l'\'el\'ement nul de $\Ext^2$ qui donne la triade $\Lpont'= 0 \fl
k(-2) \Fl{0} k$. L\`a encore cette triade a les $C$ et $H$ voulus, elle est minimale et
toute r\'esolution majeure $\Mpont'$ est \'el\'ementaire.
 Cela montre que ces triades ne sont pas pseudo-isomorphes, sinon leurs r\'esolutions
majeures seraient les m\^emes \`a un \ps\ induisant un isomorphisme sur les conoyaux
pr\`es, ce qui donnerait le m\^eme \'el\'ement du $\Ext^2$.

Par ailleurs, en vertu de 1.5, les foncteurs associ\'es sont tous deux extensions de
$\Tor_1^A(k,\tas)$ par $k(-2) \T_A \tas$ et on voit qu'on a $V(Q) = \Hom (k,Q) \oplus
(k(-2) \T Q)$ comme $A$-module. De plus, la structure de $R_A$-module est
n\'ecessairement triviale, donc les foncteurs sont isomorphes.

\tarte {j) Annexe : g\'en\'eralisation}

Nous donnons dans ce paragraphe quelques \'el\'ements permettant de g\'en\'eraliser la
d\'efinition des triades en s'affranchissant de l'hypoth\`ese de platitude. La
d\'efinition est exactement celle de 1.10 dans laquelle on ne suppose plus les complexes
ad\'equats. On parlera alors de triades ``faibles''. On a encore l'existence des
r\'esolutions majeures (1.14). La diff\'erence essentielle r\'eside dans la d\'efinition
du foncteur
$V$ associ\'e qui est d\'efini ``en le pensant dans la cat\'egorie d\'eriv\'ee'' : 

\th {Proposition-D\'efinition 1.36}. Soit  $\Lpont$ une triade faible et $\Mpont \fl
\Lpont$ une r\'esolution majeure.
 Le {\bf foncteur associ\'e} \`a $\Lpont$ est le foncteur $V_{\Lpont}$ (ou simplement
$V$) de $ \Mod_A$ dans $\Gr_{R_A}$ d\'efini pour tout $A$-module $Q$ par la formule
$V(Q)= h_0 (\Mpont
\T_A Q)$. Il est ind\'ependant du choix de la triade majeure $\Mpont$. L'application
qui \`a une triade $\Lpont$ associe le foncteur $V_{\Lpont}$ est fonctorielle.

\dem Le fait que $V$ ne d\'epende pas de $\Mpont$ r\'esulte du lemme 1.37 ci-dessous et
de la d\'emonstra\-tion de [AG] III 12.3 (car les modules $M_i$ sont plats sur $A$). La
fonctorialit\'e vient du lemme 1.38 ci-dessous.

\th {Lemme 1.37}. Soit $\Lpont$ une triade faible et soient $i : \Mpont \fl \Lpont$, $i'
:
\Mpont' \fl \Lpont$ deux r\'esolutions majeures de $\Lpont$. Alors, il existe un
morphisme de triades $g : \Mpont \fl \Mpont'$ tel que $i$ et $i'g$ soient homotopes. En
particulier les isomorphismes induits par $i$ et $i'g$ sur les c\oe urs et les conoyaux
de ces triades co\"\i ncident.

\th {Lemme 1.38}. Soit $f : \Lpont \fl \Lpont'$ un morphisme de triades faibles et $i' :
\Fpont' \fl \Lpont'$ une r\'esolution libre. Alors, il existe une r\'esolution libre
$i : \Fpont \fl \Lpont$ et un morphisme de complexes $g : \Fpont \fl \Fpont '$ tels
que les fl\`eches $fi$ et $i'g$ soient homotopes.

Ces lemmes, que nous laissons au lecteur, se d\'emontrent en travaillant dans la
cat\'egorie triangul\'ee des complexes de $R_A$-modules, cf. [RD].

\remas {1.39} \lign 
1) Le lemme 1.5 vaut, pour le nouveau foncteur $V$, sans hypoth\`ese de
platitude. \lign 2)  Attention, le foncteur $V$ d\'efini ci-dessus n'est pas,
en g\'en\'eral,
\'egal au foncteur $V_0$  na\"\i f d\'efini par  $V_0 (Q) = h_0 (\Lpont \T_A Q)$. 
 Consid\'erons par  exemple  un anneau $A$ local d'id\'eal maximal $m$ et de corps
r\'esiduel $k$ et  la triade
$0
\fl A
\Fl{p} k$ o\`u $p$ est la projection canonique. On a donc $C=0$ et $H= m$, d'o\`u
$V(Q) = m\T_A Q$ en vertu de 1.5, tandis que $V_0 (Q)$ est \'egal \`a $mQ$. 
\lign  3)  Bien entendu, dans la d\'efinition des pseudo-isomorphismes entre triades
g\'en\'erales, la condition portant sur le foncteur na\"\i f $V_0= h_0(\Lpont \T_A
\tas)$ doit
\^etre remplac\'ee par la condition analogue avec le nouveau foncteur $V$.

\titre {2. Faisceaux triadiques}

L'objectif de ce paragraphe est de d\'efinir les faisceaux triadiques qui sont
directement li\'es  aux triades majeures et de montrer que les notions de
pseudo-isomorphismes, au sens de [HMDP1] pour les faisceaux et du \S \ 1 pour les
triades,  se correspondent.

\vskip 0.3 cm

Rappelons que pour un faisceau coh\'erent $\cF$ sur $\bP^3_A$, plat sur $A$,
l'expression ``$H^i_* \cF$ commute au changement de base''  signifie 
que, pour tout $A$-module $Q$, l'homomorphisme canonique $(H^i_* \cF)
\T_A Q  \fl H^i_*(\cF
\T_A Q)$ est un isomorphisme. Cette condition est \'equivalente \`a l'exactitude \`a
droite du   foncteur
$H^i_*(\cF
\T_A \tas)$ ou encore, via la suite longue de cohomologie, \`a l'exactitude \`a gauche de
$H^{i+1}_*(\cF
\T_A \tas)$.  Rappelons aussi, cf. [AG] III
12.11 que si
$H^i_*
\cF$ commute au changement de base, on a les \'equivalences :  $H^i_* \cF$ est plat sur
$A$
$\equi\;$ $H^i(\cF \T_A \tas)$ exact $\equi\;$ 
$H^{i-1}_*
\cF$ commute au changement de base.
 
 \tarte {a) Faisceaux triadiques}

\th {Proposition 2.1}. Soit $\cN$ un faisceau coh\'erent sur $\bP^3_A$. Les conditions
suivantes sont \'equivalentes : \lign
i) Il existe une suite exacte 
$$0 \fl \cN \fl \cL_1 \fl \cL_0 \fl \cL_{-1} \fl 0$$
avec les faisceaux $\cL_i$ dissoci\'es. \lign
ii) Le faisceau $\cN$ est localement libre et il existe une suite exacte 
$$0 \fl \cM^{-1} \fl \cM^0 \fl \cM^1 \fl \cN^\vee \fl 0$$
avec les $\cM^i$ dissoci\'es. \lign
iii) Le faisceau $\cN$ est localement libre et le foncteur $H^2_* (\cN \T_A \tas)$ de
$\Mod_A$ dans $\Gr_{R_A}$ est exact \`a droite.
\lign iv) Le faisceau $\cN$ est localement libre et le foncteur $H^0_* (\cN^\vee \T_A
\tas)$ est exact \`a droite. \lign
v) Le faisceau $\cN$ est localement libre et le foncteur $H^3_* (\cN \T_A \tas)$ de
$\Mod_A$ dans $\Gr_{R_A}$ est exact \`a gauche. \lign
vi) Le faisceau $\cN$ est localement libre et le foncteur $H^1_* (\cN^\vee \T_A \tas)$
de
$\Mod_A$ dans $\Gr_{R_A}$ est exact \`a gauche.

\dem L'\'equivalence de {\sl i)} et {\sl ii)} est claire par dualit\'e. (On notera que
la condition {\sl i)} implique que le faisceau est localement libre.) De m\^eme
l'\'equivalence de {\sl iii)} et {\sl v)} (resp. de {\sl iv)} et {\sl vi)}) est
\'evidente.

\vskip 0.2 cm

\noindent {\sl iii)} $\equi$ {\sl vi)}. Il suffit de montrer les assertions analogues
pour les foncteurs $H^2 (\cN (l) \T_A \tas)$ et $H^1 (\cN^\vee(-l-4) \T_A \tas)$ pour
tout $l
\in \bZ$. En vertu de [H] 7.4, le foncteur $F= H^2 (\cN (l) \T_A \tas)$ est dual  du foncteur $F^*=H^1 (\cN^\vee(-l-4)
\T_A \tas)$ au sens
de la dualit\'e des foncteurs coh\'erents, cf. [H, \S 4], de sorte que $F$ exact \`a
droite
\'equivaut \`a $F^*$ exact \`a gauche.

\vskip 0.2 cm

\noindent {\sl ii) $\impl$ iv)}. La suite exacte de {\sl ii)} donne deux suites exactes
courtes $0 \fl \cM^{-1} \fl \cM^0 \fl \cK \fl 0$ et 
$0 \fl \cK \fl \cM^1 \fl \cN^\vee \fl 0$ avec $\cK$ localement libre. Comme les $\cM^i$
sont dissoci\'es, le foncteur $H^1_* (\cK \T \tas)$ est nul. Il en r\'esulte que le
morphisme de foncteurs $H^0_* (\cM^1 \T \tas) \fl H^0_*(\cN^\vee \T \tas)$ est un
\'epimorphisme. Comme $\cM^1$ est dissoci\'e, $H^0_* (\cM^1 \T \tas) $ est exact \`a
droite, donc aussi $ H^0_*(\cN^\vee \T \tas)$.

\vskip 0.2 cm

\noindent {\sl iv) $\impl$ ii)}. On pose $T = H^0_* (\cN^\vee)$ et on consid\`ere une
r\'esolution libre $F_1 \Fl{v} F_0 \fl T \fl 0$. Soit $G= \Ker v$. On a la suite des
faisceaux associ\'es $0 \fl \cG \fl \cF_1 \fl \cF_0 \fl \cN^\vee \fl 0$ avec les
$\cF_i$ dissoci\'es et $\cG$ localement libre. Il s'agit de montrer que $\cG$ est
dissoci\'e. On coupe cette suite en deux :
$0 \fl \cG \fl \cF_1 \fl \cS \fl 0$ et $0 \fl \cS \fl \cF_0 \fl \cN^\vee \fl 0$. La
fl\`eche 
$F_0 = H^0_*(\cF_0) \fl H^0_*(\cN^\vee)=T$ est surjective par construction. Par
ailleurs, le foncteur $H^0_* (\cN^\vee \T \tas)$ \'etant exact \`a droite, on a 
$H^0_* (\cN^\vee \T \tas)=H^0_* (\cN^\vee) \T \tas$. Comme la m\^eme
chose est vraie pour le faisceau dissoci\'e $\cF_0$, on en conclut que le morphisme de
foncteurs $H^0_*(\cF_0 \T \tas) \fl H^0_*(\cN^\vee \T \tas)$ est surjectif, ce qui
montre que $H^1_*(\cS \T \tas)= H^2_*(\cG \T\tas)$ est nul et que $H^0_*(\cS \T \tas)$
est exact
\`a droite. Le m\^eme argument avec l'autre suite exacte montre que $H^0_*(\cF_1 \T\tas)
\fl H^0_*(\cS \T\tas) $ est surjectif, donc que $H^1_*(\cG \T\tas)$ est nul.  En vertu
de [H] 7.7 on voit que $\cG$ est dissoci\'e, ce qui montre {\sl ii)}.

\rema {2.2} La condition $H^2_* (\cN \T_A \tas)$ exact \`a droite  est \'equivalente
\`a la platitude de $H^3_* \cN$ sur $A$ (cf. [AG] III 12.11 ou le fait que le foncteur
$H^3_* (\cN \T_A \tas)$ \'etant exact \`a droite, il est exact \`a gauche si et
seulement si $H^3_* \cN$ est plat sur $A$).

\th {D\'efinition 2.3}. Un faisceau sur $\bP^3_A$ v\'erifiant les conditions
\'equivalentes de 2.1 est appel\'e un faisceau {\bf triadique}. 

Un faisceau dissoci\'e
est \'evidemment triadique. Une somme directe de faisceaux tria\-diques est triadique.

\tarte {b) Triades et faisceaux triadiques}

Dans ce paragraphe nous montrons qu'il y a une correspondance entre triades majeures et
faisceaux triadiques. 

\th {Proposition 2.4}. (Triade associ\'ee \`a un faisceau triadique) Soit $\cN$ un
faisceau triadique et soit
$$0 \fl \cN \fl \cL_1 \fl \cL_0 \fl \cL_{-1} \fl 0 \leqno {(1)}$$
une suite exacte comme en 2.1. Alors, si on pose $L_i= H^0_* (\cL_i)$ :\lign
1) le complexe $\Lpont =(L_1 \fl L_0 \fl L_{-1})$ d\'eduit de $(1)$ est une triade
majeure de noyau
$N= H^0_* (\cN)$, \lign
2) si on note $H$ et $C$  le c\oe ur et le conoyau de $\Lpont$  et
$V_{\Lpont}$ le foncteur associ\'e (cf. 1.3), on a 
$$ H^1_* (\cN\T_A \tas) = h_0(\Lpont\T_A \tas)= V_{\Lpont}, \quad \hbox {en particulier,
on a }\; H= H^1_* \cN,$$
$$ H^2_* (\cN\T \tas) = C \T_A \tas, \quad  \hbox {en particulier, on a } \; C= H^2_*
\cN,
$$ 3) le  dual de $N$ sur $R_A$ : $N^\vee = \Hom_{R_A} (N,R_A)$,  est \'egal \`a $H^0_*
(\cN^\vee)$ et on a une suite exacte de $R_A$-modules :
$$0 \fl L_{-1}^\vee \fl L_0^\vee \fl L_1^\vee \fl N^\vee \fl 0,$$
4) si on a deux suites exactes du type $(1)$ correspondant au m\^eme faisceau $\cN$, les
triades associ\'ees sont isomorphes
\`a homotopie pr\`es,\lign
 5) on a 
$\Ext^1_{R_A}(N,R_A)
\simeq H^1_* \cN^\vee \simeq C^*(4)$. Si on suppose $A$ int\`egre, 
$\Ext^1_{R_A}(N,R_A)$ est nul si et seulement si $C= H^2_* (\cN)$ est de torsion.
Autrement dit, $\Lpont$ est \'el\'ementaire (cf. 1.28) si et seulement si $\cN$ est
extraverti (cf. [HMDP1] 2.6).

\dem Montrons d'abord 2). On coupe la suite exacte (1) en deux :
$0 \fl \cN \fl \cL_1 \fl \cK \fl 0$ et $0 \fl \cK \fl \cL_0 \fl \cL_1 \fl 0$. En tenant
compte du fait que $H^0_* (\cL_i \T \tas) = H^0_* (\cL_i )\T \tas$ car $\cL_i$ est
dissoci\'e on voit qu'on a $ H^2_* (\cN\T \tas)  = H^1_*(\cK \T\tas) = C \T\tas$ et 
$H^1_* (\cN\T_A \tas) = h_0(\Lpont\T_A \tas)$. Comme $\cN$ est localement libre,
$H=H^1_* (\cN)$ et
$C=H^2_* (\cN)$ sont de type fini sur $A$, de sorte que $\Lpont$ est bien une triade.

Pour le point 3) rappelons tout d'abord  que si
$A$ est un anneau noeth\'erien  et 
$M$  un
$R_A$-module gradu\'e de type fini sur $A$ on a $\Ext^i_{R_A} (M, R_A)
=0$ pour $0
\leq i
\leq 3$ (cf. [E] 18.4 : il suffit de noter que l'annulateur de $M$ contient une suite
r\'eguli\`ere de la forme $(X^r, Y^r, Z^r, T^r)$ pour $r$ assez grand).

Consid\'erons la triade  $\Lpont = (L_1 \Fl {d_1} L_0 \Fl {d_0} L_{-1})$ et posons
$B =
\Im d_1$, $Z =
\Ker d_0$, $K = \Im d_0$. On a les suites exactes
$0 \fl N \fl L_1 \fl B \fl 0$, $0 \fl B \fl Z \fl H \fl 0$, $0 \fl Z \fl L_0 \fl K
\fl 0$ et $0 \fl K \fl L_{-1} \fl C \fl 0$.
Comme $H$ et $C$ sont de type fini sur $A$, le rappel ci-dessus montre qu'on a 
$\Ext^i_{R_A} (H,R_A)= \Ext^i_{R_A} (C,R_A)=0$ pour $i \leq 3$ et on en d\'eduit
$\Ext^1_{R_A} (B,R_A)= \Ext^1_{R_A} (Z,R_A) = \Ext^2_{R_A} (K,R_A)=\Ext^3_{R_A}
(C,R_A)=0$ et $\Ext^1_{R_A} (K,R_A)=\Ext^2_{R_A} (C,R_A)=0$. En d\'efinitive on a
la suite exacte sur les mo\-du\-les duaux :
$$0 \fl L_{-1}^\vee \Fl {d_0^\vee} L_0^\vee \Fl{d_1^\vee} L_1^\vee \fl N^\vee \fl
0.$$

Pour 4),   on a les deux triades $\Lpont$ et $\Lpont'$ et il  r\'esulte de 3) que
$\Lpont^\vee$ et $\Lpont'^\vee$ sont deux r\'esolutions libres de $N^\vee$,
donc isomorphes
\`a homotopie pr\`es, donc aussi 
$\Lpont$ et $\Lpont'$.

 Montrons le point 5) :
comme $\cO_{\bP}$ n'a pas de $H^1$, on a $\Ext^1_{R_A} (N, R_A(n)) ^0 =
\Ext^1_{\cO_{\bP}} (\cN, \cO_{\bP}(n))$ et ce dernier terme n'est autre que
$H^1\cN^\vee (n)$ en vertu de [AG] 6.3 et 6.7. On a aussi $H^1_*\cN^\vee =
H^2_*\cK^\vee$ et la suite exacte 
$0 \fl H^2_*\cK^\vee \fl H^3_*\cL_{-1}^\vee \Fl {u} H^3_*\cL_{0}^\vee $. Mais, on a vu
(cf. 0.1) qu'on a $H^3_*\cL_{i}^\vee  \simeq (L_i(-4))^*$ et la suite exacte $L_0
\Fl {d_0} L_{-1} \fl C \fl 0$ montre par dualit\'e que $\Ker u$ n'est autre que
$C^* (4)$.
On suppose  $A$ int\`egre. Si $C$ est de torsion on a $C^*=0$. R\'eciproquement,
si
$C^*$ est nul on a 
$C^{**} =0$, donc la fl\`eche $C \fl C^{**}$ est nulle ce qui implique que $C$
est de torsion.

\vskip 0.3 cm

\th {Proposition 2.5}. (Faisceau triadique associ\'e \`a une triade majeure) \lign
Soit $\Lpont$ une triade majeure, $N$ son noyau, $\cN$ le faisceau associ\'e. Alors,
$\cN$ est triadique. Cette op\'eration d\'efinit un foncteur de la cat\'egorie des
triades majeures munie des morphismes de complexes \`a homotopie pr\`es dans celle des
faisceaux triadiques sur
$\bP^3_A$ et ce foncteur  est une \'equivalence de cat\'egories.

\dem 
Si on pose $\Lpont=(L_1 \Fl {d_1} L_0 \Fl {d_0}
L_{-1})$,  le fait que $H$ et $C$ soient de type fini sur $A$ implique que les
faisceaux associ\'es sont nuls, d'o\`u l'exactitude de la suite 
$$0 \fl \cN \fl \cL_1 \Fl{\wi d_1} \cL_0 \Fl{\wi d_0} \cL_{-1} \fl 0, $$ 
qui montre que $\cN$ est triadique, puisque les $L_i$ et donc les $\cL_i$ sont
dissoci\'es.

Il est clair que cette correspondance est fonctorielle et c'est encore vrai \`a
homotopie pr\`es, puisque deux morphismes de complexes homotopes induisent la m\^eme
application sur l'homologie et en particulier sur le  $h_1$. Le fait que ce
foncteur soit surjectif vient de 2.4. 

Si $f: \cN \fl \cN'$ est un morphisme de faisceaux, il induit $f^\vee : \cN'^\vee \fl
\cN^\vee$ et, en passant aux sections globales,  $H^0_* f^\vee : N'^\vee \fl
N^\vee$. Comme
$\Lpont^\vee$ et $\Lpont'^\vee$ sont des r\'esolutions projectives de $N^\vee$ et
$N'^\vee$ (cf. 2.4.3), $H^0_* f^\vee$ se rel\`eve en un morphisme de $\Lpont'^\vee$ dans
$\Lpont^\vee$ qui donne, par dualit\'e, le morphisme de triades cherch\'e.

Si deux morphismes $f,g : \Lpont \fl \Lpont'$
donnent le m\^eme morphisme $\cN \fl \cN'$, ils donnent aussi le m\^eme morphisme
$N'^\vee \fl N^\vee$ et, comme $\Lpont^\vee$ et
$\Lpont'^\vee$ sont des r\'esolutions projectives de $ N^\vee$ et $ N'^\vee$, $f^\vee$ et
$g^\vee$ sont homotopes, donc aussi $f$ et $g$.

\tarte {c) Pseudo-isomorphismes}

Rappelons la d\'efinition des pseudo-isomorphismes de faisceaux dans le cas des
faisceaux localement libres (cf. [HMDP1] 2.1 et 2.2.0) :

\th {D\'efinition 2.6}. Un morphisme $f : \cN \fl \cN'$ de faisceaux localement libres
est appel\'e un {\bf pseudo-isomorphisme}  (en abr\'eg\'e un \ps) s'il
induit :
\lign 1) un isomorphisme de foncteurs $H^1_*(\cN\T_A \tas) \fl H^1_*(\cN'\T_A \tas)$, \lign
2) un monomorphisme de foncteurs $H^2_*(\cN\T_A \tas) \fl H^2_*(\cN'\T_A \tas)$.\lign
Deux faisceaux  sont dits pseudo-isomorphes s'il existe une cha\^\i ne
de \ps\ qui les joint.

\vskip 0.3 cm

Dans la correspondance \'etudi\'ee au paragraphe pr\'ec\'edent, les \ps\ de triades et
de faisceaux se correspondent :

\th {Proposition 2.7}. Soit $u : \Lpont \fl \Lpont'$ un morphisme de triades majeures
et $\wi u : \cN \fl \cN'$ le morphisme correspondant sur les faisceaux triadiques.
Alors $u$ est un \ps\ si et seulement si $\wi u$ est un \ps\ c'est-\`a-dire si et
seulement si on a les deux conditions suivantes : \lign a) $\wi u$ induit un isomorphisme
$H^1_*(\cN)
\fl H^1_*(\cN')$, \lign b) $\wi u$ induit une injection $H^2_*(\cN) \fl H^2_*(\cN')$ \`a
quotient plat sur $A$.

\dem La premi\`ere assertion r\'esulte de 2.4 et la seconde de 1.16.

\vskip 0.3 cm

\rema {2.8} Il  r\'esulte   de 2.7 et 1.22 que si deux triades majeures
sont pseudo-isomorphes il en est est de m\^eme des faisceaux associ\'es. La
r\'eciproque est vraie, cf. 2.14, mais pas imm\'ediate. En effet, si on a une cha\^\i
ne de \ps\ de faisceaux dont les extr\'emit\'es sont triadiques, les faisceaux
interm\'ediaires ne le sont pas n\'ecessairement, de sorte qu'on ne peut pas appliquer
2.7. 

\vskip 0.3 cm

Rappelons  la caract\'erisation suivante des \ps\ de faisceaux (cf. [HMDP1] 2.4) :

\th {Proposition 2.9}. Soit $f : \cN \fl \cN'$ un morphisme de faisceaux localement
libres. Alors, $f$ est un \ps\ si et seulement si il existe  un faisceau dissoci\'e $\cL$
et un morphisme $p : \cL \fl \cN'$ tel que l'on ait une suite exacte 
$$0 \fl \cS \fl \cN \oplus \cL \,\Fl {(f,p)} \, \cN' \fl 0$$
avec $\cS$ dissoci\'e.

\th {Corollaire 2.10}. Soit $f : \cN \fl \cN'$ un pseudo-isomorphisme  de faisceaux
triadiques qui induit un {\bf isomorphisme } de foncteurs $H^2_* (\cN \T_A \tas) \fl 
H^2_* (\cN' \T_A \tas)$. Alors, $\cN$ et $\cN'$ sont
 stablement isomorphes.

\dem   En vertu de 2.9 on a une suite exacte $0 \fl \cS \fl \cN \oplus \cL \,\Fl {(f,p)}
\, \cN' \fl 0$ et il suffit de montrer que cette suite est scind\'ee. Cette suite
correspond \`a un \'el\'ements $\xi \in \Ext^1_{\cO_\bP} (\cN', \cS)$ et il s'agit de
montrer que cet \'el\'ement est nul, ou encore que le morphisme 
$\f : \Ext^1_{\cO_\bP} (\cN', \cS) \fl \Ext^1_{\cO_\bP} (\cN\oplus \cL, \cS) =
\Ext^1_{\cO_\bP} (\cN, \cS)$ induit par
$f$ est injectif. Ce morphisme est encore \'egal au morphisme 
$H^1 (\cN'^\vee \T_{\cO_\bP} \cS )\fl H^1 (\cN^\vee \T_{\cO_\bP} \cS)$ induit par $f$ et
c'est la valeur en $A$ du morphisme de foncteurs 
$$\Phi : H^1( \cN'^\vee \T_{\cO_\bP} \cS \T_A \tas) \fl H^1 ( \cN^\vee \T_{\cO_\bP}\cS
\T_A \tas).$$  Par la dualit\'e de Serre exprim\'ee en termes de foncteurs coh\'erents,
cf. [H] 7.4,
$\Phi$ est dual du morphisme $\Phi^\vee : H^2
(\cN
\T_{\cO_\bP} \cS 
\T_{\cO_\bP}\om \T_A \tas) \fl  H^2 (\cN' \T_{\cO_\bP} \cS  \T_{\cO_\bP}\om \T_A \tas) $.
Comme
$\cS$ et
$\om = \cO_{\bP}(-4)$ sont dissoci\'ees ce dernier morphisme est somme de certains
composants du morphisme $H^2_* (\cN \T_A \tas) \fl 
H^2_* (\cN' \T_A \tas)$ qui est un isomorphisme par hypoth\`ese, d'o\`u la conclusion.

\tarte {d) Le lemme de Verdier inverse pour les faisceaux triadiques}

Nous avons montr\'e dans [HMDP1] 2.11  le lemme (dit de Verdier) suivant :
\smallskip

{\sl  Soient $\cN'$ et $\cN''$ 
des faisceaux 
localement libres sur $\bP^3_A$ pseudo-isomorphes. Il existe un faisceau  localement
libre  $\cN_0$ et des
\ps\  : 
$\cN' \lf \cN_0 \fl \cN''$.}
\smallskip

Le lemme analogue avec des faisceaux triadiques n'est pas vrai, mais on a le r\'esultat
suivant en sens inverse, comme dans le cas des triades (cf. 1.22) :

\th {Th\'eor\`eme 2.11}. Soient $\cN'$ et $\cN''$ 
des faisceaux 
triadiques sur $\bP^3_A$ pseudo-isomorphes. Il existe un faisceau 
triadique 
$\cN$ et des
\ps\  : 
$\cN' \fl \cN \lf \cN''$.

La d\'emonstration de ce th\'eor\`eme n\'ecessite quelques \'etapes. Le premier lemme
montre comment construire des faisceaux triadiques par une technique de d\'esaturation.
Cette technique sera utilis\'ee au chapitre 3 pour montrer l'existence de r\'esolutions
de type \sE\ et \sN\ triadiques.

\th {Lemme 2.12}. (Lemme de d\'esaturation) Soit $0 \fl \cE \fl \cF \fl \cG \fl 0$ une
suite exacte de faisceaux coh\'erents sur $\bP^3_A$. On suppose $\cE$ localement libre,
$\cF$ dissoci\'e et $\cG$ plat sur $A$. Il existe des faisceaux $\cE'$ et $\cF'$ tels que
$\cE'$ soit dual d'un faisceau triadique, que
$\cF'$ soit dissoci\'e et qu'on ait un diagramme commutatif de suites exactes, o\`u $f$
est dual d'un \ps\  :
$$\matrix { 0& \fl & \cE'& \fl &\cF' & \fl & \cG & \fl & 0 \cr
&&\Vf{f}&&\Vf{} &&\parallel \cr
0& \fl & \cE& \fl &\cF & \fl & \cG & \fl & 0. \cr
}$$

\dem 

Notons d'abord que, comme $\cE$ est localement libre, on a $\cExt^i_{\cO_\bP} (\cG,
\cO_\bP) =0$ pour $i \geq 1$.

 Soit $G$ l'image de $F = H^0_* \cF$ dans $ H^0_* \cG$. On sait que pour $n$ grand, on
a $G_n = H^0 \cG (n)$. Soit
$n_0$ un entier, tel que cette condition soit r\'ealis\'ee pour $n \geq n_0$, et tel que
le foncteur $H^0 (\cG (n) \T_A \tas
)$ soit exact pour $n \geq n_0$ (un tel entier existe en vertu des th\'eor\`emes de
cohomologie et changement de base, cf. [AG] III 12.11). Posons $G' = \bigoplus_{n \geq
n_0} H^0 \cG (n)$. On a donc $G' \subset G$ et $G'$ est obtenu \`a partir de $G$ par
``d\'esaturation''. Soit
$F'
\fl G'$  le d\'ebut d'une r\'esolution libre minimale de
$G'$. L'inclusion de $G'$ dans $G$ se rel\`eve en un morphisme $F' \fl F$. Si $E'$ est
le noyau de $F' \fl G'$ on obtient le diagramme annonc\'e en passant aux faisceaux
associ\'es.

Il reste \`a montrer que $\cE'$ est dual d'un faisceau triadique et que $f$ est dual
d'un
\ps. On note que $\cE'$ est plat sur $A$ car $\cG$ l'est et qu'on a
$\cExt^i_{\cO_\bP} (\cE',
\cO_\bP) =0$ pour $i \geq 1$, ce qui montre d\'ej\`a que $\cE'$ est localement libre.

On v\'erifie ensuite que $f$ est dual d'un \ps. Comme $\cG$ est plat les
suites exactes du diagramme ci-dessus restent exactes en tensorisant par un $A$-module
$Q$. Le diagramme obtenu en d\'eroulant les suites de cohomologie donne la conclusion.

Montrons enfin que le foncteur $H^0_* (\cE' \T_A \tas) $ est exact \`a droite. Il suffit
de le voir pour les foncteurs $H^0 (\cE'(n) \T_A \tas) $ pour tout $n$. Comme
$\cG$ est plat sur $A$ on a la suite exacte de foncteurs 
$0 \fl \cE'(n) \T_A \tas \fl \cF'(n) \T_A \tas \fl \cG(n) \T_A \tas \fl 0$. Pour $n
\leq n_0$, comme $F'$ est une r\'esolution minimale de $G'$  on a $F'_n=0$  et on en
d\'eduit $H^0 (\cE'(n) \T_A \tas) =0$. Pour $n \geq n_0$, comme $H^0 (\cG(n) \T_A \tas)$
commute au changement de base et comme $F'_n \fl G'_n$ est surjectif on a la suite
exacte de foncteurs
$$0 \fl H^0 (\cE'(n) \T_A \tas) \fl H^0 (\cF'(n) \T_A \tas) \fl H^0 (\cG(n) \T_A \tas)
\fl 0$$
et il en r\'esulte que $ H^0 (\cE'(n) \T_A \tas) $ est exact \`a droite (car $H^0
(\cF'(n) \T_A \tas)$ l'est).

\th {Corollaire 2.13}. \lign
1)  Soit $\cE$ un faisceau localement libre. Il
existe un faisceau dual d'un faisceau triadique $\cE'$  et un morphisme
 $\cE' \fl \cE$ qui est dual d'un \ps. \lign
2) Soit $\cN$ un faisceau localement libre. Il
existe un faisceau  triadique $\cN'$  et un \ps\
 $\cN \fl \cN'$.  En particulier, la classe de
pseudo-isomorphisme  d'un faisceau localement libre
contient un faisceau triadique.

\dem Par dualit\'e il suffit de prouver l'assertion avec $\cE$. En vertu de 2.12 il
suffit d'ins\'erer $\cE$ dans une suite exacte $0 \fl \cE \fl \cF \fl \cG \fl 0$ avec
$\cF$ dissoci\'e et $\cG$ plat. Pour cela on consid\`ere $\cE^\vee$, on prend le
d\'ebut d'une r\'esolution minimale de $N= H^0_* \cE^\vee$ : $0 \fl G^\vee \fl F^\vee \fl
N
\fl 0$, avec
$F^\vee$ libre, on passe aux faisceaux associ\'es et on dualise.

\vskip 0.3 cm

Nous pouvons maintenant montrer le lemme de Verdier inverse 2.11. On commence par prendre
un faisceau localement
libre  $\cN_0$ et des
\ps\  : 
$\cN' \lf \cN_0 \fl \cN''$ (cf. [HMDP1] 2.11). On en d\'eduit par dualit\'e une fl\`eche 
$\cN'^\vee \oplus \cN''^\vee \fl \cN_0^\vee$. On rajoute un faisceau dissoci\'e $\cP$
de sorte que la fl\`eche $\cP \oplus \cN'^\vee \oplus \cN''^\vee \fl \cN_0^\vee$soit
surjective. Soit $\cE$ le noyau de cette fl\`eche. C'est un faisceau localement libre
et si on pose $\cN= \cE^\vee$ on a la suite exacte de faisceaux localement libres :
$$0 \fl \cN_0 \fl \cP^\vee \oplus \cN' \oplus \cN''\,  \Fl {(u, \pi' , \pi'')} \, \cN \fl
0. \leqno {(1)}$$
Nous allons montrer que $\pi'$ et $\pi''$ sont des \ps, ce qui, avec 2.13 \'etablira le
th\'eor\`eme.

 La suite $(1)$  reste exacte en tensorisant par un $A$-module $Q$, de sorte qu'elle
donne naissance \`a une suite exacte de foncteurs :
$$ \eqalign {\cdots \fl &H^1 (\cN_0 \T_A \tas) \Fl{^t(j',j'')} H^1 (\cN' \T_A \tas) 
\oplus H^1 (\cN''
\T_A
\tas)  \Fl {(p',p'')} H^1 (\cN \T_A \tas)  \Fl {v}\cr
& H^2 (\cN_0 \T_A \tas) \Fl {^t(k',k'')} H^2 (\cN' \T_A \tas)  \oplus H^2 (\cN'' \T_A
\tas)
 \Fl {(q',q'')} H^2 (\cN \T_A \tas) \fl \cdots}$$ 
dans laquelle les fl\`eches $j'$ et $j''$ sont des isomorphismes et $k'$ et $k''$ des
monomorphismes. On en d\'eduit d'abord qu'on a $v=0$, puis que $p'$ et $p''$ sont des
isomorphismes et enfin que $q'$ et $q''$ sont des monomorphismes, donc que $\pi'$ et
$\pi''$ sont des \ps.

\vskip 0.3 cm

Le corollaire suivant (qui va donner le th\'eor\`eme de Rao, cf. 3.9) montre que la
correspondance entre triades majeures et faisceaux triadiques se comporte bien vis \`a
vis de la relation de pseudo-isomorphisme, cf. remarque  2.8 :

\th {Corollaire 2.14}. Soient $\Lpont'$ et $\Lpont''$ deux triades majeures et $\cN'$ et
$\cN''$ les faisceaux triadiques associ\'es. Alors, $\Lpont'$ et $\Lpont''$ sont
pseudo-isomorphes si et seulement si $\cN'$ et $\cN''$  le sont.

\dem Cela r\'esulte de 2.7 et 2.11.

\tarte {e) Le cas d'un \avd}

Dans ce paragraphe nous supposons que $A$ est un \avd\ d'uniformisante $a$. Rappelons,
cf. [HMDP1] 2.6, qu'un faisceau localement libre 
 $\cN$ est dit extraverti s'il v\'erifie $H^1_* \cN^\vee=0$.

\th {Proposition 2.15}. Soit $A$ un \avd\ et soit $\cN$ un faisceau localement libre
extraverti. Alors
$\cN$ est triadique. La classe de pseudo-isomorphisme  de tout faisceau
localement libre contient un unique faisceau extraverti (donc triadique) minimal, i.e.
sans facteur direct dissoci\'e. Un tel faisceau sera dit {\bf \'el\'ementaire}. La
triade associ\'ee est \'el\'ementaire.

\dem Posons $M= H^0_* \cN^\vee$ et consid\'erons une pr\'esentation minimale de $M$~:
$0 \fl M_2 \fl L_1 \fl L_0 \Fl {p} M \fl 0$ avec $L_i$ libre. Soit $K$ le noyau de $p$.
 On a la suite exacte des faisceaux associ\'es :
$0 \fl \cM_2 \fl \cL_1 \fl \cL_0 \fl  \cN^\vee \fl 0$ et il s'agit de montrer que
$\cM_2$ est dissoci\'e. Pour cela, en vertu de [H] 7.7, il suffit de montrer que les
foncteurs $H^i_* (\cM_2 \T_A \tas)$ sont nuls pour $i=1,2$. 

Par hypoth\`ese $\cN$  v\'erifie $H^1_* \cN^\vee=0$. On a donc aussi, si $\cK$ est le
faisceau associ\'e \`a $K$, $H^2_* \cK=0$, d'o\`u la suite exacte 
$0 \fl H^3_* \cM_2 \fl H^3_* \cL_1 \fl H^3_* \cK \fl 0$. Comme $\cL_1$ est dissoci\'e, 
$H^3_*
\cL_1$ est plat, donc sans torsion et il en r\'esulte que  $H^3_* \cM_2 $ est sans
torsion, donc plat, puisque $A$ est un \avd. En vertu de [AG] III 12.11, il en r\'esulte
que $H^2_*(\cM_2 \T_A \tas)$ est exact \`a droite (c'est-\`a-dire que $H^2_*\cM_2$
commute au changement de base). Mais, par construction, la fl\`eche $p$ est surjective,
donc on a $H^1_* \cK = H^2_* \cM_2 =0$ et le foncteur $H^2_*(\cM_2 \T_A \tas)$ est nul.
Toujours par [AG] III 12.11, il en r\'esulte que $H^1_*(\cM_2 \T_A \tas)$ commute au
changement de base. Mais, comme on a $K= H^0_* \cK$, la surjectivit\'e de la fl\`eche 
$L_1 \fl K$ montre qu'on a $H^1_* \cM_2=0$, donc le foncteur correspondant est nul, cqfd.

 Le faisceau $\cN$ est donc \`a la fois extraverti et triadique. Quitte
\`a lui retirer un \'eventuel facteur dissoci\'e on peut le supposer minimal et il est
alors unique dans sa classe de pseudo-isomorphisme, cf. [HMDP1] 2.14.

\rema {2.16} On peut retrouver ce r\'esultat avec les r\'esultats du paragraphe 1. En
vertu de [HMDP1] 2.16, il suffit de montrer que toute classe de pseudo-isomorphisme
contient un faisceau triadique extraverti. En vertu de 2.13 elle contient un
$\cN$ triadique. Soit
$\Lpont$ la triade majeure associ\'ee. Il existe une triade \'el\'ementaire $\Mpont$
avec un \ps\
$\Mpont
\fl
\Lpont$ (cf. 1.32). Le faisceau triadique $\cN_0$ associ\'e \`a $\Mpont$ est tel que
$H^2_*
{\cN_0}$ est de torsion (cf.  1.32.1 et 2.4.5), de sorte que $\cN_0$ est extraverti et
triadique.

\vskip 0.5 cm

\titre {3. Courbes et triades : les th\'eor\`emes de Rao}

\tarte {a) R\'esolutions de type \sE\ et \sN\ triadiques}

\th {Th\'eor\`eme 3.1}. Soit $\cC$ une famille de courbes de $\bP^3_A$. \lign
1) Il existe
une r\'esolution de $\cC$ de type \sE\ cotriadique, i.e.
une suite exacte
$$  0 \fl \cE \fl \cF\fl
\cJ_\cC
\fl 0 $$
o\`u $\cF$  est dissoci\'e  et o\`u $\cE$ est  dual d'un faisceau triadique.
\lign 2) Il existe
une r\'esolution de $\cC$  de type \sN\ triadique, i.e., une suite exacte $$0 \fl
\cP
\fl \cN \fl \cJ_\cC
\fl 0 $$ o\`u $\cP$ est un faisceau dissoci\'e et
$\cN$ un faisceau triadique.

\dem Commen\c cons par le type \sE. On sait (cf. [HMDP1] 2.20) qu'il existe une
r\'esolution de type \sE. Pour obtenir une r\'esolution cotriadique il suffit d'appliquer
le lemme de d\'esaturation 2.12. 

\vskip 0.2 cm

Passons au type \sN. Vu l'importance de ce cas, nous allons donner deux d\'emonstrations de
l'existence de la r\'esolution, la premi\`ere par liaison \`a partir de l'existence de la
r\'esolution de type
\sE, l'autre, plus directe, en tronquant le complexe de
\v Cech.

\vskip 0.2 cm

{\it Premi\`ere d\'emonstration} 
\vskip 0.2 cm

 On effectue une liaison
\'el\'ementaire par une intersection compl\`ete globale $\cD$ de degr\'es $s$ et $t$ qui transforme
$\cC$ en
$\cC'$. On prend une r\'esolution de type \sE\ cotriadique de $\cC'$ obtenue par la
m\'ethode du lemme de d\'esaturation. Pr\'ecis\'ement, on choisit un entier $n_0$ tel que
pour
$n\geq n_0$ le foncteur $H^1 (\cE (n) \T_A \tas )$ soit nul (cf. [AG] III 12.11), on
pose  $I = \bigoplus_{n \geq n_0} H^0 \cJ_{\cC'} (n)$ et on prend un d\'ebut de
r\'esolution libre $F$ de
$I$ qui conduit \`a la suite exacte :
$0 \fl \cE \fl \cF \Fl{f} \cJ_{\cC'} \fl 0$.
Soit $$0 \fl \cO_\bP (-s-t) \fl \cO_\bP (-s) \oplus \cO_\bP (-t) \Fl {g} \cJ_\cD \fl
0$$ la r\'esolution de $\cD$. L'inclusion $\cJ_\cD \subset \cJ_{\cC'}$ fournit une suite
exacte
 $$ 0 \fl \cE' \fl \cF \oplus \cO_\bP (-s) \oplus \cO_\bP (-t) \Fl {f+g} \cJ_{\cC'} \fl
0.$$
Montrons que le faisceau $\cE'$ est  dual d'un triadique, c'est-\`a-dire que, pour tout
$n$, le foncteur 
$H^0 (\cE' (n)
\T_A
\tas )$ est exact \`a droite. On a une suite exacte 
$$0 \fl \cE \fl \cE' \fl \cO_\bP (-s) \oplus \cO_\bP (-t) \fl 0.$$ Pour $n\geq n_0$ on a
$H^1 (\cE (n) \T_A \tas )=0$ et il en r\'esulte que $H^0 (\cE' (n) \T_A \tas )$ est exact
\`a droite (car c'est vrai pour
$\cE$  et pour les dissoci\'es). 

Par ailleurs, on a une suite exacte 
$$0 \fl \cO_\bP (-s-t) \fl \cE' \fl \cF \fl \cJ_{\cC'/\cD} \fl 0$$
et, pour $n <n_0$ on a $H^0 (\cF (n) \T_A \tas )=0$, donc 
$H^0 (\cE' (n) \T_A \tas )= H^0 (\cO_\bP(n-s-t) (n) \T_A \tas )$ et ce dernier foncteur
est exact \`a droite.

Il r\'esulte alors de [HMDP1] 2.25 que $\cC$ admet une r\'esolution de type \sN\ dont
le faisceau $\cN$ est $\cE'^\vee (-s-t) \oplus \cO_\bP (-s) \oplus \cO_\bP (-t)$ qui est
triadique comme annonc\'e.

\vskip 0.3 cm

{\it Deuxi\`eme d\'emonstration}

\vskip 0.2 cm

  On consid\`ere le
complexe de \v Cech de $\cJ_C$ :
$$\bR
\G_*
\cJ_{\cC}= (0 \fl M^0 \Fl{\de^0} M^1 \Fl{\de^1} M^2 \Fl{\de^2} M^3 \fl 0)$$ vu
comme
\'el\'ement de la cat\'egorie d\'eriv\'ee
$\cD^+(R_A)$. Les
modules $M^i$ sont plats sur $A$. On
consid\`ere le complexe  tronqu\'e
 $\s_{\leq 2}(\bR\G_*\cJ_{\cC})$, cf. [RD] p. 69.
Il s'agit du complexe
$M^{^\tas}= (M^0 \fl M^1
\fl \Ker \de^2)$ et il  calcule les groupes
$H^i_*\cJ_\cC$ pour
$i=0,1,2$. Comme $M^3/\Im \de^2 = H^3_* \cJ_\cC = H^3_*
\cO_\bP$ est plat sur $A$ il en est de m\^eme de $\Ker \de^2$. Il en r\'esulte que  
$M^{^\tas}$ calcule aussi les foncteurs 
$H^i_*(\cJ_\cC\T_A\tas )$ pour $i=0,1,2$. 

 On choisit ensuite un entier $r$ assez petit pour que les foncteurs  $H^0 (\cJ_\cC (n) \T_A \tas)
= H^1
(\cJ_\cC (n)\T_A \tas)$ soient nuls pour
$n
\leq r$ (on dira que $r$ est convenable) et on consid\`ere le complexe  $M^{^\tas}_{>r} $
tronqu\'e  cette fois par rapport \`a la graduation des
$R_A$-modules. Les modules de cohomologie de ce complexe sont de type fini sur $R_A$, de sorte
qu'en vertu de [AG] III 12.3 il existe un complexe $S_r^{^\tas}$  form\'e de modules libres  de type fini 
tel que l'on ait
$M^{^\tas}_{>r} = S_r^{^\tas}$ dans la cat\'egorie d\'eriv\'ee. Ce complexe est form\'e
de modules $S_r^i$ avec $i \leq 2$, mais, attention, il
 contient
{\it a priori} des termes $S_r^i$ pour $i<0$.
La proposition suivante donne alors la conclusion voulue :

\th {Proposition 3.2}.  On reprend les
notations ci-dessus :   $\cC$ est une famille de courbes,
$M^{^\tas}= \s_{\leq 2}(\bR\G_*\cJ_{\cC})$ le complexe de \v Cech tronqu\'e,  $r$ un
entier convenable et 
$S_r^{^\tas}$ un complexe form\'e de modules
libres  de type fini  tel que l'on ait $M^{^\tas}_{>r} = S_r^{^\tas}$ dans la
cat\'egorie d\'eriv\'ee. On pose
$L_1= S^0$, $L_0= S^1$ et $L_{-1}= S^2$. \lign
1)  Le complexe ${\Lpont}_r= (L_1
\fl L_0
\fl L_{-1})$  est une triade majeure.  \lign
2)
Si $r$ et $r'$ sont deux entiers convenables avec $r'<r$ on a un \ps\ naturel $u_{r,r'}
:{\Lpont}_r
\fl {\Lpont}_{r'}$ de sorte que toutes les triades ${\Lpont}_r$ sont
pseudo-isomorphes. Si on a trois entiers convenables avec
$r''<r'<r$ on a
$u_{r,r'} u_{r',r''} = u_{r,r''}$ \`a homotopie pr\`es. \lign
3)  Soit $N$ le
noyau de ${\Lpont}_r$ et soit $\cN$ le faisceau associ\'e. Il existe une r\'esolution de type
\sN\ triadique de $\cC$ :
$0 \fl \cP \fl \cN \fl \cJ_{\cC} \fl 0$.

\dem Montrons 1). Comme le complexe ${\Lpont}_r$  est form\'e de $R_A$-modules libres de
type fini, il s'agit de voir que ses groupes d'homologie $h_0$ et $h_{-1}$ sont de type
fini sur
$A$. Pour $h_0$ c'est clair car les complexes  consid\'er\'es calculent tous 
$H^1_*(\cJ_\cC \T_A \tas)$. Le groupe $h_{-1}$ n'est
autre que 
$(H^2_* \cJ_\cC)_{>r}$, donc de type fini sur $A$. 

Pour le point 2), si on a $r' <r$ on a une inclusion $j_{r,r'} : M^{^\tas}_{>r} \subset
M^{^\tas}_{>r'}$ qui est un \qis\ et induit un \qis\  $f_{r,r'} : S_r^{^\tas} \fl S_{r'}^{^\tas}$
puis, en se limitant aux termes $S^0,S^1,S^2$, le \ps\ $u_{r,r'}$ annonc\'e. Comme on a 
$j_{r,r'} j_{r',r''} = j_{r,r''}$ on en d\'eduit la m\^eme relation pour $f$ et $u$, \`a
homotopie pr\`es.

Passons au point 3).  Le complexe $\bR
\G_*
\cJ_{\cC}$  calcule la cohomologie de $\cJ_\cC$ de sorte que les complexes  
$M^{^\tas}$, $M^{^\tas}_{>r}$ et $S_r^{^\tas}$ calculent encore $H^0_* \cJ_\cC$. En
passant \`a $\Lpont$, c'est-\`a-dire en supprimant les termes de degr\'e $<0$ de
$S_r^{^\tas}$, on a encore une surjection 
$h_1 \Lpont =N
\fl H^0_* \cJ_\cC =I_\cC$. On en d\'eduit une suite exacte de faisceaux $0 \fl \cP \fl
\cN
\fl
\cJ_\cC
\fl 0$. Comme l'homologie de
$\Lpont$ en degr\'es $0$ et $-1$ est de type fini sur $A$ on a la suite exacte $0 \fl
\cN \fl \cL_1 \fl
\cL_0 \fl \cL_{-1} \fl 0$ qui montre que les groupes d'homologie de $\Lpont \T_A \tas$ sont
aussi les $H^i(\cN \T_A \tas)$. Il en r\'esulte qu'on a $H^i(\cP \T_A \tas)=0$ pour $i=1,2$ ce
qui prouve que $\cP$ est dissoci\'e en vertu de [H] 7.9. On a donc bien construit une
r\'esolution de type \sN\ triadique de $\cC$.

\remas {3.3} \lign
1) La m\'ethode de construction  de la r\'esolution de type \sE\ par d\'esaturation et la
deuxi\`eme m\'ethode de construction de la r\'esolution de type \sN\ sont  voisines car
 toutes deux font appel \`a un proc\'ed\'e de troncature.  C'est d'ailleurs ce qui distingue les
r\'esolutions obtenues ainsi des r\'esolutions (triadiques) de type \sE\ et
\sN\ g\'en\'erales. En effet, dans le cas g\'en\'eral on voit que
$H^1_* \cE$ est un quotient (de type fini sur $A$) du $R_A$-module $H^0_* \cJ_\cC$ et $H^2_* \cN$
un sous-$R_A$-module (de type fini sur $A$) de
$H^2_*
\cJ_\cC$, tandis que dans les constructions ci-dessus ce sont plus pr\'ecis\'ement des modules
tronqu\'es de la forme  $(H^0_* \cJ_\cC)_{>r}$ et $(H^2_* \cJ_\cC)_{\leq r}$.
\lign 2) Il y a encore une autre m\'ethode pour montrer l'existence des r\'esolutions de
type
\sN\  triadiques. Elle consiste \`a  faire le m\^eme travail que dans 2.12
ci-dessus mais de mani\`ere duale en d\'esaturant le module 
 $\Om_0= H^0_*(\bP_T,
\om)$ au lieu de d\'esaturer $I_\cC$.

\th {Proposition 3.4}. 
Soient $r_1$ et $r_2$  deux r\'esolutions de type \sE\ de $\cC$. Il existe une
r\'esolution  $r'$ de type \sE\
 cotriadique avec des morphismes  $r'
\fl r_1$ et
$r'\fl r_2$.

\dem Posons, pour $k=1,2$, $(r_k) = (0 \fl \cE_k \fl \cF_k \fl \cJ_\cC \fl 0)$ et soit
$I_k$ l'image du module libre $F_k$ dans $H^0_* \cJ_\cC$ de sorte que  $I_1$ et $I_2$
sont contenus dans $I_\cC$. 

On applique encore la m\'ethode du lemme de d\'esaturation : on choisit un $J =
\bigoplus_{n
\geq n_0} H^0 \cJ_\cC (n)$ avec $n_0$ assez grand pour que le foncteur $H^1$ soit
nul  et pour que $J$ soit contenu dans les $I_k$ (il
suffit que l'on ait
$H^1 \cE_k(n) =0$ pour $n \geq n_0$), puis on prend le d\'ebut d'une r\'esolution
minimale
$F'
\fl J
\fl 0$ et la suite de faisceaux associ\'ee \`a $0 \fl E' \fl F' \fl J \fl 0$ convient.

\th {Proposition 3.5}.  Soient $r_1$ et $r_2$  
deux r\'esolutions de type \sN\ de $\cC$. Il existe une
r\'esolution  $r'$ de type \sN\ 
triadique avec des morphismes  $r_1 \fl r'$ et $r_2 \fl r'$.

\dem On effectue  une liaison par une intersection compl\`ete
$\cD$ contenant $\cC$ qui donne une famille $\cC'$ de courbes (cf. [HMDP1] 1.4). Les
r\'esolutions de type
\sN\ de $\cC$ donnent des r\'esolutions de type \sE\ de $\cC'$ (cf. [HMDP1] 2.24). On
applique alors 3.4 \`a ces r\'esolutions et on obtient une r\'esolution cotriadique de
$\cC'$ qui  domine les autres, puis, par liaison (cf. [HMDP1] 2.25), la r\'esolution de
type
\sN\ triadique de
$\cC$ cherch\'ee.

\tarte {b) Triade associ\'ee \`a une famille de courbes}

La proposition suivante permet d'associer une triade  \`a une famille de
courbes.  On notera, par rapport au cas d'un corps, que la triade en question n'est d\'efinie ici
qu'\`a pseudo-isomorphisme pr\`es.

\th {Proposition-d\'efinition 3.6}. Soit $\cC$ une famille de courbes de $\bP^3_A$ munie d'une
r\'esolution de type \sN\ triadique : $0 \fl \cP \fl \cN \fl \cJ_\cC \fl 0$ et soit 
$$0 \fl \cN \fl  \cL_1 \fl \cL_0 \fl \cL_{-1} \fl 0 \leqno{(1)}$$
une suite exacte associ\'ee au faisceau $\cN$ (cf. 2.1). Alors, si on pose
$L_i= H^0_* (\cL_i)$ le complexe $\Lpont =(L_1 \fl L_0 \fl L_{-1})$ d\'eduit de
(1) est une triade majeure qui est appel\'ee une {\bf triade
 de Rao} de $\cC$. Elle est bien d\'efinie
\`a  \ps\ pr\`es. \lign
 Le foncteur associ\'e
\`a toute triade de Rao de $\cC$ est le foncteur $H^1_* (\cJ_\cC \T_A \tas)$ : il est appel\'e 
{\bf foncteur de Rao} de $\cC$. \lign
On note
$\cT(\cC)$ la {\bf classe de triades} (pour la relation de pseudo-isomorphisme)
qui contient les triades de Rao de $\cC$.  

\dem Il est clair que $\Lpont$ est une triade majeure. Le fait que deux triades de Rao sont
pseudo-isomorphes
provient de [HMDP1] 2.18.3 et de 2.14 ci-dessus. L'assertion sur le
foncteur associ\'e r\'esulte de l'isomorphisme 
$H^1_* (\cN \T_A \tas) \simeq H^1_* (\cJ_\cC \T_A \tas)$ et de 2.4.2.

\vskip 0.3 cm

On peut  d\'efinir aussi par dualit\'e une ``cotriade'' de Rao de $\cC$ :

\th {Proposition-d\'efinition 3.7}. Soit $\cC$ une famille de courbes de $\bP^3_A$ munie d'une
r\'esolution de type \sE\ cotriadique : $0 \fl \cE \fl \cF \fl \cJ_\cC \fl 0$ et soit 
$$0 \fl \cM_1 \fl \cM_0 \fl \cM_{-1} \fl \cE \fl 0 \leqno{(2)}$$
une suite exacte associ\'ee au faisceau $\cE$ (cf. 2.1). Alors, si on pose 
$M_i = H^3_* \cM_i$, le complexe $\Mpont = (M_1 \fl M_0 \fl M_{-1})$ est  appel\'e une
cotriade de Rao de $\cC$. Elle est bien d\'efinie au dual d'un \ps\
 pr\`es.

Le th\'eor\`eme suivant explicite le lien entre triade et cotriade associ\'ees \`a une famille
de courbes :

\th {Th\'eor\`eme 3.8}. Soit $\cC$ une famille de courbes, $\Lpont$ (resp. $\Mpont$) une triade
(resp. une cotriade) de Rao de $\cC$. Alors, il existe un \ps\ fort (de complexes) $\f :
\Lpont
\fl
\Mpont$.

\dem On utilise les notations de \S\ 0 {\it d)}. On pose $X= \bP^3_A$. Consid\'erons des
r\'esolutions de type \sE\ et
\sN\ de $\cC$ comme en 3.6 et 3.7 ci-dessus avec les suites exactes (1) et (2). La suite
(1) montre que le complexe
$\cL_{\tas} = (\cL_1 \fl
\cL_0
\fl
\cL_{-1})$ (que nous consid\'ererons comme un complexe en degr\'es $-1, 0, 1$ en
inversant les indices) est
\'egal
\`a
$\cN[1]$ dans la cat\'egorie d\'eriv\'ee
$\cD(X)$.  Si
$ \cN \fl
\cI_{\tas}$ est une r\'esolution injective on obtient un morphisme $\cL_{\tas}[-1] \fl 
\cI_\tas$
qui est un \qis.  Cela fournit, par application du foncteur $\G_*$ un morphisme 
$$ \a : \Lpont [-1] = \G_* \cL_\tas [-1] \fl \G_* \cI_\tas = \bR \G_* \cN$$
dans la cat\'egorie d\'eriv\'ee $\cD(R_A)$.

   Par ailleurs, il r\'esulte
de [HMDP1] 2.17.3 qu'il existe un faisceau dissoci\'e
$\cL$ avec une suite exacte 
$0 \fl \cE \fl \cL \fl \cN \fl 0$. Le triangle associ\'e \`a cette suite exacte donne un
morphisme $\beta_0 : \cN \fl \cE [1]$ dans la cat\'egorie d\'eriv\'ee $\cD (X)$ 
puis, par
application de $\G_*$, un morphisme dans $\cD (R_A)$ :
$$\beta : \bR \G_* \cN \fl \bR \G_* \cE [1]$$
 qui induit, puisque $\cL$ est dissoci\'e, un \'epimorphisme de foncteurs
$H^0_* (\cN \T_A \tas) \fl H^1_* (\cE \T_A \tas)$, un isomorphisme
$H^1_* (\cN \T_A \tas) \fl H^2_* (\cE \T_A \tas)$, un monomorphisme
$H^2_* (\cN \T_A \tas) \fl H^3_* (\cE \T_A \tas)$.

On a, ensuite, l'isomorphisme de dualit\'e vu en 0.5 :
$$\theta : \bR \G_* \cE [1] \fl \bR\Homgr_{A} (\bR\G_* \cE^\vee (-4) , A) [-2].$$

De la m\^eme fa\c con que ci-dessus, et avec les m\^emes conventions de degr\'es, la
suite duale de (2) donne un morphisme
$$\g_0 : \G_* \cM_{\tas}^\vee (-4) [-1] \fl \bR\G_* \cE^\vee (-4),$$ ou encore, en
tenant compte de la dualit\'e de Serre : $H^0_* \cM_i^\vee(-4) \simeq (H^3_* \cM_i)^*$ :
$$\g_0 : \Mpont ^* [-1] \fl \bR\G_* \cE^\vee (-4),$$
 et, par
application de
$\bR\Homgr_{A}(\tas, A)$, on obtient
$$\g : \bR\Homgr_{A} (\bR\G_* \cE^\vee (-4), A) [-2] \fl \bR\Homgr_{A}
(\Mpont^*[-1] , A) [-2] =
\Mpont [-1]$$
la derni\`ere \'egalit\'e provenant du fait que $\Mpont^*$ est form\'e de $R_A$-modules libres,
donc est acyclique pour $\Homgr$.

En d\'efinitive, en d\'ecalant de $1$,  on a  un morphisme $\f =\g \theta \beta \a :
\Lpont \fl \Mpont$ dans $\cD(R_A)$. {\it A priori}, ce morphisme est dans la cat\'egorie
d\'eriv\'ee mais, comme $\Lpont$ est un complexe de $R_A$-modules libres de type fini,
$\f$ est repr\'esent\'e par un morphisme de complexes, unique \`a homotopie pr\`es (cf.
[RD] I 4.7 p. 46 en changeant injectif en projectif et en renversant le sens des
fl\`eches). De plus, ce morphisme est un \ps\ fort. En effet, cela r\'esulte du fait que
$\theta$ est un isomorphisme, que
$\a$ et $\g$ induisent des isomorphismes sur les foncteurs $h_1, h_0$ et $h_{-1}$ et
des propri\'et\'es de
$\beta$ \'enonc\'ees ci-dessus.

\tarte {c) Le th\'eor\`eme de Rao pour les triades : fibres de $\Psi_A$}

 Un anneau $A$ \'etant donn\'e, les r\'esultats du paragraphe
pr\'ec\'edent permettent de d\'efinir une application
$\Psi_A$ qui
\`a une famille de courbes $\cC$ associe la classe de triades $\cT(\cC)$. Les
th\'eor\`emes ``de Rao''    vont permettre de pr\'eciser l'image et
les fibres de
$\Psi_A$ lorsque l'anneau $A$ est local et son corps r\'esiduel  infini.

\vskip 0.3 cm

Dans [HMDP1] 3.2 nous avons montr\'e le th\'eor\`eme suivant (g\'en\'eralisation aux
familles de courbes du th\'eor\`eme de Rao, deuxi\`eme forme) :

\th {Th\'eor\`eme de Rao pour les faisceaux}. On suppose $A$ local \`a corps r\'esiduel
infini. Soient $\cC$ et
$\cC'$ deux familles plates de courbes param\'etr\'ees par  $A$, munies de r\'esolutions
de type
\sN, avec des faisceaux $\cN$,
$\cN'$. Alors, $\cC$ et
$\cC'$ sont dans la m\^eme classe de biliaison si et seulement si $\cN$ et $\cN'$
sont pseudo-isomorphes, \`a d\'ecalage pr\`es.

En termes de triades, ce r\'esultat permet de calculer les fibres de $\Psi_A$, ce qui
g\'en\'eralise la premi\`ere forme du th\'eor\`eme de Rao :

\th {Th\'eor\`eme 3.9}. (Th\'eor\`eme de Rao pour les triades)  On suppose $A$ local
\`a corps r\'esiduel infini.  Soient
$\cC$ et
$\cC'$ deux familles plates de courbes param\'etr\'ees par  $A$ et soient $\Lpont$ et
$\Lpont'$ des triades de Rao de $\cC$. Alors,
$\cC$ et
$\cC'$ sont dans la m\^eme classe de biliaison si et seulement si les triades $\Lpont$ et
$\Lpont'$  sont pseudo-isomorphes,  \`a d\'ecalage pr\`es (i.e. s'il existe un entier
$h$ tel que l'on ait
$\cT(\cC) = \cT(\cC') (h)$).

\dem Cela r\'esulte du th\'eor\`eme de Rao rappel\'e ci-dessus, de 3.6 et de 2.14.

\vskip 0.3 cm

\tarte {d) Le th\'eor\`eme de Rao pour les triades : image de $\Psi_A$}

Gr\^ace aux r\'esultats de [HMDP2]  exprim\'es en termes de triades on
peut montrer la g\'en\'eralisation de l'assertion de surjectivit\'e (toujours \`a
d\'ecalage pr\`es) du th\'eor\`eme de Rao (i.e. de l'application $\Psi_A$) :

\th {Th\'eor\`eme 3.10}. On suppose $A$ local \`a corps r\'esiduel infini.  Soit
$\Lpont$ une triade. Il existe une famille de courbes
$\cC$ telle que la triade associ\'ee \`a $\cC$ soit pseudo-isomorphe \`a $\Lpont$, \`a
d\'ecalage pr\`es.

\dem On peut supposer que $\Lpont$ est une triade majeure. On prend le faisceau
triadique associ\'e $\cN$ (cf. 2.5). En vertu de [HMDP2] 2.7 il existe une courbe $\cC$
avec une r\'esolution de type \sN\ de la forme $0 \fl \cP \fl \cN \fl \cJ_\cC \fl 0$ et
cette courbe convient.

\vskip 0.3 cm

Plus pr\'ecis\'ement, le th\'eor\`eme suivant donne exactement l'image de $\Psi_A$ i.e.
les d\'ecalages possibles :

\th {Th\'eor\`eme 3.11}. On suppose $A$ local \`a corps r\'esiduel infini. Soit $\Lpont$
une triade, $\cT$ sa classe de pseudo-isomorphisme, $\cN$ un faisceau triadique
associ\'e et soit
$\cN_0$ le faisceau extraverti minimal (unique) de la classe de pseudo-isomorphisme de
$\cN$ (cf. [HMDP1] 2.14). 
Soit  $q=
q_{\cN_0}$ la fonction $q$ de $\cN_0$ (cf. [HMDP2] 2.4) et soit  $h_0= \sum_{n \in \bZ}
nq(n) +
\deg
\cN_0$. \lign 1) Il
existe une famille de courbes 
$\cC_0$ et une r\'esolution 
$ 0 \fl \cP_0 \fl \cN_0 \fl \cJ_{\cC_0} (h_0) \fl 0$ avec $\cP_0$ dissoci\'e. La
classe de triade associ\'ee \`a $\cC_0$ est \'egale \`a $\cT(-h_0)$.  \lign 
2) Si $\cC_1$ est
une famille de courbes telle que $\cT(\cC_1) = \cT (-h)$,   on a
$h
\geq h_0$.  Si
$d$ et
$g$ (resp.
$d_0$ et
$g_0$) sont respectivement le degr\'e et le genre de $\cC_1$ (resp. $\cC_0$) on a $d \geq
d_0$ et
$g
\geq g_0$.\lign
3) R\'eciproquement, pour tout $h \geq h_0$ il existe une famille de courbes $\cC_1$ avec
$\cT(\cC_1) = \cT(-h)$. \lign 
4) Si de plus on a
$h=h_0$,
$\cC_1$ et $\cC_0$ sont jointes par une d\'eformation  \`a cohomologie uniforme 
(cf. [HMDP2] 2.8) et triade constante. 
\lign On dit que
$\cC_0$ est une {\bf famille minimale} de courbes. \lign
Les autres familles  de
courbes de la classe de biliaison s'obtiennent \`a partir de $\cC_0$ par des biliaisons
\'el\'ementaires  suivies d'une d\'eformation \`a cohomologie uniforme et triade
constante.

\dem Cela r\'esulte de [HMDP2] 2.9 et 2.10 (Th.  de Lazarsfeld-Rao). Dans  le point 4),
comme les familles $\cC_0$ et $\cC_1$ ont des r\'esolutions de type \sN\ avec les 
m\^emes faisceaux $\cP$ et $\cN$ les triades associ\'ees sont les m\^emes \`a \ps\
pr\`es.

\remas {3.12} \lign
1)  Ce th\'eor\`eme d\'ecrit exactement les familles de courbes minimales
associ\'ees aux triades. Il g\'en\'eralise donc les r\'esultats de [MDP1,3] qui
traitaient ce probl\`eme pour les courbes sur un corps, \`a partir des modules de Rao.
On notera que si deux triades sont pseudo-isomorphes elles
correspondent \`a la m\^eme    classe de biliaison de familles de courbes et, en
particulier, ont m\^eme famille minimale, \`a d\'eformation pr\`es.
\lign
2) Lorsque $A$ est un \avd\  le faisceau $\cN_0$ extraverti minimal est triadique et
on peut le calculer explicitement \`a partir d'une  triade donn\'ee (cf. 1.32 et les
exemples du \S 5). Il est caract\'eris\'e par le fait que
$H^2_*
\cN$ est de torsion. Dans le cas d'un anneau local quelconque la d\'etermination du
faisceau extraverti minimal est plus probl\'ematique.

\tarte {e)  Le cas de la liaison impaire}

On a, gr\^ace \`a la dualit\'e, une caract\'erisation de la liaison impaire :

\th {Th\'eor\`eme  3.13}. On suppose $A$ local \`a corps r\'esiduel infini. Soient $\cC$
et $\cC'$ deux familles plates de courbes param\'etr\'ees par  $A$, et soient $\Lpont$ et
$\Lpont'$ des triades de Rao de $\cC$ et $\cC'$. Alors, $\cC$ et
$\cC'$ sont li\'ees par un nombre impair de liaisons \'el\'ementaires si et
seulement si les triades
$\Lpont$ et
$\Lpont'$ sont duales \`a \ps\ pr\`es et \`a d\'ecalage pr\`es,
 i.e. s'il existe un entier $h$ tel que $\Lpont
(h)$ soit pseudo-isomorphe \`a  une triade duale de $\Lpont'$.

\dem Supposons $\cC$ et $\cC'$ jointes par un nombre impair de liaisons et montrons que les
triades associ\'ees sont duales, \`a d\'ecalage pr\`es.  En vertu de 3.9 il suffit de traiter le
cas d'une liaison par des surfaces de degr\'es
$s$ et
$t$. On consid\`ere des r\'esolutions de type
\sN\ et
\sE\ de $\cC$ et de type \sN\ de $\cC'$, avec les faisceaux respectifs $\cN, \cE, \cN'$. On sait,
 cf. [HMDP1] 2.24 qu'on peut prendre
$\cE =
\cN'^\vee (-s-t)$. On
note $\Lpont$ et $\Lpont'$ les triades associ\'ees \`a $\cN$ et $\cN'$ et $\Mpont$ la cotriade
associ\'ee \`a $\cE$.  En vertu de  0.2 on a 
$\Mpont = (\Lpont'(s+t-4))^*$ au sens des complexes. Le th\'eor\`eme de comparaison 3.8
ci-dessus fournit un \ps\ fort
$\f :
\Lpont
\fl \Mpont$ et on conclut par 1.27.4.

R\'eciproquement, si les triades sont duales, on lie $\cC$ \`a une famille $\cC_1$. Le sens
direct montre que les triades associ\'ees \`a $\cC_1$ et $\cC'$ sont pseudo-isomorphes et on
conclut avec 3.9 appliqu\'e \`a $\cC_1$ et $\cC'$.

\rema {3.14} Si le corps r\'esiduel $k$ de $A$ est fini les assertions d'existence 3.9, 3.10,
3.11 et 3.13 restent vraises \`a condition de remplacer \'eventuellement $A$ par un anneau local
$B$ fini et
\'etale sur $A$.

\titre {4. Courbes, faisceaux et triades : dictionnaire}

\tarte {a) Propri\'et\'es des triades et des foncteurs triadiques}

Les propositions suivantes  caract\'erisent deux types
importants de triades :

\th {Proposition 4.1}. Soit $\Lpont = (L_1 \Fl {d_1} L_0 \Fl {d_0} L_{-1})$ une triade,
$H$ son c\oe ur, $C$ son conoyau et $V$ le foncteur associ\'e. Les propri\'et\'es
suivantes sont
\'equivalentes :\lign 
1) $V$ est exact \`a droite, \lign
2) on a un isomorphisme de foncteurs
$V \simeq H\T_A \tas$,\lign
3) $C$ est plat sur $A$, \lign
 4)  $\Lpont$ est pseudo-isomorphe \`a une triade de la forme $L_1 \Fl {d_1} L_0
\Fl {d_0} 0 $. \lign
5) $\Lpont$ est pseudo-isomorphe \`a  une triade majeure (resp.  mineure) de la forme
$L_1 \Fl {d_1} L_0
\Fl {d_0} 0 $. 
On dit alors que la triade est exacte \`a droite ou {\bf modulaire}.

\dem
 $1) \impl 2)$. 
Soit $Q$
 un $A$-module quelconque, on le r\'esout par des $A$-modules libres :
$  F_1 \fl F_0 \fl Q \fl 0$
et, comme $V$ est exact \`a droite, on en d\'eduit la suite exacte $V(F_1) \fl V(F_0) \fl 
V(Q)
\fl 0$. Mais, comme
$V$ est additif on a $V(F_i) = V(A) \T_A F_i$ d'o\`u la conclusion puisque l'on
a $H=V(A)$.

$2) \equi 3)$. La platitude de $C$ est \'equivalente \`a la nullit\'e de
$\Tor^1_A(C,Q)$ pout tout $Q$ et la conclusion r\'esulte du lemme 1.5.

 $ 3) \impl 4)$. En vertu de 1.30 et 1.16 on se ram\`ene au cas o\`u $C$ est nul. Il
suffit alors de remplacer $\Lpont$ par $L_1 \Fl{d_1} \Ker d_0$ (on notera que, comme
$d_0 : L_0
\fl L_{-1}$ est surjective et $L_{-1}$ plat sur $A$, il en est de m\^eme de $\Ker d_0$).

$4) \impl 5)$. Cela vient de 1.14 et 1.24.

$5) \impl 2)$ est clair et $2) \impl 1)$  r\'esulte de l'exactitude \`a droite du
produit tensoriel.

\th {Proposition 4.2}. Soit $\Lpont = (L_1 \Fl {d_1} L_0 \Fl {d_0} L_{-1})$ une triade,
$H$ son c\oe ur et $V$ le foncteur associ\'e. Les propri\'et\'es
suivantes sont \'equivalentes :
\lign 1) $V$ est exact \`a gauche, \lign
2) il existe un $R_A$-module gradu\'e $D$, de type fini sur $A$, tel que l'on
ait un isomorphisme de foncteurs
$V \simeq \Homgr_{A}(D,\tas)$,\lign
3)  $E=\Coker d_1$ est plat sur $A$, \lign
3') il existe une triade $\Lpont'= (L'_1 \Fl {d'_1} L'_0 \Fl {d'_0}
L'_{-1})$ pseudo-isomorphe \`a $\Lpont$ telle que 
 $E'=\Coker d'_1$ soit plat sur $A$, \lign
4) 
 $\Ker d_1$ commute au changement de base, i.e., pour tout $A$-module de type
fini
$Q$, la fl\`eche naturelle $$(\Ker d_1) \T_A Q \fl \Ker (d_1 \T_A Q) \qquad
\hbox {est un isomorphisme},$$
 5)   $\Lpont $ est pseudo-isomorphe \`a une triade de la forme $0 \Fl {d_1} L_0
\Fl {d_0} L_{-1} $, \lign
6)$\Lpont $ est pseudo-isomorphe \`a une triade mineure $\Ppont $ de la forme $0 \Fl
{d_1} P_0
\Fl {d_0} P_{-1}$. \lign
On dit alors que la triade est exacte \`a gauche ou {\bf repr\'esentable}. \lign
Si, de plus, $A$ est un anneau principal on peut remplacer dans 3) la platitude de
$\Coker d_1$ par celle de $H= h_0 (\Lpont)$.

\dem
$1) \impl 2)$. On consid\`ere une r\'esolution injective de $Q$ : $0 \fl Q \fl I_0 \fl
I_1$. Comme le foncteur $V$ est exact \`a gauche on obtient la suite exacte $0 \fl V(Q)
\fl V(I_0) \fl V(I_1)$. Posons $D = h_0 (\Lpont^*)$ et $C' = h_{-1} (\Lpont^*)$. On montre
alors, comme pour 1.5, mais en renversant le sens des fl\`eches, qu'on a, pour tout
$A$-module
$M$, la suite exacte (cf. aussi [H] 4.5)
$$0 \fl \Ext^1_A(C',M) \fl V(M) \fl \Homgr_{A}(D,M) \fl \Ext^2_A(C',M).$$
En appliquant cette suite \`a $M= I_k$ on voit que l'on a $V(I_k) = \Homgr_{A} (D,I_k)$
(car
$I_k$ est injectif de sorte que les $\Ext^i$ sont nuls) et on en d\'eduit $V(Q) =
\Homgr_{A} (D,Q)$.

$2) \impl 1)$ est clair.

Soit $\Lpont$ un complexe de $\cT$. On a vu, cf. 1.4, que, si $d : E \fl L_{-1}$ est la
fl\`eche induite par $d_0$, on a $V(Q) = \Ker d\T Q$. Si on a une
injection
$Q'
\fl Q$ de
$A$-modules l'injectivit\'e de
$V(Q')
\fl V(Q)$ est alors
\'equivalente \`a celle de $E \T_A Q' \fl E\T_A Q$, donc \`a la platitude de $E$. Ceci
montre $1) \impl 3)$ et $3' \impl 1)$.  Bien entendu $3) \impl 3')$ est clair.

Montrons  l'\'equivalence de 3) et 4). Si on pose $B= \Im d_1$, la platitude de
$E$ implique  celle de
$B$, donc \'equivaut aux relations ${\rm Tor}_A^1\, (E,Q)= {\rm Tor}_A^1\, (B,Q) =0$ pour
tout
$A$-module  $Q$ et l'application de
$\tas
\T_A Q$ \`a la suite exacte $0 \fl \Ker d_1 \fl L_1 \fl B \fl 0$ montre  que ces
conditions
\'equivalent aussi \`a la commutation du noyau  au produit tensoriel.

$3) \impl 5)$ est clair (il suffit de prendre le complexe $0 \fl E \Fl {d} L_{-1}$ qui est
ad\'equat puisque $E$ est plat).

$5) \impl 6)$ se d\'emontre en tronquant les complexes comme en 1.24 et $6)
\impl 3')$ est clair.

\vskip 0.2 cm

 Pour l'assertion suppl\'ementaire, notons qu'on a la suite exacte $0 \fl H \fl E \fl \Im
d_0
\fl 0$. Si l'anneau est principal, le module $\Im d_0$ qui est inclus dans $L_{-1}$ est sans
torsion, donc plat, de sorte que la platitude de $E$ et celle de $H$ sont \'equivalentes.

\th {D\'efinition 4.3}. Une triade est dite exacte si elle est \`a la fois exacte
\`a gauche et \`a droite. Une telle triade est modulaire et d\'efinie par un
$R_A$-module
$H$ plat sur $A$.

\remas {4.4} \lign
1) Comme les propri\'et\'es 4.1 et 4.2 s'expriment en termes de foncteurs $V$, il est
clair que si une triade est pseudo-isomorphe
\`a une triade modulaire (resp. repr\'esentable, resp. exacte) elle est modulaire (resp.
repr\'esentable, resp. exacte). Si $\cT$ est une classe de triades pour la relation de
pseudo-isomorphisme on dira que $\cT$ est modulaire (resp. repr\'esentable, resp.
exacte) si l'une quelconque des triades de $\cT$ l'est.
\lign   2) Si $\Lpont$ et $\Mpont$ sont deux triades duales, $\Lpont$ est modulaire  si et
seulement si $\Mpont$ est repr\'esentable (utiliser les caract\'erisations 4)
 et 6) de 4.1 et 4.2).

\vskip 0.3 cm

La proposition suivante montre que, dans le cas des triades modulaires 
ou repr\'esenta\-bles,
le foncteur d\'etermine la triade, si bien que la difficult\'e
\'evoqu\'ee en 1.35.{\it c} ne se produit pas.

\th {Proposition 4.5}. Soient $\Lpont$ et $\Lpont'$ deux triades modulaires (resp.
repr\'esentables) et $V$ et $V'$ les foncteurs associ\'es. On suppose $V$ et $V'$
isomorphes. Alors $\Lpont$ et $\Lpont'$  sont pseudo-isomorphes.

\dem Dans le cas modulaire, posons $H= V(A)$ et $H'= V'(A)$. Ces modules sont isomorphes et
$V$ et
$V'$ sont tous deux isomorphes au foncteur  $H \T_A \tas$. On peut supposer que 
$\Lpont = (L_1
\Fl{d_1} L_0
\fl 0)$ et $\Lpont' = (L'_1 \Fl{d'_1} L'_0 \fl 0)$  sont deux triades majeures. Les
fl\`eches $d_1$ et $d'_1$ donnent des pr\'esentations libres de $H$. On a alors une
fl\`eche $\Lpont \fl \Lpont'$ qui est un pseudo-isomorphisme, d'o\`u le r\'esultat.

\vskip 0.3 cm

Le cas repr\'esentable r\'esulte du cas modulaire par dualit\'e (cf. 4.4.2 et 1.27).

\tarte {b) Le dictionnaire courbes-faisceaux-triades}

L'\'enonc\'e suivant \'etablit le lien entre les propri\'et\'es des familles de
courbes (\`a sp\'ecialit\'e constante, etc.), pour lesquelles on renvoie \`a [MDP1] VI,
 les propri\'et\'es des triades associ\'ees (exactitude \`a droite, etc.), et celles
des faisceaux $\cE$ et $\cN$ des r\'esolutions de type \sE\ et \sN.

\th {Proposition 4.6}. Soit $\cC$ une famille de courbes, $\Lpont$ une
triade associ\'e. Soient $0 \fl \cE \fl  \cF \fl \cJ_\cC \fl 0$  une r\'esolution de
type \sE\  et
$0 \fl \cP \fl  \cN \fl \cJ_\cC \fl 0$ une r\'esolution de type \sN\
triadique  de $\cJ_\cC$.  Les propri\'et\'es suivantes sont
\'equivalentes~:
\lign
1) $\cC$ est \`a sp\'ecialit\'e constante, \lign
1') $H^2_* \cJ_\cC$ est plat sur $A$, \lign
1'') $H^1_* \cJ_\cC$ commute au changement de base, \lign
2)  $\Lpont$ est une
triade exacte \`a droite, \lign
3) $H^2_* \cN$ est plat sur $A$, \lign
3') $H^1_* \cN$ commute au changement de base, \lign
4) $H^2_* \cE$ commute au changement de base, \lign
4') $H^3_* \cE$ est plat sur $A$. \lign
De plus, si $A$ est un \avd\ et si le faisceau $\cN$ est extraverti, la condition 3)
est \'equivalente \`a $H^2_* \cN=0$.

\dem L'\'equivalence de {\sl 1)} et {\sl 1')} est claire  par d\'efinition de la
sp\'ecialit\'e constante. 
 
Pour toute famille de courbes
$\cC$, comme
$H^3_*
\cJ_\cC$ est plat (c'est
$H^3_*
\cO_\bP$) et commute au changement de base, $H^2_* \cJ_\cC$ commute au changement 
de base de sorte qu'on a l'\'equivalence de {\sl 1')} et {\sl 1'')}. Comme le foncteur
$V$ de la triade associ\'ee \`a $\cC$ n'est autre que $H^1_* (\cJ_\cC \T_A \tas)
$ (cf. 3.6), l'\'equivalence de {\sl 1'')} et {\sl 2)} est claire.

On a un isomorphisme de foncteurs $H^1_* (\cN \T_A \tas) \simeq H^1_* (\cJ_\cC \T_A
\tas)$ d'o\`u l'\'equivalence de {\sl 3')} et {\sl 1'')}. Comme $\cN$ est triadique,
$H^2_* \cN$ commute au changement de base, et cette propri\'et\'e est encore
\'equivalente
\`a {\sl 3)}.

Enfin l'\'equivalence de {\sl 1'')} et {\sl 4)} est claire et {\sl 4')} en r\'esulte
aussit\^ot (car un $H^3$ commute toujours au changement de base).

Dans le cas extraverti, comme $H^2_* \cN$ est de torsion, il est plat si et seulement
si il est nul.

\th {Proposition 4.7}. Soit $\cC$ une famille de courbes, $\Lpont$ une
triade associ\'ee. Soient $0 \fl \cE \fl  \cF \fl \cJ_\cC \fl 0$ une r\'esolution de
type \sE\ cotriadique  et
$0 \fl \cP \fl  \cN \fl \cJ_\cC \fl 0$ 
une r\'esolution de type
\sN\  de $\cJ_\cC$.  Les propri\'et\'es suivantes sont
\'equivalentes~:
\lign
1) $\cC$ est \`a postulation constante \lign
1') $H^0_* \cJ_\cC$ commute au changement de base, \lign
1'') le foncteur $H^1_* (\cJ_\cC \T_A \tas)$ est exact \`a gauche, \lign
2)  $\Lpont$ est une
triade exacte \`a gauche, \lign
3) $H^0_* \cN$ commute au changement de base, \lign
4) $H^1_* \cE$ commute au changement de base.

\dem L'\'equivalence de {\sl 1)} et {\sl 1')} r\'esulte de la d\'efinition de la
postulation cons\-tante. On voit aussit\^ot en tensorisant une suite exacte $0 \fl Q' \fl
Q \fl Q'' \fl 0$, que cette condition \'equivaut \`a {\sl 1'')} et donc \`a {\sl 2)}
par d\'efinition de la triade associ\'ee (cf. 3.6).

 Pour {\sl 3)}  on a la
suite exacte de foncteurs :
$$ 0 \fl H^0_* (\cP \T_A \tas) \fl H^0_* (\cN \T_A \tas) \fl H^0_* (\cJ_\cC \T_A \tas)
\fl 0 $$ et, comme $H^0_* (\cP \T_A \tas)$ est exact, on voit  facilement que {\sl 1')}
et {\sl 3)} sont \'equivalentes.

 Pour {\sl 4)}  on a la
suite exacte de foncteurs :
$$ 0 \fl H^0_* (\cE \T_A \tas) \fl H^0_* (\cF \T_A \tas) \fl H^0_* (\cJ_\cC \T_A \tas)
\Fl{\f}  H^1_* (\cE \T_A \tas)  \fl 0 $$ 
avec  $H^0_* (\cE \T_A \tas)$ exact (car $\cF$ est dissoci\'e) et, comme $\cE$ est
dual d'un triadique, 
$ H^0_* (\cE
\T_A \tas) $  exact \`a droite et $ H^1_* (\cE \T_A \tas) $ exact \`a gauche.  On coupe
la suite exacte en deux en introduisant $F= \Ker \f$. La premi\`ere suite exacte montre
que $F$ est exact \`a droite et on conclut \`a l'\'equivalence de {\sl 1')} et de {\sl
4)} gr\^ace \`a la deuxi\`eme suite exacte.

\th {Corollaire 4.8}.  Soit $\cC$ une famille de courbes, $\Lpont$ une
triade associ\'ee. Soient $0 \fl \cE \fl  \cF \fl \cJ_\cC \fl 0$  une r\'esolution de
type \sE\ cotriadique  et
$0 \fl \cP \fl  \cN \fl \cJ_\cC \fl 0$ une r\'esolution de type
\sN\  de $\cJ_\cC$.  Les propri\'et\'es suivantes sont
\'equivalentes~:
\lign
1) $\cC$ est \`a cohomologie constante \lign
2)  $\Lpont$ est une
triade exacte, \lign
3) $H^0_* \cN$ est plat et commute au changement de base, \lign
4) $H^2_* \cE$ est plat et commute au changement de base.

\vfill\eject

\titre {5. Construction de triades, application \`a la construction de familles de
courbes}

\vskip 0.3 cm

 Ce paragraphe, qui est essentiellement constitu\'e d'exemples, a pour 
objectif essentiel de montrer comment on peut, en utilisant les triades, construire des familles de
courbes gauches, param\'etr\'ees par un \avd.\footnote {$(^5)$} {Une autre question
fondamentale, que nous aborderons dans une \'etude ult\'erieure, est de montrer, au contraire,
l'inexistence de certaines familles de courbes.}

\vskip 0.3 cm

Un bon fil conducteur pour le
parcourir est l'exemple des courbes de degr\'e $4$ et genre $0$ \'evoqu\'e dans
l'introduction. Dans ce cas, on sait, cf. [MDP4] ou [El], que le sch\'ema de Hilbert
$H_{4,0}$ est r\'eunion de deux sous-sch\'emas \`a cohomologie constante $H_1= H_{\g_1,
\rho_1}$ et $H_2= H_{\g_2,
\rho_2}$. Ces  sous-sch\'emas sont tous deux  de dimension
$16$ et leurs adh\'erences sont les composantes irr\'eductibles de $H_{4,0}$.
 La  courbe g\'en\'erale $C_1$ de $H_1$ est une courbe de bidegr\'e $(3,1)$ sur une
quadrique, de module de Rao $k(-1)$. La courbe g\'en\'erale $C_2$ de $H_2$ est
la r\'eunion disjointe d'une cubique plane et d'une droite, avec un module de Rao
 du type de $R/(X,Y,Z,T^3)$, de dimensions $1,1,1,$ en degr\'es $0,1,2$. Par semi-continuit\'e,
il est clair que $H_1$ est ouvert et $H_2$ ferm\'e et  la question est de savoir si $H_2$ est
faiblement adh\'erent \`a
$H_1$, i.e. si on a 
$\ov {H_1}
\cap H_2
\neq \vide$, ou encore, s'il existe une famille de courbes $\cC$, param\'etr\'ee par un
\avd\ $A$ dont le point g\'en\'erique est dans $H_1$ et le point sp\'ecial dans $H_2$.
Nous allons voir que c'est bien le cas (cf. 5.21). 

\vskip 0.3 cm

On sait, en vertu de 3.6, qu'\`a  toute famille de courbes $\cC$ est associ\'ee une
triade, dont le foncteur associ\'e mesure la variation du module de Rao dans $\cC$ et le
th\'eor\`eme de semi-continuit\'e affirme que les dimensions du module de Rao au
point g\'en\'erique sont inf\'erieures \`a celles  du module de Rao au point sp\'ecial.
La proposition 5.8  ci-dessous  pr\'ecise la relation entre ces modules en montrant
que, pour toute triade sur un
\avd, la valeur du foncteur
$V$ associ\'e
 au point g\'en\'erique est,
\`a d\'eformation pr\`es, un sous-quotient de la valeur de $V$ au point ferm\'e.

Cette constatation \'etant faite, notre objectif est de parcourir le chemin en sens
inverse.

Il s'agit d'abord, partant d'un  sous-quotient, de construire
une triade correspondant \`a ce sous-quotient. En g\'en\'eral la construction se fait en
deux \'etapes : on d\'etermine le c\oe ur $H$ et le conoyau $C$ de la triade
souhait\'ee, puis on construit une extension $u$ de longueur $2$ de $C$ par $H$. On reconstitue
alors ais\'ement la triade
\`a partir de ces
\'el\'ements, cf. \S a).  Dans chaque \'etape il n'y a pas
unicit\'e des constructions et des choix vont devoir
\^etre faits, avec des incidences sur les courbes obtenues. 

 Ensuite,
avec l'algorithme d\'evelopp\'e dans [HMDP2], on construit les familles  de courbes
associ\'ees
\`a cette triade, en commen\c cant par les minimales et en effectuant des biliaisons si
n\'ecessaire.

Nous donnons plusieurs exemples construits selon cette proc\'edure et une premi\`ere
analyse, sommaire, montrant comment les diverses \'etapes doivent \^etre r\'ealis\'ees
pour obtenir une famille de courbes de degr\'e et genre donn\'es et dont la variation de la
cohomologie est prescrite.

\vskip 0.3 cm

\noindent {\bf Notation.} Dans ce qui suit on utilisera souvent la notation chiffr\'ee
 suivante : le module
$\bigoplus _{i=1}^r R_A(-n_i)^{\a_i}$ sera not\'e $n_1^{\a_1}, \cdots, n_r^{\a_r}$.
Pour \'eviter le risque de confusion avec le module nul on notera $\ul 0$ le module
$R_A$.

\vskip 0.3 cm

\tarte {a) Construction d'une triade \`a partir de son c\oe ur, de son conoyau et d'une
 extension de longueur $2$ de ces modules}

Dans ce paragraphe $A$ d\'esigne un anneau noeth\'erien.

\vskip 0.3 cm

Si $\Lpont = L_1 \Fl {d_1} L_0
\Fl {d_0} L_{-1}$ est une triade, on    lui  associe son c\oe ur
$H$, son conoyau $C$  et l'extension de $R_A$-modules
$$0 \fl H \fl E \Fl{d} L_{-1} \fl C \fl 0,$$
dont la classe est  un
\'el\'ement $u$ de $\Ext^2_{R_A} (C,H)$, cf. 1.1. 
Nous montrons
ci-dessous comment on peut, inversement, construire explicitement, \`a partir de $C,H,u$,
une triade qui redonne ces \'el\'ements. Dans le cas o\`u $A$ est un \avd\ et $\Lpont$ une
triade majeure \'el\'ementaire la construction ci-dessous redonne la triade initiale \`a \ps\
pr\`es.

\vskip 0.3 cm

On consid\`ere  une r\'esolution libre gradu\'ee de $C$ :
$$ \cdots \fl P_3 \Fl {\de_2} P_2 \Fl {\de_1} P_1 \Fl {\de_0} P_{0} \fl C \fl 0 \leqno {(1)}$$
et une r\'esolution libre gradu\'ee de $H$ :
$$\cdots \fl Q_1 \Fl {\g} Q_0 \Fl {p} H \fl 0.  \leqno {(2)}$$
On d\'esigne par $\Ppont$ et $\Qpont$ les complexes obtenus \`a partir de ces suites exactes en
oubliant les termes $C$ et $H$. Un \'el\'ement $u$ de  $\Ext_{R_A}^2 (C,H)$ correspond  \`a un
morphisme
$\wh u : P_2
\fl H$ v\'erifiant $\wh u \de_2 =0$ (modulo ceux qui proviennent de $P_1$). L'homomorphisme $\wh
u$ induit un morphisme de complexes $\upont : \Ppont [2]\fl \Qpont$, avec des fl\`eches
$u_0 : P_2 \fl Q_0$, $u_1 : P_3 \fl Q_1$ etc. Soit $\Kpont$  le c\^one  du morphisme $\upont
$ et
$\Lpont$ le complexe obtenu en tronquant $\Kpont$ en degr\'es $-1,0,1$. On a donc 
$L_{-1}= P_{0}$, $L_0= P_1 \oplus Q_0$, $L_1= P_2 \oplus Q_1$ et  les fl\`eches $d_0 = (\de_0 \
0)$ et $d_1 = \pmatrix{ \de_1& 0 \cr u_0& \g \cr}$, avec le r\'esultat suivant :

\th {Proposition 5.1}. \lign 
1)  Le complexe $\Lpont$ construit ci-dessus est une triade majeure
qui a 
 pour conoyau $C$, pour
c\oe ur
$H$ et dont l'\'el\'ement  de $\Ext^2_{R_A} (C,H)$ associ\'e est $u$. 
Cette triade est bien d\'efinie \`a un \ps\ fort pr\`es. On la note
$\Lpont (C,H,u)$. \lign
2) On suppose que est $A$ un \avd. Soit  $\Lpont$  une triade majeure
\'el\'ementaire et soient $C$, $H$ et $u$ le conoyau, le c\oe ur et l'extension associ\'es
\`a $\Lpont$. Alors, $\Lpont$ et $\Lpont (C,H,u)$ sont pseudo-isomorphes.

\dem 1) Par d\'efinition du c\^one on a la suite exacte de complexes :
$0 \fl \Qpont \fl \Kpont \fl \Ppont \fl 0$ 
et la suite exacte d'homologie 
 associ\'ee donne les isomorphismes cherch\'es. Des choix diff\'erents des r\'esolutions de $C$ et
$H$ ainsi que des rel\`evements $u_i$ d\'efinissent le m\^eme complexe 
$\Kpont$  
\`a
\qis\ pr\`es, de sorte que $\Lpont$, qui est obtenu par troncature, est bien d\'efini \`a  un \ps\
fort pr\`es.  

Le point 2) r\'esulte de 1.34.2.

\remas {5.2} \lign
1) On peut d\'ecrire explicitement  l'extension de $C$ par $H$ associ\'ee \`a
$u$. On consid\`ere le
module
$E$, somme amalgam\'ee de
$H$ et $P_1$ au-dessus de $P_2$, c'est-\`a-dire le module donn\'e par la suite exacte :
$$P_2 \Fl {(\wh u,\de_1)} H \oplus P_1 \fl E \fl 0. \leqno {(3)}$$
 L'extension de
longueur
$2$ d\'efinie par $u$ est alors :
$$ 0 \fl H \fl E \fl P_{0} \fl C \fl 0$$
et on obtient une triade pseudo-isomorphe \`a  $\Lpont (C,H,u)$ en juxtaposant $L_{-1}= P_{0}$ et
une r\'esolution
$L_1
\fl L_0
\fl E \fl 0$. \lign
2) Au sens de la remarque pr\'ec\'edente, la construction de la triade au moyen du  c\^one de
$\upont$ effectu\'ee ci-dessus consiste, pour r\'esoudre $E$,
\`a r\'esoudre $H \oplus P_1$ par $Q_0 \oplus P_1$. Dans certains cas cette construction n'est
pas la plus \'economique (en ce sens qu'elle donne des modules $L_i$ de rangs inutilement grands).
On a en particulier le r\'esultat suivant :

\th {Proposition 5.3}. Avec les notations pr\'ec\'edentes,  on suppose que $u$ admet un
rel\`evement $\wh u$ surjectif et on appelle $\ov u$ l'homomorphisme de $\Ker \de_0$ dans $H$
induit par $\wh u$. On a alors une suite exacte
$0
\fl
\Ker
\ov u
\fl P_1
\fl E
\fl 0$. On obtient une r\'esolution  $L_1 \fl L_0 \fl E \fl 0$ en prenant $L_0= P_1$ et pour $L_1$
le premier module d'une r\'esolution de
$\Ker
\ov u$. La triade $L_1 \fl L_0 \fl P_{0}$ obtenue est pseudo-isomorphe \`a $\Lpont (C,H,u)$.

\dem  Cela r\'esulte de la surjectivit\'e de $\ov u$ et de la suite exacte $(3)$ ci-dessus.

\tarte  {b) Le cas d'un \avd}

Dans ce paragraphe $A$ d\'esigne  un anneau de
va\-lua\-tion discr\`ete d'uniformisante $a$ et de corps r\'esiduel $k$ et on commence par montrer
un r\'esultat qui compare les groupes d'extensions  sur $R$ et sur
$R_A$.

\th {Proposition 5.4}. \lign
1)  Soient $M'$ et $M''$ deux $R_A$-modules gradu\'es annul\'es par $a$
(de sorte que ce sont  aussi des $R$-modules). On a une suite exacte :
$$ \eqalign {\cdots &\fl \Ext^1_{R} (M'',M') \Fl{j_1}  \Ext^1_{R_A} (M'',M') \Fl {\psi_1}
\Hom_R(M'',M')
\cr &\fl
\Ext^2_{R} (M'',M') \Fl {j_2} \Ext^2_{R_A} (M'',M') \Fl{\psi_2} \Ext^1_{R} (M'',M') \fl
\cdots\cr}$$ 2) Si, de plus, $A$ est une $k$-alg\`ebre, les fl\`eches $j_k$, pour $k \geq 1$,
admettent des r\'etractions et on a pour tout $k \geq 1$ une suite exacte scind\'ee :
$$0 \fl \Ext^k_{R}
(M'',M') \Fl {j_k} \Ext^k_{R_A} (M'',M') \Fl{\psi_k} \Ext^{k-1}_{R} (M'',M')  \fl 0.$$

\dem 1) On consid\`ere une r\'esolution libre gradu\'ee de $M''$ sur $R_A$ :
$$Ê\cdots \fl L_2 \Fl {\de_1} L_1 \Fl {\de_0} L_0 \Fl{p} M'' \fl 0$$
et on pose $K = \Im \de_0$. En tensorisant par $k$ on obtient un complexe
$$Ê\cdots \fl \ov L_2 \Fl {\ov\de_1} \ov L_1 \Fl {\ov\de_0} \ov L_0 $$
dont les seuls groupes d'homologie non nuls sont  $M'' \T_A k$ en degr\'e
$0$ et
$\Tor_1^A (M'',k)$ en degr\'e $1$,  qui sont tous deux isomorphes \`a  $M''$  puisque $M''$ est
annul\'e par $a$. Si on pose $\ov K = K
\T_A k$ on a donc une suite exacte :
$$ 0 \fl M'' \fl \ov K \fl \ov L_0 \Fl{\ov p} M'' \fl 0. \leqno {(4)}$$
On pose $Z = \Ker \ov p$ et on note alors les faits suivants : 

\noindent 1) Pour $i \geq 1$ on a un isomorphisme 
$$\Ext^i_{R_A} (K,M') \simeq \Ext^{i+1}_{R_A} (M'',M').$$
2) Comme $\cdots \ov L_2 \Fl {\ov\de_1} \ov L_1 \fl \ov K \fl 0$ est une r\'esolution de $\ov K$
et comme on a, pour $i \geq 1$, un isomorphisme
$$\Hom_{R_A} (L_i, M') \simeq \Hom_{R} (\ov L_i, M')$$
(car $M'$ est annul\'e par $a$), on en d\'eduit, pour $i \geq 0$ des isomorphismes :
$$\Ext^i_{R_A} (K,M')  \simeq \Ext^i_{R} (\ov K,M') .$$
3) La partie droite de la suite exacte $(4)$ donne  pour $i \geq 1$ les isomorphismes :
$$\Ext^{i+1}_{R} (M'',M') \simeq \Ext^{i}_{R} (Z,M').$$
La conclusion du point 1) s'obtient alors en \'ecrivant la suite d'homologie associ\'ee \`a la
partie gauche de
$(4)$ (avec une variante pour $k=1$ que nous laissons au lecteur).

2) La fl\`eche $j_k$ s'interpr\`ete en termes d'extensions : elle consiste simplement \`a
consid\'erer un $R$-module $M$ extension de $M''$ par $M'$ comme un $R_A$-module annul\'e par
$a$. Si $A$ est une $k$-alg\`ebre cette fl\`eche admet une r\'etraction \'evidente qui consiste
\`a consid\'erer un $R_A$-module comme un  $R$-module par restriction des scalaires, ce qui donne
le point 2).

\vskip 0.3 cm

On reprend maintenant les notations du \S {\it a)} dans le cas d'un anneau de valuation
discr\`ete.

 Soit $u \in \Ext^2_{R_A} (C,H)$. La fonctorialit\'e   du foncteur $\Ext^2$ appliqu\'ee aux
morphismes canoniques $H \fl H \T_A k$ et $\Tor_1^A(C,k) \fl C$ donne un homomorphisme
$$ \f : \Ext^2_{R_A} (C,H) \fl \Ext^2_{R_A} (\Tor_1^A(C,k),H\T_A k).$$
Comme les deux modules $H\T_A k $ et $\Tor_1^A(C,k)$ sont annul\'es  par $a$ on a aussi le
morphisme 
$$\psi_2 : \Ext^2_{R_A} (\Tor_1^A(C,k),H\T_A k) \fl \Ext^1_{R} (\Tor_1^A(C,k),H\T_A k)$$
d\'efini en 5.4 ci-dessus. Si on pose   $\theta = \psi_2\, \f$  on v\'erifie ais\'ement 
 le r\'esultat suivant :

\th {Proposition 5.5}. Soient $C,H$ des $R_A$-modules, $u \in \Ext^2_{R_A} (C,H)$ et $\Lpont=
\Lpont (C,H,u)$ la triade d\'efinie en 5.1. On note $V$ le foncteur associ\'e \`a $\Lpont$.
Alors, la classe de $V(k)$ comme extension de $\Tor_1^A(C,k)$  par $H\T_A k $  (cf. 1.5) est
\'egale \`a
$\theta (u)$.

\tarte {c) Sous-quotient associ\'e \`a une triade}

\def \tri {\triangleleft}

\th {D\'efinition 5.6}.
 Soit $k$ un corps et soient $M$ et
$M_0$ deux $R$-modules gradu\'es. 
On dit que $M$ est un {\bf sous-quotient} de
$M_0$ si $M$ est  un quotient d'un sous-module gradu\'e de $M_0$ (ou, ce qui revient au
m\^eme, un sous-module d'un quotient de $M_0$). On note alors $M \triangleleft M_0$.

Se donner $M$ comme sous-quotient de $M_0$ revient \`a se donner un drapeau de sous-modules : $M_1
\subset J \subset M_0$ avec $M = J/M_1$ ou encore un
diagramme commutatif de suites exactes  du type  suivant :
$$ \matrix { &&0&&0 \cr
&& \vf && \vf\cr && M_1&=&M_1 \cr && \Vf{l}&&\vf \cr
0 & \fl& J &\Fl {i}& M_0 &\Fl {p}& M_{-1}& \fl& 0 \cr
&& \Vf{\pi}&&\vf&&\parallel \cr
0& \fl&M&\fl& M_0/M_1& \fl& M_{-1}& \fl& 0 \cr
&&\vf &&\vf\cr &&0&&0 \cr
} \leqno {(M \tri M_0)}$$
Il revient encore au m\^eme de se donner la suite exacte verticale $0 \fl M_1 \fl J \fl M \fl 0$
 et le module $M_0$ comme  extension de
$J$ par
$M_{-1}$.

\th {D\'efinition 5.7}.
 Soit $A$   un anneau de
va\-lua\-tion discr\`ete  de corps r\'esiduel $k$. On appelle {\bf d\'eformation de sous-quotient
} la donn\'ee d'un $R$-module $M_0$, d'un sous-quotient
$M$  de
$M_0$  et d'un $R_A$-module
$M_A$ qui est une d\'eformation plate de $M$ (i.e. un $R_A$-module, plat sur $A$ tel que l'on ait
$M_A
\T_A k = M$). Une telle donn\'ee sera not\'ee $M_A \sim M \tri M_0$. 

\th {Proposition 5.8}. Soient $A$ un  anneau de
valuation discr\`ete, $K$ son corps des fractions et $k$ son corps
r\'esiduel. Soit
$\Lpont$ une
triade sur $A$  (resp. une triade modulaire, resp. une triade
repr\'esentable), $V$ le foncteur associ\'e,
$M_0= V(k)$. Soient $H$ le c\oe ur de $\Lpont$ et $C$ son conoyau. On note $H_\tau$ le
sous-module de $A$-torsion de $H$. Alors, le
$R_A$-module gradu\'e 
$M_A= H/H_\tau$ est plat sur
$A$ et on a $M_A \T_A K = V(K)$ et $M_A \T_A k = M$, o\`u $M$ est le
sous-quotient  (resp. le quotient, resp. le sous-module) de $M_0$ d\'efini par le drapeau 
$$M_1 = H_\tau
\T_A k \subset J= H
\T_A k \subset M_0$$
avec comme quotient $M_{-1} = M_0/J=
\Tor_1^A (C,k)$.
Ê\lign 
On note $\cD (\Lpont)$ la d\'eformation de sous-quotient ainsi
obtenue et on dit que la  triade joint les modules $M_0$ et
$M$ sur $A$. 

\dem Comme $A$ est un \avd, $M_A$ est un $R_A$-module gradu\'e,
plat sur
$A$. Comme $K$ est plat sur $A$ on a bien $M_A \T_A K = H \T_A K = V(K)$. On
a la suite exacte (cf. 1.5)
$$0 \fl H\T_A k \fl V(k) \fl \Tor_A^1 (C,k) \fl 0$$
qui s'\'ecrit encore $0 \fl J \fl M_0 \fl M_{-1} \fl 0$.

 Par ailleurs, la suite exacte $0
\fl H_\tau \fl H \fl H/H_\tau \fl 0$ dont le quotient est plat donne, en tensorisant par
$k$, la suite $0 \fl M_1 \fl J \fl M \fl 0$, d'o\`u un diagramme du type $(M \tri M_0)$ qui donne
les propri\'et\'es annonc\'ees.

L'assertion sur les cas  modulaire et repr\'esentable r\'esulte des \'equivalences $\Lpont$
modulaire
$\equi$
$C$ plat et $\Lpont$ repr\'esentable $\equi$ $H$ plat (cf. 4.1 et 4.2).

\vskip 0.3 cm

On obtient aussit\^ot, avec 5.8 le corollaire suivant, annonc\'e dans l'introduction :

\th {Corollaire 5.9}. Soient $A$ un  anneau de
valuation discr\`ete, $K$ son corps des fractions et $k$ son corps
r\'esiduel. Soit $\cC$ une famille de courbes param\'etr\'ee par $A$. Alors, le module
de Rao  $H^1_* \cJ_{\cC \T_A K}$  au point g\'en\'erique est le point g\'en\'erique
  d'une d\'eformation
plate sur
$A$ d'un sous-quotient du module de Rao  $H^1_* \cJ_{\cC \T_A k}$ au point ferm\'e.

\tarte {d) 
Construction de triades \`a partir de d\'eformations de sous-quotients : analyse des conditions
n\'eces\-saires}

\`A partir de ce paragraphe l'anneau $A$ est un \avd\ d'uniformi\-san\-te $a$ et de corps
r\'esiduel
$k$.

On consid\'ere une d\'eformation de sous-quotient $M_A \sim M \tri M_0$, le sous-quotient $M \tri
M_0$ \'etant donn\'e par un diagramme $(M \tri M_0)$  comme  ci-dessus et on  cherche \`a
construire  une triade 
$\Lpont$
\'el\'ementaire (i.e., dont le conoyau  est de torsion)
 qui redonne cette d\'eformation par 5.8. 
 Si $H$ et $C$ d\'esignent le c\oe ur  et le conoyau
de 
$\Lpont$  on a, par 5.8, des conditions n\'ecessaires sur 
$C$ et $H$ :

\th {Conditions n\'ecessaires  5.10}. Avec les notations ci-dessus, on a les conditions  :
$$M_A \simeq H / H_\tau, \quad J= H \T_A k, \quad M_{-1} =
\Tor_1^A (C,k).$$ On en d\'eduit $M= H/H_\tau \T_A k$ et $M_1 = H_\tau \T_A k$. Si $C$ et $H$
v\'erifient ces conditions on dira qu'ils sont compatibles avec la d\'eformation de sous-quotient
$M_A \sim M \tri M_0$.

\remas {5.11} \lign
1) Il est essentiel de noter qu'en g\'en\'eral, les conditions  5.10 ne d\'eterminent pas
enti\`erement
$C$ et
$H$.\lign
2) La condition sur $C$ donne une  indication sur la structure de $A$-module de $C$.
Comme $C$ est de torsion on a, pour chaque degr\'e $n$, $C_n = \bigoplus_{i=1}^{r_n} A/(a^{n_i})$
avec
$n_i >0$ et on a alors 
$\Tor_1^A (C_n,k) \simeq C_n \T_A k = k^{r_n}= M_{-1,n}$ qui donne $r_n$ (mais pas les $n_i$). 
 La d\'etermination des structures de $R_A$-modules sur $C$ 
compatibles avec la donn\'ee de la structure de $R$-module de $M_{-1} =\Tor_1^A (C,k)$ n'est pas
\'evidente en g\'en\'eral.\lign
3) Du c\^ot\'e de $H$, on a d\'ej\`a $H/H_\tau = M_A$ et, comme $A$-module, on a $H= H_\tau \oplus
H/H_\tau$.   Si on pose, pour un degr\'e $n$,  $H_{\tau,n} =\bigoplus_{i=1}^{r_n} A/(a^{n_i})$, la
formule $H_\tau \T_A k= M_1$ permet de d\'eterminer ${r_n}$. Du point de vue
de la structure de $R_A$-module on a  la suite exacte 
$0 \fl H_\tau \fl H \fl H/H_\tau \fl 0$
qui se r\'eduit modulo $a$ en $0 \fl M_1 \fl J \fl M \fl 0$, mais cela ne d\'etermine pas
compl\`etement cette structure. Toutefois, lorsque $H_\tau$ est annul\'e par $a$ on peut calculer
$H$
\`a l'aide du lemme \'el\'ementaire suivant :

\th {Lemme 5.12}. Soit $H$ un $R_A$-module gradu\'e. On a une suite exacte :
$$ 0 \fl aH_\tau \fl H \fl (H/H_\tau) \ti_{(H/H_\tau) \T_A k} (H \T_A k)\fl 0.$$

Examinons  l'exemple $(4,0)$ \'evoqu\'e dans l'introduction de ce chapitre.  On consid\`ere
$M=k(-1)$ comme sous-quotient de
$M_0= R/(X,Y,Z,T^3)$ comme suit : $M_0$ est engendr\'e comme espace vectoriel par $1,t,t^2$,
avec $t^3=0$ et on prend $J= <t> \simeq R(-1)/(X,Y,Z,T^2)$, $M_{-1} = k$ et $M_1 = <t^2>
\simeq k(-2)$. Dans ce cas, il n'y a pas le choix pour la d\'eformation $M_A$ : on a $M_A = A(-1)$
avec une structure de $R_A$-module triviale. On peut alors d\'eterminer tous les modules $C$ et
$H$ compatibles avec ces donn\'ees :

\th {Proposition 5.13}. Avec les notations ci-dessus,  $C$ et $H$ sont compatibles avec la
d\'eformation de sous-quotient $M_A \sim M \tri M_0$ si et seulement si  on a  
 $C=
A/(a^n)$ avec $n>0$ et  une structure de $R_A$-module triviale et $H \simeq R_A(-1)/ (X,Y,Z,
a^mT,T^2)$ avec $m>0$.

\dem
L'assertion sur $C$ r\'esulte de la formule  $\Tor_1^A (C,k) = k$. 
 Pour $H$ on a, comme $A$-module,
$H= H_1
\oplus H_2$ avec $H_1 \simeq H/H_\tau \simeq A(-1)$ et $H_2 = H_\tau \simeq
(A/(a^m))(-2)$ avec
$m >0$. Il reste \`a pr\'eciser la structure de $R_A$-module de $H$. Soient $e$ et $f$ les
g\'en\'erateurs de $H_1$ et  $H_2$. Vu la r\'eduction des multiplications dans $J = H \T_A
k$, et quitte \`a faire un changement de variables dans $R_A$, on peut supposer
$Xe=Ye=Ze=0$ et $Te= f$. On en d\'eduit le r\'esultat.

\vskip 0.5 cm

Retournons au cas g\'en\'eral et supposons qu'on ait des modules
$C$ et
$H$  compatibles avec la d\'eformation $M_A \sim M \tri M_0$. Pour avoir une triade qui
donne cette d\'eformation il reste,  comme on l'a vu au paragraphe {\it a)}, 
\`a trouver un
\'el\'ement  convenable $u$ de $\Ext^2_{R_A} (C,H)$. Comme la donn\'ee de
$H$ et de $C$ redonne d\'ej\`a
$M_{-1}$, $M_A$,  ainsi que la suite exacte $0 \fl M_1 \fl J \fl M \fl 0$, il reste \`a choisir
$u$ pour que
 l'extension  de $M_{-1}$ par $J$ associ\'ee \`a la triade soit bien $M_0$.  Avec 5.5 et 5.10 on
obtient :

\th {Proposition 5.14}. Soient $C$ et $H$ des modules  compatibles avec la d\'eformation de
sous-quotient
$M_A \sim M \tri M_0$, soit $u \in \Ext^2_{R_A} (C,H)$ et soit
$\Lpont =
\Lpont (C,H,u)$ la triade associ\'ee. Alors, la d\'eformation de sous-quotient $\cD (\Lpont)$ est
\'egale \`a $M_A \sim M \tri M_0$ si et seulement si l'\'el\'ement  $\theta (u) \in \Ext^1_{R}
(\Tor_1^A(C,k),H\T_A k)$ (cf. 5.5) est la classe de l'extension
$M_0$ de
$M_{-1}$ par $J$.

Dans le cas de l'exemple $(4,0)$ nous verrons aux paragraphes  {\it e)} et {\it h)} plusieurs
cons\-tructions de triades associ\'ees au sous-quotient donn\'e en 5.13.

\tarte {e) Construction de triades \`a partir de sous-quotients : la construction
triviale}

Soit $A$ un \avd\ d'uniformisante $a$, de corps des fractions $K$ et
de corps r\'esiduel $k$. On suppose de plus que $A$ est une {\bf $k$-alg\`ebre}. Nous
allons  construire une triade sur $A$ qui joint un module
$M_0$ et un de ses sous-quotients $M$. Cette construction est tr\`es simple   et valable dans tous
les cas mais elle ne donnera pas, en g\'en\'eral, les courbes de plus petit degr\'e, cf. 5.18.3 et
\S {\it i)}).

\vskip 0.3 cm

Soit
$M_0$ un
$R$-module gradu\'e de longueur finie et
$M$ un sous-quotient de $M_0$  donn\'e par un diagramme
commutatif
$(M
\tri M_0)$ comme en 5.6. Avec les notations de ce diagramme on pose $j=il$.

\th {Proposition 5.15}.  Avec les notations ci-dessus, le complexe 
$$\Mpont = \big(M_1 \T_k A \Fl {j \T a} M_0 \T_k A \Fl {p \T a} M_{-1} \T_k A \big)$$ est 
une
triade. Elle est   form\'ee de
$R_A$-modules gradu\'es  libres et de type fini sur $A$ et le foncteur  associ\'e $V$, v\'erifie
$V(K) = M\T_k K$ et $V(k) = M_0$. \lign Cette triade sera appel\'ee la {\bf triade triviale}
associ\'ee au sous-quotient $M \tri M_0$ (cf. 5.16). Si $M$  est un quotient (resp. un
sous-module) de $M_0$ la triade $\Mpont$ est  modulaire (resp.  repr\'esentable).

\dem  C'est une v\'erification imm\'ediate. Pour la derni\`ere assertion,  on a
$M_{-1}=0$ (resp.
$M_1 =0$) si
$M$ est un quotient (resp. un sous-module) de $M_0$. En vertu de 4.1 et 4.2 la triade
est  alors modulaire (resp. repr\'esentable).
\vskip 0.3 cm

\remas {5.16} \lign
1) Expliquons l'appellation ``triviale'' pour la triade construite ci-dessus. On v\'erifie
d'abord  les formules suivantes :
$C \simeq M_{-1}$ (consid\'er\'e comme $R_A$-module annul\'e par $a$), 
$H_\tau  \simeq M_1$ (consid\'er\'e comme $R_A$-module annul\'e par $a$), ce qui correspond aux
choix triviaux pour les structures de $A$ et de $R_A$-modules de $C$ et $H_\tau$ (cf. 5.11). On a
ensuite
$M_A
\simeq M\T_k A$ : $M_A$ est une d\'eformation triviale du sous-quotient $M$ de $M_0$. On d\'eduit
de ces formules la structure de $H$ comme produit fibr\'e  $H  \simeq M_A \times_M J$ (cf.
5.12).\lign L'\'el\'ement $u \in \Ext^2_{R_A} (C,H)$  correspond \`a l'extension 
$0 \fl H \fl E \fl M_{-1}\T_A k \fl C\fl 0$, cf. 1.1.
On v\'erifie que l'on a $E_\tau  \simeq H_\tau \simeq M_1$ et $E \T_A k \simeq M_0$, de
sorte que, comme $E_\tau$ est annul\'e par $a$, $E$ est un
produit fibr\'e en vertu de 5.12.  On a aussi   $(E/E_\tau) \T_A k \simeq
M_0/M_1$, donc
$E/E_\tau$   est une d\'eformation plate de
$M_0/M_1$ et on v\'erifie que cette d\'eformation est, elle aussi, triviale : $E/E_\tau =
(M_0/M_1)
\T_k A$.
\lign On peut aussi comprendre l'\'el\'ement $u$ de la mani\`ere suivante. On part de
l'extension 
$0 \fl J \fl M_0 \fl M_{-1} \fl 0$ (dont la classe doit \^etre \'egale \`a $\theta (u)$) et
l'on en d\'eduit trivialement, en tensorisant par $A$ sur $k$,  une extension de $R_A$-modules
de
$M_{-1} \T_k A$ par $J \T_A k$. En juxtaposant cette extension avec la pr\'esentation
\'evidente de $C= M_{-1}$ :
$0 \fl M_{-1} \T_k A \Fl{1 \T a} M_{-1} \T_k A  \fl C \fl 0$ on obtient une extension de $C$
par $J \T_k A$. Comme $H$ est le produit fibr\'e $ (M \T_k A) \times_M J$, on a une fl\`eche
canonique de $J \T_k A$ dans $H$ et l'extension $u$ s'en d\'eduit par fonctorialit\'e.

Le nom de
triviale donn\'ee \`a la triade de 5.15 correspond donc au fait que, pour tous ses \'el\'ements
constitutifs pour qui un choix se pr\'esente  \`a partir des conditions n\'ecessaires 5.10 et
5.14, ce choix est effectu\'e de mani\`ere triviale.
 \lign
2) L'existence de la construction triviale montre qu'il y a toujours une triade qui joint un
module \`a un sous-quotient. Une question  plus d\'elicate est l'existence d'une triade qui
corresponde \`a  une d\'eformation de sous-quotient $M_A \sim M \tri M_0$ donn\'ee (pas
n\'ecessairement triviale). Nous montrerons dans un travail ult\'erieur qu'une telle triade
n'existe que si une certaine condition sur la d\'eformation (la nullit\'e d'un  \'el\'ement
de $\Ext^3_R (M_{-1}, M_1)$) est v\'erifi\'ee.
 \lign
3) La construction triviale
 n'est pas en g\'en\'eral la seule possible qui corresponde \`a un sous-quotient donn\'e.  Nous
verrons ci-dessous d'autres exemples de constructions, soit avec  des foncteurs triadiques
diff\'erents de celui de la triade triviale (cf. 5.19), soit avec le m\^eme foncteur, mais des
triades distinctes (cf. 5.20).

\expls {5.17}  Dans tous les exemples on travaille sur une $k$-alg\`ebre qui est un \avd\
$A$  d'uniformisante
$a$, de corps des fractions $K$ 
 et de corps r\'esiduel $k$.

\vskip 0.3 cm

 1) Soit $m$ l'id\'eal maximal
$(X,Y,Z,T)$ de
$R$. L'exemple le plus simple de triade s'obtient en prenant pour $M_0$ le $R$-module
$R/m=k$ et pour $M$ le module nul, vu comme quotient de $M_0$. On a alors la triade
modulaire  d\'efinie par le module $H = A/(a)$ consid\'er\'e comme $R_A$-module (en
degr\'e $0$ et avec des multiplications triviales par $X,Y,Z,T$). Les valeurs du foncteur
$V$ aux points g\'en\'erique et ferm\'e de $\Spec A$ sont bien $0$ et $k$. 

Une  r\'esolution majeure de cette triade est
$ R_A(-1)^4 \oplus R_A \Fl {d_1} R_A \fl 0$
avec $d_1 = (X\, Y\, Z\, T\, a)$. 

\vskip 0.3 cm

2) Un deuxi\`eme exemple consiste \`a prendre les m\^emes modules $k$ et $0$ mais
en consid\'erant cette fois $0$ comme sous-module de $k$. La triade est alors
donn\'ee par  $0 \fl A \Fl {a} A$ et on a encore les valeurs prescrites
de $V$. On notera qu'on a $V(A) = 0$ de
sorte que $V$, qui n'est pas nul, n'est pas de la forme ($V(A) \T_A \tas$). 

La triade majeure associ\'ee s'\'ecrit, en notation condens\'ee :
$$   1^4,2^6 \Fl {d_1}  \ul 0, 1^4 \Fl {d_0}
\ul 0$$  $${\rm avec} \qquad d_0 = (a\, X\, Y\, Z\, T) \qquad {\rm et} \qquad
d_1=\pmatrix {U&0&
\cr -aI_4& V\cr}$$
o\`u $I_4$ est la matrice identit\'e d'ordre $4$, $U$ la matrice $(X,\, Y,\, Z,\, T)$
et o\`u l'on a pos\'e
$$V= \pmatrix {Y&Z&T&0&0&0 \cr -X&0&0&Z&T&0 \cr 0&-X&0&-Y&0&T \cr 0&0&-X&0&-Y&-Z
\cr}.$$

3) Dans cet exemple, qui part du sous-quotient correspondant \`a une \'eventuelle famille
de courbes de $H_{4,0}$, la triade obtenue n'est  ni modulaire, ni
repr\'esentable.

 On
consid\`ere le module $M_0=R/(X,Y,Z,T^3)$ qui est de dimension $1$ en degr\'es $0,1,2$ et
son sous-quotient $M=k(-1)$, concentr\'e en degr\'e $1$. La m\'ethode de 5.15 ci-dessus
donne une triade mineure qui joint ces modules sur $A$ :
$$ R_A(-2)/ (X,Y,Z,T) \Fl {\de_1} R_A/(X,Y,Z,T^3) \Fl{\de_0} R_A/(X,Y,Z,T)$$ 
avec $\de_1 (1) = at^2$ et $\de_0(1) = a$ (o\`u $t$ d\'esigne l'image de $T$ dans le
quotient). 

On obtient une r\'esolution majeure  $\Lpont = (L_1 \Fl
{d_1} L_0 \Fl{d_0} L_{-1})
$ de la triade  pr\'ec\'edente   en posant $L_1= 1^3,2^7 ,3$ ;
$L_0= \ul 0,1^4$ ; $L_{-1}=\ul 0$  avec les fl\`eches 
$$ d_1 = \pmatrix { X& Y & Z& 0&0&0&0&0&0&aT^2&T^3 \cr
-a&0&0&Z&T&0&0&0&Y&0&0 \cr
0&-a&0 &0&0&Z&T&0&-X&0&0\cr
0&0&-a&-X&0&-Y&0&T&0&0&0\cr
0&0&0&0&-X&0&-Y&-Z&0&-a^2T&-aT^2 \cr
}$$
et $d_0= (a\ X\ Y \ Z \ T)$.

Nous exhiberons plus loin une triade, joignant les m\^emes modules et d\'efinissant le
m\^eme foncteur que $\Lpont$ mais qui ne lui sera pas pseudo-isomorphe.

\tarte {f) Calcul des familles de courbes obtenues \`a partir des triades triviales}

\expls {5.18} Nous d\'eterminons ici les familles minimales de courbes associ\'ees (entre
autres) aux exemples de 5.17. Nous utilisons pour cela les r\'esultats de [HMDP2] o\`u nous
avons construit une fonction $q$ qui g\'en\'eralise au cas des familles de courbes celle de
[MDP1]. Dans tout ce qui suit c'est de cette nouvelle fonction $q$ qu'il sera question.

1) On consid\`ere la triade
modulaire  de 5.17.1 qui joint le module $k$ au module nul (vu comme  quotient de $k$).
On voit ais\'ement que le noyau $N$  est l'image de la matrice $s : R_A(-1)^4 \oplus
R_A(-2) ^6
\fl R_A \oplus R_A(-1)^4$ avec
$$
s = \pmatrix { U & 0\cr
a \,I_4 & V \cr
}$$
o\`u $I_4$ est la matrice identit\'e d'ordre $4$, o\`u $U= (X\, Y\, Z\, T)$ et o\`u $V$ est
la matrice des relations entre $X,Y,Z,T$ d\'ej\`a vue en 5.17.2.

Le calcul de la famille minimale associ\'ee \`a cette triade a \'et\'e effectu\'e dans
[HMDP2] 3.2. Ces courbes sont  de degr\'e
$6$ et de genre
$3$, la famille est
\`a sp\'ecialit\'e constante. La courbe g\'en\'erique $C$  de la famille $\cC$ est une
courbe ACM tandis que la courbe sp\'eciale $C_0$ est une courbe de la classe de biliaison
de deux droites disjointes qui a m\^eme cohomologie qu'une courbe de bidegr\'e $(2,4)$
trac\'ee sur une quadrique lisse, en particulier son module de Rao est concentr\'e en
degr\'e $2$. Cet exemple est le premier de ceux \'evoqu\'es dans l'introduction.

\vskip 0.3 cm

1') Cet exemple g\'en\'eralise le pr\'ec\'edent. On part d'un module $M_0$ de Koszul
(cf. [MDP1] IV.6 ou [MDP3] V.2), $M_0 = R/(f_1,f_2,f_3,f_4)$, o\`u les $f_i$ sont des
polyn\^omes homog\`enes de degr\'es
$n_i$ avec $0<n_1 \leq n_2 \leq n_3 \leq n_4$ et forment une suite r\'eguli\`ere de $R$. On
consid\`ere le module $H = R_A/(a,f_1,f_2,f_3,f_4)$ et la triade modulaire d\'efinie par
$H$. La r\'esolution de $H$ est donn\'ee par le complexe de Koszul
$$ \cdots \fl \bigoplus_{i=1}^4 R_A(-n_i) \oplus \bigoplus_{1 \leq i < j \leq 4} R_A(-n_i-n_j)
\Fl {s} R_A \oplus  \bigoplus_{i=1}^4 R_A(-n_i) \Fl {a \oplus U} R_A \fl H \fl 0$$
avec $U =  (f_1\, f_2\, f_3\, f_4)$ et $s$ la matrice par blocs comme en 5.17.2 mais o\`u,
dans $V$, les $f_i$ remplacent $X,Y,Z,T$. Cette  r\'esolution donne \`a la fois  la triade
et la fl\`eche
$s$ qui permet le calcul des invariants (cf. [HMDP2] 3.1). On v\'erifie d'abord qu'on a
$b_0 = n_1-1$. Il est clair qu'on a
$b_0
\geq n_1-1$ et il y a deux cas : si
$n_1 <n_2$ on a $\a_{n_1} = 1$ et $\beta_{n_1}=0$, donc $n_1-1 \geq b_0$ par [HMDP2] {\it
loc. cit.} ; si
$n_1=n_2$ on a  $\a_{n_1} = \beta_{n_1}=1$ mais le module engendr\'e par les colonnes de la
matrice $s_1(t)$ n'est pas libre  de rang $1$  (car la suite est r\'eguli\`ere). 

L'analyse des valeurs de $\a_n$ et $\beta_n$ se fait alors comme dans [MDP1] ou [MDP3]
 et donne les valeurs de la fonction $q$ ou plut\^ot celles de $q\sh$ : on a,
en posant comme dans {\it loc. cit.}, $\mu = \sup (n_1+n_4, n_2+n_3)$, 
$$q\sh (n) = \cases {0,& si $n<n_1+n_2$ ; \cr 1,& si $n_1+n_2 \leq n < n_1+n_3$ ; \cr
2,& si $n_1+n_3 \leq n < \mu$ ; \cr
3,& si $n\geq \mu$. \cr}$$

En particulier on obtient une famille minimale de courbes param\'etr\'ee par $A$ avec la
r\'esolution
$$0 \fl \cO_{\bP_A}(-n_1-n_2) \oplus \cO_{\bP_A} (-n_1-n_3) \oplus \cO_{\bP_A} (-\mu)
\fl \cN \fl \cJ_{\cC} (h) \fl 0$$ o\`u $\cN$ est d\'efini par la suite  $0 \fl \cN \fl 
\cO_{\bP_A}
\bigoplus \
\bigoplus_{i=1}^4 \cO_{\bP_A}(-n_i) \Fl {(a, U)} \cO_{\bP_A} $. Comme dans l'exemple
1), la famille est
\`a sp\'ecialit\'e constante, la courbe g\'en\'erique est ACM et la courbe sp\'eciale
est dans la classe de biliaison du module de Koszul $M_0$. On notera que cette derni\`ere
courbe n'est jamais une courbe minimale.

2) On consid\`ere la triade vue en 5.17.2 qui joint encore le module nul et le module $k$,
mais avec $0$ vu maintenant  comme sous-module de $k$.  Le calcul de la r\'esolution de
$N$ donne la matrice
$s : 2^6, 3^6 \fl 1^4, 2^6$ :
$$s= \pmatrix {V&0\cr a I_4& V' \cr}$$
o\`u $V$ et $V'$ sont les matrices du complexe de Koszul associ\'e \`a la suite
$(X,Y,Z,T)$.  Le calcul de la fonction $q$ a \'et\'e effectu\'e dans [HMDP2] 3.3. La
famille minimale de courbes associ\'ee
\`a cette triade  est encore une famille de courbes de degr\'e
$6$ et genre $3$ comme dans l'exemple 1), mais, cette fois-ci, il s'agit d'une famille \`a
postulation constante dont la courbe g\'en\'erique est ACM et dont la courbe sp\'eciale est
dans la classe de biliaison de deux droites disjointes et a m\^eme cohomologie que la r\'eunion
d'une quartique plane et de deux droites qui coupent chacune la quartique transversalement en
un point. En particulier le module de Rao a un unique terme non nul en degr\'e $1$ (et non
plus en degr\'e
$2$ comme dans l'exemple 1). Cet exemple est le deuxi\`eme de ceux \'evoqu\'es dans
l'introduction. Pour une
\'etude plus approfondie du sch\'ema de Hilbert
$H_{6,3}$ cf. [AA].

\vskip 0.3 cm

On notera que si on fait une biliaison $(6,+2)$ \`a partir de la famille 1) et une
biliaison $(4,+3)$ \`a partir de la famille 2), on obtient des familles de courbes lisses
de degr\'e $18$ et genre $39$, dont les fibres sp\'eciales ont pour modules $k(-4)$, et qui
sont limites de courbes ACM de deux fa\c cons diff\'erentes. On a donc, dans le sch\'ema de
Hilbert des courbes lisses
$H_{18,39}^0$, deux composantes irr\'eductibles qui ont une intersection non vide. On
retrouve ainsi l'exemple de Sernesi, cf. [S] ou [MDP1], X 5.8 p. 193.

\vskip 0.3 cm

3) On reprend la triade $\Lpont$  de l'exemple 5.17.3 qui joint le
module
 $M_0 = R/(X,Y,Z,T^3)$  et son sous-quotient $M=k(-1)$ et dont on esp\`ere qu'elle peut
fournir une famille de $H_{4,0}$.  

On calcule (en utilisant Macaulay) une r\'esolution minimale du noyau $N$ de la
triade. On obtient la matrice 
$s : 2^3,3^8,4^3\fl 1^3,
2^7,3$ avec
$$s={\pmatrix{ \scriptstyle -Z&\scriptstyle 0&\scriptstyle -Y&\scriptstyle
\scriptstyle 0&\scriptstyle 0&\scriptstyle 0&\scriptstyle 0&\scriptstyle 0&\scriptstyle 
aT^2&\scriptstyle 0&\scriptstyle 0&\scriptstyle T^3&\scriptstyle 0&\scriptstyle 0\cr \scriptstyle 
0&\scriptstyle -Z&\scriptstyle X&\scriptstyle 0&\scriptstyle 0&\scriptstyle 
0&\scriptstyle
0&\scriptstyle 0&\scriptstyle 0&\scriptstyle aT^2&\scriptstyle 0&\scriptstyle 
0&\scriptstyle
T^3&\scriptstyle 0\cr \scriptstyle  X&\scriptstyle Y&\scriptstyle 0&\scriptstyle 0&\scriptstyle
0&\scriptstyle 0&\scriptstyle 0&\scriptstyle 0&\scriptstyle 0&\scriptstyle 
0&\scriptstyle aT^2&\scriptstyle 0&\scriptstyle 0&\scriptstyle T^3 \cr \scriptstyle  a&\scriptstyle
0&\scriptstyle 0&\scriptstyle T&\scriptstyle 0&\scriptstyle 0&\scriptstyle -Y&\scriptstyle
0&\scriptstyle 0&\scriptstyle 0&\scriptstyle 0&\scriptstyle 0&\scriptstyle 0&\scriptstyle 0
\cr \scriptstyle  0&\scriptstyle 0&\scriptstyle 0&\scriptstyle -Z&\scriptstyle 0&\scriptstyle
-Y&\scriptstyle 0&\scriptstyle 0&\scriptstyle -a^2T&\scriptstyle 0&\scriptstyle
0&\scriptstyle -aT^2&\scriptstyle 0&\scriptstyle 0 \cr \scriptstyle  0&\scriptstyle a&\scriptstyle
0&\scriptstyle 0&\scriptstyle T&\scriptstyle 0&\scriptstyle X&\scriptstyle 0&\scriptstyle
0&\scriptstyle 0&\scriptstyle 0&\scriptstyle 0&\scriptstyle 0&\scriptstyle 0\cr \scriptstyle 
0&\scriptstyle 0&\scriptstyle 0&\scriptstyle 0&\scriptstyle -Z&\scriptstyle X&\scriptstyle
0&\scriptstyle 0&\scriptstyle 0&\scriptstyle -a^2T&\scriptstyle 0&\scriptstyle
0&\scriptstyle -aT^2&\scriptstyle 0 \cr \scriptstyle  0&\scriptstyle 0&\scriptstyle 0&\scriptstyle
X&\scriptstyle Y&\scriptstyle 0&\scriptstyle 0&\scriptstyle 0&\scriptstyle 0&\scriptstyle
0&\scriptstyle -a^2T&\scriptstyle 0&\scriptstyle 0&\scriptstyle -aT^2 \cr \scriptstyle  0&\scriptstyle
0&\scriptstyle a&\scriptstyle 0&\scriptstyle 0&\scriptstyle T&\scriptstyle Z&\scriptstyle
0&\scriptstyle 0&\scriptstyle 0&\scriptstyle 0&\scriptstyle 0&\scriptstyle 0&\scriptstyle 0
\cr \scriptstyle  0&\scriptstyle 0&\scriptstyle 0&\scriptstyle 0&\scriptstyle 0&\scriptstyle
0&\scriptstyle 0&\scriptstyle T&\scriptstyle -X&\scriptstyle -Y&\scriptstyle -Z&\scriptstyle
0&\scriptstyle 0&\scriptstyle 0 \cr \scriptstyle  0&\scriptstyle 0&\scriptstyle 0&\scriptstyle
0&\scriptstyle 0&\scriptstyle 0&\scriptstyle 0&\scriptstyle -a&\scriptstyle 0&\scriptstyle
0&\scriptstyle 0&\scriptstyle -X&\scriptstyle -Y&\scriptstyle -Z \cr \scriptstyle  }}$$

 On trouve dans ce cas $\a_2= \beta_2 =2$ et $b_0=1$, puis $\a_3 = \beta_3 = 6$,
donc la fonction
$q$ est donn\'ee par
$q(2)=1$,
$q(3)=4$,
$q(4)=1$ et on obtient
 une famille minimale de courbes de degr\'e $12$ et genre $16$. On constate donc que la
construction triviale effectu\'ee en 5.15 ne donne pas la famille de courbes de degr\'e
$4$ et genre $0$ escompt\'ee.

\tarte {g) Construction de triades modulaires : variation autour du module $H$}

Ce paragraphe reprend la probl\'ematique entam\'ee au \S {\it d)}. Rappelons (cf. [MDP1] I
3.6) que deux modules gradu\'es sont dits de m\^eme type si on passe de l'un \`a l'autre
par un changement de variables lin\'eaire dans $R$.   Les exemples suivants,  pourtant
donn\'es 
 dans le cas modulaire o\`u le module
$H$ d\'etermine la triade \`a \ps\ pr\`es (cf. 4.5), montrent
successivement :

-- qu'\'etant donn\'es un module $M_0$ et deux quotients $M$  et $M'$ de m\^eme
type, les triades triviales joignant $M_0$ et $M$ (resp. $M_0$ et $M'$) ne sont pas, en
g\'en\'eral, pseudo-isomorphes,

-- que lorsque le  quotient $M$ est fix\'e, il n'y a
pas unicit\'e (\`a
\ps\ pr\`es) de la triade joignant 
$M_0$ et  $M$.

\expls {5.19} Consid\'erons  le
module
$M_0 = R/(X, Y, Z^2, T^3)$ (de dimensions $1,2,2,1$ en degr\'es $0,1,2,3$).  Tout quotient
$M$ de $M_0$ qui est un module monog\`ene de dimensions $1,1$ en degr\'es $0,1$ est de la
forme
$R/(X,Y, \la Z + \mu T, Z^2, ZT, T^2)$ et tous ces modules sont de m\^eme type.

 Il y  a  deux
cas, ($\mu
\neq 0$  et
$\mu = 0$) qui correspondent au fait  que le noyau $M_1$ de la projection de
$M_0$ sur $M$ est monog\`ene ou non. Dans le premier cas on peut prendre, par exemple, $M_1 =
<\! z+t
\! >\
\subset M_0$, donc $M_0/M_1 = R/ (X,Y,Z+T,Z^2) = M$. Dans le second cas on a
$M'_1 = <\! z, t^2\! >\ \subset M_0$, donc $ M_0/M'_1=  R/ (X, Y, Z, T^2) = M'$. 

Les triades triviales associ\'ees \`a  ces quotients
  sont respectivement d\'etermin\'ees par les modules 
$$H= R_A/ (X, Y, a(Z+T), Z^2, 
T^3)\quad {\rm et} \quad H' =R_A/ (X, Y, aZ, Z^2, aT^2, T^3).$$
Il est facile de calculer les courbes minimales correspondantes (en d\'eterminant des
r\'esolutions de ces modules par Macaulay). Dans le cas de $H$ on obtient une famille
  de courbes de
degr\'e $7$ et genre $4$, dans le cas de $H'$ une famille de courbes de degr\'e $11$ et genre $13$. On constate donc que la triade
d\'efinie par
$H$ est  meilleure que celle associ\'ee \`a $H'$  en ce sens qu'elle donne des familles de
courbes de plus bas degr\'e. Nous montrerons  dans un travail ult\'erieur que cela est li\'e
au fait que le module $ \Tor_A^1(H,k)$ (ici isomorphe \`a $M_1 = H_\tau \T_A k$ car
$H_\tau$ est annul\'e par $a$) est monog\`ene, ce qui n'est pas le cas de
$
\Tor_A^1(H',k) \simeq M'_1$. 

\vskip 0.3 cm

Il y a toutefois un moyen d'obtenir, en partant de ce m\^eme sous-module $M'_1= <\! z,t^2\!
>$ et du quotient $M'=M_0/M'_1$, une triade donnant elle aussi une famille de courbes de
degr\'e
$7$ et genre
$4$, mais il faut utiliser cette fois une triade non triviale. Pr\'ecis\'ement, on
consid\`ere le
$R_A$-module 
$H'' = R_A/(X, Y, aZ, Z^2+aT^2, T^3)$ qui est, en
quelque sorte, une variante compact\'ee de $H$. La triade d\'efinie
par  $H''$ joint
$M'= M_0/M'_1$ et
$M_0$ comme celle donn\'ee par $H'$,  mais ces triades sont distinctes car les foncteurs
associ\'es le sont ($H'$ est un quotient propre de
$H''$).  On v\'erifie facilement que la famille minimale associ\'ee \`a $H''$ est encore de
degr\'e
$7$ et genre $4$. La diff\'erence entre ces deux triades  se voit au niveau des
sous-modules de torsion
$H'_\tau$ et
$H''_\tau$. Ces modules sont tous deux engendr\'es par les images de $Z$ et $T^2$, mais
$H'_\tau$ est annul\'e par
$a$ tandis que $H''_\tau$
est annul\'e seulement par
$a^2$. Cela peut encore se traduire, comme ci-dessus, par le fait que le module
$\Tor_A^1(H'',k)$  est monog\`ene, tandis que $\Tor_A^1(H',k)$ ne l'est pas.

\tarte {h) Constructions de triades \`a partir d'un sous-quotient : bis}

Ce paragraphe reprend la probl\'ematique du \S {\it d)} au cran suivant : 
on suppose qu'on a d\'etermin\'e (\`a partir d'un sous-quotient $M$ de $M_0$)  deux
$R_A$-modules
$C$ et
$H$ compatibles. On suppose
$C$ de torsion et on cherche une triade majeure \'el\'ementaire $\Lpont = L_1 \Fl {d_1} L_0
\Fl {d_0} L_{-1}$ qui admette
$C$ et
$H$ comme conoyau et c\oe ur. On  a vu au paragraphe {\it a)}  qu'un
\'el\'ement
$u$ de
$\Ext_{R_A}^2 (C,H)$ fournit une telle triade, pourvu que l'image $\theta (u)$ de cet \'el\'ement
dans
$\Ext^1_{R} (\Tor_1^A(C,k),H\T_A k)$ (cf. 5.5) corresponde
\`a l'extension
$M_0$ de
$M_{-1}$ par $J$.

Comme au paragraphe {\it a)} on consid\`ere  une r\'esolution libre gradu\'ee de $C$ :
$$ \cdots \fl P_3 \Fl {\de_2} P_2 \Fl {\de_1} P_1 \Fl {\de_0} P_{0} \fl C \fl 0 \leqno
{(*)}$$
 et  l'\'el\'ement $u$ de  $\Ext_{R_A}^2 (C,H)$ correspond  \`a un
morphisme
$\wh u : P_2
\fl H$ v\'erifiant $\wh u \de_2 =0$ (modulo ceux qui proviennent de $P_1$).

\expl{5.20}  On reprend l'exemple 1.35.c. 
On a $C= k$, $H= k(-2)$ et la r\'esolution de $C$ est donn\'ee par le  complexe de Koszul.
On a
$P_2= R_A(-1)^4
\oplus R_A(-2)^6$, donc $\wh u : P_2 \fl H$ correspond \`a un \'el\'ement de $k^6$. Comme les
fl\`eches de cette r\'esolution sont
\`a coefficients dans l'id\'eal $(a,X,Y,Z,T)= (a,U)$ la condition $\wh u\de_2=0$ est
automatique. Il y a deux cas oppos\'es (qui correspondent \`a ceux envisag\'es en
1.35.c) :
$\wh u=0$ et $\wh u \neq 0$, donc surjectif.

\vskip 0.3 cm

1) Dans le cas $\wh u=0$ on  trouve la triade suivante (en notation chiffr\'ee) :
$$1^4, 2^6 \oplus 2 , 3^4 \Fl {d_1} \ul 0, 1^4 \oplus 2 \Fl {d_0} \ul 0, \quad {\rm avec}\quad d_0
= (a, U, 0) \quad {\rm et}$$
$$d_1 = \pmatrix {U&0&0&0 \cr -aI_4&V&0&0 \cr 0 &0 &a&U \cr}$$
o\`u $U,V$ sont les matrices usuelles de Koszul.
C'est la triade majeure associ\'ee
\`a la deuxi\`eme triade de 1.35.c. Le noyau de cette triade est de rang
$10$.

\vskip 0.3 cm

La famille minimale de courbes associ\'ee a une r\'esolution  :
$$ 2^2,3^4,4^3 \fl [1^4, 2^7,3^4 \fl \ul 0,1^4,2 \fl \ul 0] $$
o\`u la  notation  correspond \`a une r\'esolution 
$0 \fl P \fl N \fl I_{\cC}(h) \fl 0$ avec 
$P = 2^2,3^4,4^3 $ et 
pour $N$  le noyau d'une triade $1^4, 2^7,3^4 \fl \ul 0,1^4,2 \fl \ul 0$. Le d\'ecalage
$h$ est ici \'egal \`a $4$ et la r\'esolution est \'equivalente, du point de vue de la
caract\'eristique d'Euler, \`a celle obtenue en simplifiant tous les chiffres possibles,
c'est \`a dire (apr\`es d\'ecalage) 
$8^3 \fl 6^4$. On obtient une famille de courbes de degr\'e $24$ et genre $65$.

\vskip 0.3 cm

2) Dans le cas surjectif  on prend $L_0= P_1$ et $d_0= \de_0 = (a,U)$, on
appelle
$e_1,
\cdots,
e_4; \e_1,\cdots, \e_6$ la base de $P_2$ et on choisit $\wh u : P_2 \fl H$ qui envoie $\e_1$
sur
$1$ et les autres sur $0$. On a alors la suite exacte 
$$0 \fl \Im \de_2 \fl \Ker \wh u \fl \Ker \ov u \fl 0$$
et il est clair que $\Ker \wh u$ est engendr\'e par $e_1, \cdots, e_4, \e_2, \cdots,
\e_6, X\e_1, Y\e_1,Z\e_1,T\e_1$. Pour $\Ker \ov u$ les deux derniers vecteurs sont inutiles
car ils sont dans l'image de $\de_2$ modulo $\e_2, \cdots ,\e_6$. 
On peut donc prendre $L_1= 1^4,2^5,3^2$ avec la fl\`eche $d_1$ obtenue
en composant l'injection canonique de $\Ker u$ dans $P_2$ avec
$$\de_1 = \pmatrix {U&0\cr -aI_4&V \cr}$$
c'est-\`a-dire en rempla\c cant  dans  $\de_1$ la colonne de
degr\'e $2$ $^t(0,Y,-X,0,0)$ par les deux colonnes de degr\'e $3$ : $^t(0,XY,-X^2,0,0)$ et
$^t(0,Y^2,-XY,0,0)$. Cette triade est la triade majeure associ\'ee
\`a la premi\`ere triade de 1.35.c. Le noyau de $d_1$ est cette fois-ci de rang $7$.

 La famille minimale de courbes associ\'ee a pour chiffres 
$$2^2, 3^2, 4^2 \fl [1^4, 2^5, 3^2 \fl \ul 0, 1^4 \fl \ul 0].$$
Le d\'ecalage est $2$ et en simplifiant et d\'ecalant on obtient une r\'esolution du
type
$6^2
\fl 4^3$ soit une courbe de degr\'e
$12$ et genre $17$.

\vskip 0.3 cm

Deux remarques s'imposent \`a propos de cet exemple :

a) Le cas de l'extension non triviale est meilleur que l'autre
(au sens o\`u il donne un plus petit d\'ecalage).

b) La famille $12,17$ ainsi obtenue est la meilleure possible car la courbe minimale pour
le module $k \oplus k(-2)$ est une $12,17$ , cf. [MDP1] IV 6.10.

\tarte {i) L'exemple des courbes de degr\'e $4$ et genre $0$ : construction de la triade}

Nous reprenons ici l'exemple qui nous a servi de fil conducteur tout au long de ce
paragraphe. Si on part du module $M_0= R/(X,Y,Z,T^3)$ et de son sous-quotient $k(-1)$ on
a vu en 5.13 qu'une triade reliant ces modules a pour invariants $C = R_A/(a^p,X,Y,Z,T)$
et
$H = R_A(-1)/(X,Y,Z, a^nT, T^2)$. On prendra d\'esormais $p=n=1$. On a trouv\'e en 5.17.3
une triade admettant ces \'el\'ements, mais la courbe minimale \'etait de degr\'e $12$ et
genre $16$. On cherche donc \`a construire une autre triade dont la famille minimale soit de plus
petit degr\'e.

 La r\'esolution de $C$ est
toujours donn\'ee par le complexe de Koszul, comme en 5.20. On a donc, en particulier
$P_2 = 1^4,2^6$ et on s'int\'eresse aux $\wh u : P_2 \fl H$.

 On note tout d'abord que la contrainte de la proposition
5.14 impose \`a $\wh u$ d'\^etre surjectif. En effet, l'\'el\'ement $\theta
(u)$ de
$\Ext^1_{R_A} (M_{-1}, J)$ doit correspondre \`a $M_0$ et on v\'erifie que cela
implique que l'homomorphisme  $v : P_2 \T_A k \fl H\T_A k$ induit par $\wh u$ est surjectif,
donc aussi
$\wh u$ par Nakayama.

On reprend les notations de 5.20 pour la base de $P_2$, $e_1,
\cdots,
e_4; \e_1,\cdots, \e_6$. Il est facile de d\'eterminer les $\wh u : P_2 \fl H$ surjectifs. On
pose
$\wh u(e_i) = a_i$ et $\wh u(\e_i) = b_i$. En
\'ecrivant $\wh u\de_2=0$ et en travaillant dans $H=k[a,T]/(aT,T^2)$ on montre  que, si
 $\wh u$ est surjectif, on doit avoir $a_4 \in k^*$. \`A changement de base
pr\`es on peut supposer les autres $a_i$ nuls et il y a alors deux cas que nous \'etudions
ci-dessous : 

1)  le cas o\`u tous les
$b_i$ sont nuls,

2) le cas o\`u l'un des $b_i$ est non nul et  on peut alors supposer que c'est $b_1$.

\vskip 0.3 cm

1) Le premier exemple  va redonner la triade triviale, donc, cf. 5.18.3, une courbe
$(12,16)$. On d\'efinit   $\wh u$ en envoyant $e_4$ sur $\ov 1$ dans $H$  et les autres sur
$0$. On a bien $\wh u\de_2=0$ (car la ligne de
$\de_2$ qui correspond \`a $e_4$ est $(0,0,-X,0,-Y,-Z, 0,0,0,0)$ qui est nulle dans $H$).
Il est facile de calculer $\Ker  \wh u$, qui est engendr\'e par $e_1,e_2,e_3,
\e_1,\cdots,\e_6, Xe_4,Ye_4,Ze_4$, $aTe_4,T^2e_4$. Dans
$\Ker \ov u$ (i.e. modulo $\Im \de_2$) les vecteurs $Xe_4,Ye_4,Ze_4$ sont inutiles (car
ils sont dans $\Im \de_2$ modulo $e_1,e_2,e_3$). On prend donc 
$L_1=R_1= 1^3, 2^6, 2,3$ et on obtient $d_1$ en composant l'injection $\Ker \wh u \subset 
P_2$
avec $\de_1$ ce qui donne exactement la matrice obtenue en 5.17. On
note que $\Ker \ov u$ est de rang $11$ et le noyau $N$  de rang $7$.

\vskip 0.3 cm
2) Le deuxi\`eme exemple  va donner une triade dont la courbe
minimale sera la $(4,0)$ tant convoit\'ee.

 On prend ici $\wh u(e_4) =
\ov 1$, $\wh u(\e_1) = -\ov T$ et les autres nuls. On v\'erifie $\wh u\de_2=0$. Pour le $\e_1$ la
ligne correspondante de $\de_2$ est $(-a, 0,\cdots,0,T,Z)$ et $aT, T^2,ZT$ sont bien nuls
dans $H$. On montre que le noyau de $\ov u$ est engendr\'e par les images des $e_i$ ($i=
1,2,3$), des $\e_j$ ($j= 2,...,6$) et par $Te_4+\e_1$. On en d\'eduit $R_1= L_1=
1^3,2^6$, d'o\`u la triade 
$$\Lpont' = (1^3, 2^6 \Fl {d'_1} \ul 0, 1^4 \Fl {d'_0} \ul 0) \quad {\rm avec } \quad
d'_0=d_0 = (a,X,Y,Z,T) \quad {\rm et}
$$ 
$$ d'_1 = \pmatrix { X& Y & Z& 0&0&0&0&0&T^2 \cr
-a&0&0&Z&T&0&0&0&Y \cr
0&-a&0 &0&0&Z&T&0&-X\cr
0&0&-a&-X&0&-Y&0&T&0\cr
0&0&0&0&-X&0&-Y&-Z&-aT \cr
}$$ On voit que
$\Ker
\ov u$ est de rang $9$ et $N$ de rang $5$.

\tarte {j) L'exemple des $(4,0)$ : construction de la famille de courbes}

On calcule (en utilisant Macaulay) une r\'esolution minimale du noyau $N$ de la triade
$\Lpont'$ construite en {\it i)}. On obtient ainsi la matrice   $s: 2^2, 3^6,4^2
\fl 1^3,2^6$ avec 
$$s = \pmatrix{
\scriptstyle -Z&\scriptstyle 0&\scriptstyle 0&\scriptstyle 0&\scriptstyle 0&\scriptstyle
-YT&\scriptstyle -Y^2&\scriptstyle aT^2-XY&\scriptstyle 0&\scriptstyle T^3\cr
\scriptstyle  0&\scriptstyle -Z&\scriptstyle 0&\scriptstyle 0&\scriptstyle
0&\scriptstyle -XT&\scriptstyle aT^2+XY&\scriptstyle X^2&\scriptstyle T^3&\scriptstyle
0\cr \scriptstyle  X&\scriptstyle Y&\scriptstyle 0&\scriptstyle 0&\scriptstyle
T^2&\scriptstyle 0&\scriptstyle 0&\scriptstyle 0&\scriptstyle 0&\scriptstyle 0\cr
\scriptstyle  a&\scriptstyle 0&\scriptstyle T&\scriptstyle 0&\scriptstyle
-Y&\scriptstyle 0&\scriptstyle 0&\scriptstyle 0&\scriptstyle 0&\scriptstyle 0\cr
\scriptstyle  0&\scriptstyle 0&\scriptstyle -Z&\scriptstyle 0&\scriptstyle
0&\scriptstyle aY&\scriptstyle 0&\scriptstyle -a^2T&\scriptstyle -Y^2&\scriptstyle
-aT^2-XY \cr \scriptstyle  0&\scriptstyle a&\scriptstyle 0&\scriptstyle T&\scriptstyle
X&\scriptstyle 0&\scriptstyle 0&\scriptstyle 0&\scriptstyle 0&\scriptstyle 0\cr
\scriptstyle  0&\scriptstyle 0&\scriptstyle 0&\scriptstyle -Z&\scriptstyle
0&\scriptstyle -aX&\scriptstyle -a^2T&\scriptstyle 0&\scriptstyle -aT^2+XY&\scriptstyle
X^2 \cr \scriptstyle  0&\scriptstyle 0&\scriptstyle X&\scriptstyle Y&\scriptstyle
-aT&\scriptstyle 0&\scriptstyle 0&\scriptstyle 0&\scriptstyle 0&\scriptstyle 0 \cr
\scriptstyle  0&\scriptstyle 0&\scriptstyle 0&\scriptstyle 0&\scriptstyle
-Z&\scriptstyle 0&\scriptstyle -aY&\scriptstyle -aX&\scriptstyle -YT&\scriptstyle -XT
\cr \scriptstyle  }$$

On trouve $\a_2=2, \beta_2 =1$ et $b_0=1$, puis $\a_3 =5, \beta_3 =
4$, donc la fonction
$q$ est donn\'ee par $q(2) = 1$ et
$q(3)=3$ et on obtient
 une famille de courbes 
 $\cC$ avec la
r\'esolution :
$$ 0 \fl \cO_{\bP_A^3} (-2) \oplus \cO_{\bP_A^3} (-3)^3 \fl \cN \fl \cJ_{\cC}  \fl 0$$
$${\rm avec} \quad 0 \fl \cN \fl \cO_{\bP_A^3}(-1)^3 \oplus \cO_{\bP_A^3}(-2)^6 \fl
\cO_{\bP_A^3}  \oplus \cO_{\bP_A^3}(-1)^4 \fl \cO_{\bP_A^3} \fl 0.$$
Il s'agit d'une famille de courbes de degr\'e $4$ et genre $0$, dont la courbe
g\'en\'erique a la cohomologie d'une courbe de bidegr\'e $(1,3)$ sur une quadrique et la
courbe sp\'eciale celle de la r\'eunion disjointe d'une cubique plane et d'une droite. Une
cons\'equence de l'existence de cette famille est le th\'eor\`eme suivant (il s'agit du
troisi\`eme exemple de l'introduction) :

\th {Th\'eor\`eme 5.21}. Le sch\'ema de Hilbert $H_{4,0}$ des courbes localement
Cohen-Macaulay de degr\'e $4$ et genre $0$ de $\bP^3$ est connexe.

\dem On sait, cf. [MDP4] ou [El], que le sch\'ema de Hilbert est r\'eunion de deux
sous-sch\'emas irr\'eductibles $H_1$ (ouvert) et $H_2$ (ferm\'e), tous deux de dimension
$16$ avec les cohomologies des deux types
\'evoqu\'es ci-dessus. L'existence de la famille de courbes construite ci-dessus atteste
qu'il y a un point $t$ de $H_2$ dans l'adh\'erence de $H_1$. Alors, $H_{4,0}$ est l'union
des deux connexes $\ov {H_1}$ et $H_2$  qui ont une intersection non vide. Il est donc
connexe. Bien entendu, comme les dimensions  sont \'egales $H_2$ est
seulement sous-adh\'erent \`a $H_1$ (i.e. $\ov {H_1} \cap H_2 \neq \vide$ et non $\ov {H_1}
\subset H_2$ ). Les courbes de $H_2$ r\'eunions disjointes d'une cubique plane et d'une
droite ne sont pas dans l'adh\'erence de $H_1$ (car ce sont des points lisses du sch\'ema de
Hilbert). Les courbes de $H_2$ qui sont dans l'adh\'erence de
$H_1$ sont non r\'eduites, cf. [MDP4] 0.6 et 5.22 ci-dessous.

\rema {5.22} On peut expliciter une famille du type ci-dessus : on prend la r\'eunion des
courbes d'id\'eaux 
$$J_a = (X^2,XY,Y^2, X-a Y) \qquad {\rm et} \qquad J= (X^2,XZ,Z^2,XY-ZT).$$
 Elle est d\'efinie par l'id\'eal
$I_a$ de $R_A$ :
$$I_a = (X^2, XYZ, Y^2Z^2, XY^3-Y^2ZT, XZ^2-a YZ^2, -XZT-a XY^2+a YZT)$$
et on v\'erifie qu'il s'agit bien d'une famille plate de courbes de degr\'e $4$ et genre $0$
de la forme annonc\'ee.

\rema {5.23} Les deux triades consid\'er\'ees ci-dessus (cf. f)) qui joignent toutes deux
les modules $k(-1)$ et $R/(X,Y,Z,T^3)$ ne sont pas pseudo-isomorphes car les courbes
mi\-ni\-males associ\'ees ne sont pas les m\^emes (cf. 3.11). En revanche on peut
montrer en utilisant les m\'ethodes \'evoqu\'ees en 1.36 que les foncteurs associ\'es
sont les m\^emes ce qui fournit un nouvel exemple du type de 1.35.c.

\rema {5.24}  Dans le cas g\'en\'eral on ignore si le sch\'ema de Hilbert $H_{d,g}$ est connexe.
C'est vrai pour $d=2$ (il est irr\'eductible) et Nollet l'a montr\'e pour $d=3$, cf. [N]. La
notion de triade peut \^etre  une voie d'acc\`es \`a ce probl\`eme. Ainsi,  S.
A\"
\i t-Amrane a montr\'e par cette m\'ethode que le sch\'ema de Hilbert des courbes de degr\'e $d$
et de genre
$(d-3)(d-4) /2$ est connexe, g\'en\'eralisant le cas de $H_{4,0}$, cf. [AA].

\titre {R\'ef\'erences bibliographiques}

[AA] A\"\i t-Amrane S., Sur le sch\'ema de Hilbert $H_{d, (d-3)(d-4)/2}$, en
pr\'eparation.

[AG] Hartshorne R., Algebraic geometry, Graduate texts in Mathematics 52, Springer
Verlag, 1977.

[BB] Ballico E. et Bolondi G., The variety of module structures, Arch. der Math.
54, 1990, 397-408. 

[BBM] Ballico E., Bolondi G. et Migliore J., The Lazarsfeld-Rao problem for liaison
classes of two-codimensional subschemes of $\bP^n$, Amer. J. of Math., 113, 117-128,
1991.

[E] Eisenbud D., Commutative algebra, Graduate texts in Mathematics, Springer, 1995.

[El] Ellia Ph., On the cohomology of projective space curves, Bolletino U.M.I. 7,
9-A,  593-607, 1995.

[G] Ginouillac S.,  Sur les sch\'emas des modules de Rao de longueur $3$, note
CRAS, t. 320, S\'erie I,  1327-1330, 1995.

[H] Hartshorne R., Coherent functors, \`a para\^\i tre, Advances in Math.

[HMDP1] Hartshorne R.,  Martin-Deschamps M. et  Perrin D., Un th\'eor\`eme de Rao pour
les familles de courbes gauches, rapport de recherche LMENS-97-15, 1997.

[HMDP2] Hartshorne R., Martin-Deschamps M. et Perrin D., Construction de familles
minimales de courbes gauches,  rapport de recherche LMENS-97-29, 1997.

[Ho] Horrocks G., Vector bundles on the punctured spectrum of a local ring, Proc. London
Math. Soc. 14,  689-713, 1964.

[LR] Lazarsfeld R. et Rao A. P., Linkage of general curves of large degree, Lecture
notes 997, Springer Verlag, 1983, 267-289.

[M] Matsumura H., Commutative Ring Theory, Cambridge University Press 8,
1989.

[MDP 1]  Martin-Deschamps M. et  Perrin D., Sur la classification des
courbes gauches, Ast\'erisque, Vol. 184-185, 1990.

[MDP 2] Martin-Deschamps M. et  Perrin D., Courbes gauches et Modules de Rao,
J. reine angew. Math. 439 (1993), 103-145.

[MDP3]  Martin-Deschamps M. et  Perrin D., Construction de courbes lisses :
un th\'eor\`eme \`a la Bertini, rapport de recherche du LMENS 94-14, 1994.

[MDP4] Martin-Deschamps M. et  Perrin D., Le sch\'ema de Hilbert des courbes
localement de Cohen-Macaulay n'est (presque) jamais r\'eduit, Ann. scient. \'Ec. Norm.
Sup., $4^e$ s\'erie, t. 29, 757-785, 1996.

[N] Nollet S., The Hilbert scheme of degree three curves,  Ann. scient. \'Ec. Norm.
Sup., $4^e$ s\'erie, t. 30, 367-384, 1997.

[R] Rao A.P., Liaison among curves in $\bP^3$, Invent. Math., 50, 1979,
205-217.

[RD] Hartshorne R., Residues and duality, Lecture Notes in Math. 20, Springer Verlag,
1966.

[S] Sernesi E., Un esempio di curva ostruita in $\bP^3$, Seminario di variabili complesse,
Universita di Bologna (1981).

[V] Verdier J.-L., Cat\'egories d\'eriv\'ees, \'etat 0, {\it in} SGA $4 {1 \over 2}$,
Lecture Notes in Math. 569, Springer Verlag, 1977.

\bye